\documentclass{amsart}
\usepackage{mathabx}
\usepackage{amssymb}
\usepackage{mathrsfs}
\usepackage{stmaryrd}
\usepackage{bbm}
\usepackage{oldgerm}
\usepackage[english]{babel}
\usepackage[T1]{fontenc}
\usepackage[latin1]{inputenc}
\usepackage[all]{xy}

\newtheorem{teo}[subsection]{Theorem}
\newtheorem{prop}[subsection]{Proposition}
\newtheorem{cor}[subsection]{Corollary}
\newtheorem{lem}[subsection]{Lemma}
\newtheorem{conj}[subsection]{Conjecture}

\theoremstyle{definition}

\newtheorem{defi}[subsection]{Definition}
\newtheorem{rema}[subsection]{Remark}
\newtheorem{remas}[subsection]{Remarks}

\numberwithin{equation}{subsection}

\newcommand{\mQ}{{\mathbb Q}}
\newcommand{\mN}{{\mathbb N}}
\newcommand{\mZ}{{\mathbb Z}}
\newcommand{\mF}{{\mathbb F}}
\newcommand{\mA}{{\mathbb A}}

\newcommand{\mX}{{\mathbb X}}
\newcommand{\mY}{{\mathbb Y}}

\newcommand{\bD}{{\bf D}}

\newcommand{\bV}{{\bf V}}

\newcommand{\rH}{{\rm H}}
\newcommand{\rT}{{\rm T}}

\newcommand{\Spec}{{\rm Spec}}

\newcommand{\Sym}{{\rm Sym}}

\newcommand{\ord}{{\rm ord}}
\newcommand{\gp}{{\rm gp}}
\newcommand{\sat}{{\rm sat}}

\newcommand{\idemp}{{\rm idemp}}
\newcommand{\ab}{{\rm ab}}
\newcommand{\Gr}{{\rm Gr}}

\newcommand{\id}{{\rm id}}

\newcommand{\rR}{{\rm R}}
\newcommand{\pr}{{\rm pr}}

\newcommand{\res}{{\rm res}}
\newcommand{\alg}{{\rm alg}}
\newcommand{\rk}{{\rm rk}}
\newcommand{\fs}{{\rm fs}}
\newcommand{\ob}{{\rm Ob}}
\newcommand{\Sch}{{\rm Sch}}

\newcommand{\rsw}{{\rm rsw}}
\newcommand{\Hom}{{\rm Hom}}
\newcommand{\End}{{\rm End}}
\newcommand{\Aut}{{\rm Aut}}
\newcommand{\ZR}{{\rm ZR}}

\newcommand{\Mon}{{\bf Mon}}
\newcommand{\LS}{{\bf LS}}
\newcommand{\MLS}{{\bf MLS}}
\newcommand{\SNCP}{{\bf SNCP}}
\newcommand{\timesl}{\times^{\log}}
\newcommand{\chE}{{\Check{E}}}

\newcommand{\chf}{{\Check{f}}}

\newcommand{\chpi}{{\Check{\pi}}}

\newcommand{\chphi}{{\Check{\phi}}}
\newcommand{\chpsi}{{\Check{\psi}}}

\newcommand{\oF}{{\overline{F}}}
\newcommand{\oK}{{\overline{K}}}
\newcommand{\oL}{{\overline{L}}}
\newcommand{\oM}{\overline{M}}
\newcommand{\oN}{\overline{N}}

\newcommand{\oS}{{\overline{S}}}
\newcommand{\oT}{{\overline{T}}}

\newcommand{\oa}{{\overline{a}}}
\newcommand{\ok}{{\overline{k}}}
\newcommand{\of}{{\overline{f}}}
\newcommand{\os}{{\overline{s}}}
\newcommand{\ot}{{\overline{t}}}
\newcommand{\ou}{{\overline{u}}}

\newcommand{\ox}{{\overline{x}}}
\newcommand{\oy}{{\overline{y}}}
\newcommand{\oz}{{\overline{z}}}
\newcommand{\ocM}{{\bar{\cM}}}
\newcommand{\oeta}{{\overline{\eta}}}

\newcommand{\oxi}{{\overline{\xi}}}
\newcommand{\oLambda}{{\overline{\Lambda}}}

\newcommand{\uZ}{{\underline{Z}}}

\newcommand{\uGamma}{{\underline{\Gamma}}}

\newcommand{\cA}{{\mathscr A}}
\newcommand{\cB}{{\mathscr B}}
\newcommand{\cC}{{\mathscr C}}
\newcommand{\cD}{{\mathscr D}}
\newcommand{\cE}{{\mathscr E}}
\newcommand{\cF}{{\mathscr F}}
\newcommand{\cG}{{\mathscr G}}
\newcommand{\cI}{{\mathscr I}}
\newcommand{\cJ}{{\mathscr J}}
\newcommand{\cK}{{\mathscr K}}
\newcommand{\cL}{{\mathscr L}}
\newcommand{\co}{{\mathscr O}}

\newcommand{\cH}{{\mathscr H}}
\newcommand{\cM}{{\mathscr M}}
\newcommand{\cN}{{\mathscr N}}
\newcommand{\cP}{{\mathscr P}}

\newcommand{\cHom}{{\mathscr Hom}}
\newcommand{\cEnd}{{\mathscr End}}
\newcommand{\cDiv}{{\mathscr Div}}

\newcommand{\fF}{{\mathfrak F}}
\newcommand{\fS}{{\mathfrak S}}

\newcommand{\fm}{{\mathfrak m}}

\newcommand{\fp}{{\mathfrak p}}
\newcommand{\fq}{{\mathfrak q}}

\newcommand{\hLS}{\widehat{\LS}}

\newcommand{\tB}{{\widetilde{B}}}
\newcommand{\tC}{{\widetilde{C}}}

\newcommand{\tV}{{\widetilde{V}}}
\newcommand{\tM}{{\widetilde{M}}}

\newcommand{\tW}{{\widetilde{W}}}

\newcommand{\tZ}{{\widetilde{Z}}}
\newcommand{\tf}{{\widetilde{f}}}

\begin{document}

\title{Ramification and cleanliness}
\author{Ahmed Abbes and Takeshi Saito}
\address{A.A. CNRS \& Institut des Hautes \'Etudes Scientifiques, 35 route de Chartres, 91440 Bures-sur-Yvette, France}
\email{abbes@ihes.fr}
\address{T.S. Department of Mathematical Sciences, 
University of Tokyo, Tokyo 153-8914, Japan}
\email{t-saito@ms.u-tokyo.ac.jp}

\maketitle

\subsection*{Abstract} This article is devoted to studying the ramification of Galois torsors and 
of $\ell$-adic sheaves in characteristic $p>0$ (with $\ell\not=p$).  
Let $k$ be a perfect field of characteristic $p>0$, $X$ a smooth, separated and quasi-compact $k$-scheme, 
$D$ a simple normal crossing divisor on $X$, $U=X-D$, $\Lambda$ a finite local $\mZ_\ell$-algebra and 
$\cF$ a locally constant constructible sheaf of $\Lambda$-modules on $U$.
We introduce a {\em boundedness} condition on the ramification of $\cF$ along $D$, 
and study its main properties, in particular, some specialization properties that lead 
to the fundamental notion of {\em cleanliness} and to the definition of the {\em characteristic cycle} of $\cF$. 
The cleanliness condition extends the one introduced by Kato for rank $1$ sheaves. Roughly speaking, 
it means that the ramification of $\cF$ along $D$ is controlled by its ramification at the generic points of $D$. 
Under this condition, we propose a conjectural Riemann-Roch type formula for $\cF$. Some cases of this formula
have been previously proved by Kato and by the second author (T.S.).

\section{Introduction}

\subsection{}\label{intro1} 
The purpose of this article is to study the ramification of Galois torsors and
of $\ell$-adic sheaves in characteristic $p>0$ (with $\ell\not=p$),
developing the project started in \cite{as1,as2,as3,as4,saito1}. More precisely, this work is a sequel
to \cite{saito1}, though it can be read independently. The {\em leitmotiv} of our approach, in particular in 
this work, is to eliminate the ramification by blow-ups.  

\subsection{}\label{intro2}
Let $k$ be a perfect field of characteristic $p>0$, $X$ a smooth, separated and quasi-compact $k$-scheme, 
$D$ a simple normal crossing divisor on $X$ and $U=X-D$; we say for short that $(X,D)$ is an {\em snc-pair} over $k$.  
We fix a prime number $\ell$ different from $p$ and a finite local $\mZ_\ell$-algebra $\Lambda$.  
Let $\cF$ be a locally constant constructible sheaf of $\Lambda$-modules on $U$. The main problems 
in ramification theory are the following~:
\begin{itemize}
\item[(A)] to describe the ramification of $\cF$ along $D$;
\item[(B)] to give a Riemann-Roch type formula for $\cF$, that is, to compute the Euler-Poincaré characteristic
with compact support of $\cF$ on $U$ in terms of its  invariants of ramification (provided by (A)).
\end{itemize} 
In \cite{as3}, we gave cohomological answers to both problems that rely on the notion of {\em characteristic class}
of $\cF$. In this article, we develop a more geometric approach to problem (A) and give a conjectural formula for (B),
based on the finer notion of {\em characteristic cycle} of $\cF$. For this purpose, we start by studying the ramification of 
Galois torsors over $U$, that is, torsors over $U$ for the étale topology, under finite constant groups. 

\subsection{}\label{intro3} 
Our approach is based on a geometric construction introduced in \cite{as3,as4,saito1}.
Let $D_1,\dots,D_m$ be the irreducible components of $D$ and 
let $(X\times_kX)'$ be the blow-up of $X\times_kX$ along $D_i\times_kD_i$ for all $1\leq i\leq m$. 
We define the {\em framed self-product} $X\Asterisk_kX$ of $(X,D)$ over $k$
as the open subscheme of $(X\times_kX)'$ obtained by removing the 
strict transforms of $D\times_kX$ and $X\times_kD$ (called the logarithmic self-product in \cite{saito1}). 
We give in \ref{fram16} an equivalent definition using logarithmic geometry, that extends to more general situations. 
The diagonal morphism $\delta_X\colon X\rightarrow X\times_kX$ lifts uniquely to a morphism
$\delta\colon X\rightarrow X\Asterisk_kX$, called the {\em framed diagonal} of $(X,D)$ 
(and the logarithmic diagonal in \cite{saito1}).   
We consider $X\Asterisk_kX$ as an $X$-scheme by the second projection. 

Let $R$ be an effective rational divisor on $X$ with support in $D$ 
(i.e., a sum of non-negative rational multiples of the irreducible components of $D$).
We define in \ref{log4} the {\em dilatation $(X\Asterisk_kX)^{(R)}$  of $X\Asterisk_kX$ along $\delta$ of thickening $R$}.
It is an affine scheme over $X\Asterisk_kX$ that fits into a canonical Cartesian diagram 
\begin{equation}\label{intro3a}
\xymatrix{
U\ar[r]^-(0.5){\delta_U}\ar[d]_j&{U\times_kU}\ar[d]^-(0.5){j^{(R)}}\\
X\ar[r]^-(0.5){\delta^{(R)}}&{(X\Asterisk_kX)^{(R)}}}
\end{equation}
where $j^{(R)}$ is a canonical open immersion, $\delta^{(R)}$ is the unique morphism lifting $\delta$,
$j$ is the canonical injection and $\delta_U$ is the diagonal morphism. 
If $R$ has integral coefficients, then $(X\Asterisk_kX)^{(R)}$ is a dilatation 
in the sense of Raynaud. More precisely, $(X\Asterisk_kX)^{(R)}$ is the maximal open subscheme of the blow-up
of $X\Asterisk_kX$ along $\delta(R)$, where the exceptional divisor is equal to the pull-back of $R$ by 
the second projection to $X$  (cf. \ref{tub1}). 

\subsection{}\label{intro4} 
Let $V$ be a Galois torsor over $U$ of group $G$ and let 
$R$ be an effective rational divisor on $X$ with support in $D$.
We introduce a fundamental {\em boundedness} property of the ramification of $V/U$ along $D$.
We consider $V\times_kV$ as a Galois torsor over $U\times_kU$ of group $G\times G$, and  
denote by $W$ the quotient of $V\times_kV$ by $\Delta(G)$, where $\Delta\colon G\rightarrow G\times G$ 
is the diagonal homomorphism. 
The diagonal morphism $\delta_V\colon V\rightarrow V\times_kV$ induces a morphism
$\varepsilon_U\colon U\rightarrow W$ lifting the diagonal morphism $\delta_U\colon U\rightarrow U\times_kU$.
Note that $W$ represents the sheaf of isomorphisms of $G$-torsors
from $U\times_kV$ to $V\times_kU$ over $U\times_kU$, and that $\varepsilon_U$ corresponds to the 
identity isomorphism of $V$ (identified with the pull-backs of $U\times_kV$ and $V\times_kU$
by $\delta_U$).
We denote by $Z$ the integral closure of $(X\Asterisk_kX)^{(R)}$ in $W$,
by $\pi\colon Z\rightarrow (X\Asterisk_kX)^{(R)}$ the canonical morphism and by $\varepsilon\colon X\rightarrow Z$
the morphism induced by $\varepsilon_U\colon U\rightarrow W$. We have $\pi\circ \varepsilon=\delta^{(R)}$.
\begin{equation}\label{intro4a}
\xymatrix{
&W\ar[r]\ar[d]&Z\ar[d]^{\pi}&\\
U\ar@/^1pc/[ru]^{\varepsilon_U}\ar[r]^(0.4){\delta_U}&{U\times_kU}\ar[r]&{(X\Asterisk_kX)^{(R)}}&
X\ar@/_1pc/[ul]_{\varepsilon}\ar[l]_-(0.4){\delta^{(R)}}}
\end{equation}
Let $x\in X$. We say that the ramification of $V/U$ at $x$ is {\em bounded} by $R+$ if the morphism $\pi$ is étale at 
$\varepsilon(x)$, and that the ramification of $V/U$ along $D$ is {\em bounded} by $R+$ if $\pi$ is étale over an open neighborhood of $\varepsilon(X)$. We establish several properties of this notion. 
First, we prove that it satisfies descent for faithfully flat and log-smooth morphisms \eqref{ram-cov32}. 
The second property plays a key role in this article~: 
if $R$ has integral coefficients, we prove that the ramification of $V/U$ along $D$ 
is bounded by $R+$ if and only if 
there exists an open neighborhood $Z_0$ of $\varepsilon(X)$ in $Z$ which is étale over $(X\Asterisk_kX)^{(R)}$
and such that $\pi(Z_0)$ contains $(X\Asterisk_kX)^{(R)}\times_XR$ \eqref{ram-cov4}. 
Third, we relate this notion to its analogue for finite separable extensions 
of local fields (with possibly imperfect residue fields) defined in \cite{as1,as2}~: let 
$\xi$ be a generic point of $D$, $\oxi$ a geometric point of $X$ above $\xi$, $S$ the strict localization of $X$ at $\oxi$,
$K$ the fraction field of $\Gamma(S,\co_S)$ and $r$ the multiplicity of $R$ at $\xi$. 
We put $V\times_U\Spec(K)=\Spec(L)$, where $L=\prod_{i=1}^nL_i$ is a finite product of finite separable extensions of $K$. 
We prove in \ref{ram-cov8} that the ramification of $V/U$ at $\xi$ is bounded by $R+$ if and only if, for every $1\leq i\leq n$, the logarithmic ramification of $L_i/K$ is bounded by $r+$ in the sense of \cite{as1,as2}.

\subsection{}\label{intro6} 
Let $V$ be a Galois torsor over $U$ of group $G$, $Y$ the integral closure of $X$ in $V$,
and $R$ an effective rational divisor on $X$ with support in $D$. Assume that the following conditions are satisfied~:
\begin{itemize}
\item[(i)] for every geometric point $\oy$ of $Y$, 
the inertia group $I_\oy\subset G$ of $\oy$ has a normal $p$-Sylow subgroup~; 
\item[(ii)] for every generic point $\xi$ of $D$, the ramification of $V/U$ at $\xi$ is bounded by $R+$.
\end{itemize}
Then we prove that the ramification of $V/U$ along $D$ is bounded by $R+$ \eqref{ram-cov6}. This result is an analogue 
of the Zariski-Nagata purity theorem (\cite{sga2} X 3.4).

\subsection{}\label{intro7} 
Let $V$ be a Galois torsor over $U$ of group $G$. We define the {\em conductor} 
of $V/U$ relatively to $X$ to be the minimum effective rational divisor $R$ on $X$ with support in 
$D$ such that for every generic point $\xi$ of $D$, the ramification of $V/U$ at $\xi$ is bounded by $R+$. 
This terminology may be slightly misleading as the ramification of $V/U$ along $D$ may not be bounded by $R+$
in general. However, we prove in \ref{ram-cov18}, as a consequence of \ref{intro6}, 
that under a strong form of resolution of singularities, 
there exists an snc-pair $(X',D')$ over $k$ and a proper morphism $f\colon X'\rightarrow X$ 
inducing an isomorphism $X'-D'\stackrel{\sim}{\rightarrow} U$, such that  
if we denote by $R'$ the conductor of $V/U$ relatively to $X'$,
the ramification of $V/U$ along $D'$ is bounded by $R'+$.

\subsection{}\label{intro8}
Let $\cF$ be a locally constant constructible sheaf of $\Lambda$-modules on $U$,
$R$ an effective rational divisor on $X$ with support in $D$, $x\in X$ and $\ox$ a 
geometric point of $X$ above $x$. Recall that $\Lambda$ is a finite local $\mZ_\ell$-algebra \eqref{intro2}. 
We denote by $\pr_1,\pr_2\colon U\times_kU\rightarrow U$
the canonical projections and put 
\begin{equation}\label{intro8a}
\cH(\cF)=\cHom(\pr_2^*\cF,\pr_1^*\cF).
\end{equation} 
We prove in \ref{ram2} that the base change morphism
\begin{equation}\label{intro8b}
\alpha\colon \delta^{(R)*}j^{(R)}_*(\cH(\cF))\rightarrow j_*\delta_U^*(\cH(\cF))=j_*(\cEnd(\cF))
\end{equation}
relatively to the Cartesian diagram \eqref{intro3a} is injective. Furthermore, the following conditions are equivalent~:
\begin{itemize}
\item[(i)] The stalk $\alpha_\ox$ of the morphism $\alpha$ at $\ox$ is an isomorphism. 
\item[(ii)] There exists a Galois torsor $V$ over $U$ trivializing $\cF$ 
such the ramification of $V/U$ at $x$ is bounded by $R+$.
\end{itemize}
We give also other useful equivalent conditions. 
We say that the {\em ramification of $\cF$ at $\ox$ is bounded by $R+$} 
if $\cF$ satisfies these equivalent conditions. We say that the {\em ramification of $\cF$ along $D$ is bounded by $R+$} if
the ramification of $\cF$ at $\ox$ is bounded by $R+$ for every geometric point $\ox$ of $X$.
We establish several properties of this notion similar to those for Galois torsors. 
In particular, we relate it to the analogue notion for Galois representations 
of local fields (with possibly imperfect residue fields) \eqref{ram-loc2}. 

\subsection{}\label{intro9}
Let $\cF$ be a locally constant constructible sheaf of $\Lambda$-modules on $U$.
We define the {\em conductor of $\cF$ relatively to $X$} to be 
the minimum of the set of effective rational divisors $R$ on $X$ with support in $D$ 
such that for every geometric point $\oxi$ of $X$ above a generic point of $D$,
the ramification of $\cF$ at $\oxi$ is bounded by $R+$. 
As for Galois torsors, this terminology may be slightly misleading as the ramification of $\cF$ along $D$ 
may not be bounded by $R+$ in general. However, we prove that under a strong form of resolution of singularities, 
there exists an snc-pair $(X',D')$ over $k$ and a proper morphism $f\colon X'\rightarrow X$ 
inducing an isomorphism $X'-D'\stackrel{\sim}{\rightarrow} U$, such that  
if we denote by $R'$ the conductor of $\cF$ relatively to $X'$,
the ramification of $\cF$ along $D'$ is bounded by $R'+$ \eqref{ram10}.

\subsection{}\label{intro10}
The last part of this article is devoted to studying important specialization properties that lead 
to the fundamental notion of {\em cleanliness} and to the definition of the {\em characteristic cycle}. 
Let $R$ be an effective divisor on $X$ with support in $D$.\footnote{We consider rational divisors on
$X$ with support in $D$ and integral coefficients as Cartier divisors on $X$.} 
We prove \eqref{tub26} that $(X\Asterisk_kX)^{(R)}$ is smooth over $X$ and that 
\begin{equation}\label{intro10a}
E^{(R)}=(X\Asterisk_k X)^{(R)}\times_XR
\end{equation}
is canonically isomorphic to the twisted logarithmic tangent bundle 
\[
\bV(\Omega^1_{X/k}(\log D)\otimes_{\co_X}\co_{X}(R))\times_{X}R
\]
over $R$ (cf. \ref{not15} for the convention on vector bundles). We denote by $\chE^{(R)}$ the dual vector bundle.
Consider the following commutative diagram with Cartesian squares.
\begin{equation}\label{intro10b}
\xymatrix{
{E^{(R)}}\ar[r]\ar[d]&{(X\Asterisk_kX)^{(R)}}\ar[d]^{\pr_2}&{U\times_kU}\ar[d]\ar[l]_-(0.4){j^{(R)}}\\
{R}\ar[r]&X&U\ar[l]}
\end{equation}
Let $\cG$ be a sheaf of $\Lambda$-modules on $U\times_kU$. 
We call {\em $R$-specialization} of $\cG$ and denote by $\nu_R(\cG,X)$, 
the sheaf on $E^{(R)}$ defined by 
\begin{equation}\label{intro10c}
\nu_R(\cG, X)=j^{(R)}_*(\cG)|E^{(R)}.
\end{equation}

Let $\cF$ be a locally constant constructible sheaf of $\Lambda$-modules on $U$
such that its ramification along $D$ is bounded by $R+$ and let $\cH(\cF)$ be the sheaf on $U\times_kU$
defined in \eqref{intro8a}. We prove in \ref{ram12} that $\nu_R(\cH(\cF), X)$  
is {\em additive}, which means that its restrictions to the fibers of  $E^{(R)}$ over $R$ 
are invariant by translation (cf. \ref{add1}). 
This important property was first proved in (\cite{saito1}  2.25); we give a new proof in \ref{ram12}.

We fix a non-trivial additive character $\psi\colon \mF_p\rightarrow \Lambda^\times$
and denote by $S\subset \chE^{(R)}$ the
support of the Fourier-Deligne transform of $\nu_R(\cH(\cF), X)$ relatively to $\psi$ (cf. \ref{add3} and \ref{supp}). 
The additivity of $\nu_R(\cH(\cF), X)$ is equivalent to the fact 
that, for every $x\in R$, the set $S\cap \chE^{(R)}_x$ is finite \eqref{add4}. 
We call $S$ the {\em Fourier dual support} of $\nu_R(\cH(\cF), X)$. 
We prove in fact that $S$ is the underlying space of a closed subscheme of $\chE^{(R)}$
which is finite over $R$ \eqref{ram15}. Note that $S$ is a priori a constructible subset of $\chE^{(R)}$
and that it is not obvious that it is closed in $\chE^{(R)}$. We say that $\nu_R(\cH(\cF), X)$ is {\em non-degenerate} 
if $S$ does not meet the zero section of $\chE^{(R)}$ over $R$. 

\subsection{}\label{intro11}
We need in the following to recall a few facts from the ramification theory 
of local fields with imperfect residue fields developed in \cite{as1}, \cite{as2} and \cite{saito1}.
We refer to § \ref{rrt} for a more detailed review. 
Let $K$ be a discrete valuation field, $\co_K$ the valuation ring of $K$, $F$ the residue field of $\co_K$, 
$\oK$ a separable closure of $K$ and $\cG$ the Galois group of $\oK/K$. 
We assume that $\co_K$ is henselian and that $F$ has characteristic $p$. 
In (\cite{as1} 3.12), we defined a decreasing filtration $\cG_{\log}^r$ $(r\in \mQ_{\geq 0})$ of $\cG$
by closed normal subgroups, called the {\em logarithmic ramification filtration}. 
For a rational number $r\geq 0$, 
we put 
\begin{eqnarray*}
\cG_{\log}^{r+}&=&\overline{\bigcup_{s>r}\cG_{\log}^{s}},\\
\Gr_{\log}^r(\cG)&=&\cG_{\log}^r/\cG_{\log}^{r+}.
\end{eqnarray*}
This filtration satisfies the following properties, among others~: 
\begin{itemize}
\item[(i)] The group $\cP=\cG_{\log}^{0+}$ is the wild inertia subgroup of $\cG$, i.e., 
the $p$-Sylow subgroup of the inertia subgroup $\cG_{\log}^0$ (\cite{as1} 3.15). 
\item[(ii)] For every rational number $r>0$, the group $\Gr_{\log}^r(\cG)$ is abelian 
and is contained in the centre of the pro-$p$-group $\cP/\cG^{r+}_{\log}$ (\cite{as2} Theorem 1). 
\end{itemize}
Further properties are stated below. 

For any finite discrete $\Lambda$-representation $M$ of $\cG$, we have a canonical {\em slope decomposition} 
\begin{equation}\label{intro11a}
M=\oplus_{r\in \mQ_{\geq 0}} M^{(r)},
\end{equation}
characterised by the following properties (cf. \ref{rrt3})~: $M^{(0)}=M^{\cP}$ and  
for every $r>0$, 
\begin{equation}\label{intro11b}
(M^{(r)})^{\cG^r_{\log}}=0\ \ \ {\rm and}\ \ \ (M^{(r)})^{\cG^{r+}_{\log}}=M^{(r)}.
\end{equation}
The values $r\geq 0$ for which $M^{(r)}\not=0$ are called the {\em slopes} of $M$. 
We say that $M$ is {\em isoclinic} if it has only one slope. If $M$ is isoclinic of slope $r>0$, 
we have a canonical {\em central character decomposition} 
\begin{equation}\label{intro11c}
M=\oplus_{\chi } M_\chi,
\end{equation}
where the sum runs over finite characters 
$\chi\colon \Gr^r_{\log}\cG\rightarrow \Lambda_\chi^\times$ 
such that $\Lambda_\chi$ is a finite étale $\Lambda$-algebra (cf. \ref{rrt6}). 

We assume moreover that $K$ has characteristic $p$ and that $F$ is of finite type over $k$.
Let $\Omega^1_{\co_K}(\log)$ be the $\co_K$-module of logarithmic $1$-differential forms of $\co_K$ and
$\Omega^1_{F}(\log)=\Omega^1_{\co_K}(\log)\otimes_{\co_K}F$ (cf. \ref{rrt10}). We have a canonical exact
sequence $0\rightarrow \Omega^1_F\rightarrow \Omega^1_{F}(\log)\rightarrow F\rightarrow 0$.
We denote by $\co_{\oK}$ the integral closure of $\co_{K}$ in $\oK$, by $\oF$ the residue field of $\co_{\oK}$,
by $\ord$ the valuation of $\oK$ normalized by $\ord(K^\times)=\mZ$ and, for any rational number $r$, 
by $\fm_{\oK}^r$ (resp. $\fm_{\oK}^{r+}$) the $\co_{\oK}$-module of elements $x\in \oK$ such that 
$\ord(x)\geq r$ (resp. $\ord(x)>r$). 

The additivity property presented in \ref{intro10} is the geometric incarnation of an important property of the logarithmic 
ramification filtration proved in (\cite{saito1} 1.24), 
namely, for any rational number $r>0$, the group $\Gr^r_{\log}\cG$ 
is an $\mF_p$-vector space, and we have a canonical injective homomorphism 
\begin{equation}\label{intro11d}
\rsw\colon \Hom_{\mZ}(\Gr^r_{\log}\cG,\mF_p)\rightarrow 
\Hom_{\oF}(\fm_{\oK}^{r}/\fm_{\oK}^{r+}, \Omega^1_F(\log)\otimes_F\oF),
\end{equation}
called the {\em refined Swan conductor} (cf. \ref{rrt12}).

\subsection{}\label{intro12}
Let $\cF$ be a locally constant constructible sheaf of $\Lambda$-modules on $U$, 
$\xi$ a generic point of $D$, $X_{(\xi)}$ the henselization of $X$ at $\xi$, 
$\eta_\xi$ the generic point of $X_{(\xi)}$, 
$\oeta_\xi$ a geometric generic point of $X_{(\xi)}$ and $\cG_{\xi}$ the Galois group
of $\oeta_\xi$ over $\eta_\xi$.  
We say that $\cF$ is {\em isoclinic} at $\xi$ if the representation $\cF_{\oeta_\xi}$
of $\cG_{\xi}$ is isoclinic, and that $\cF$ is {\em isoclinic along $D$} if it is isoclinic at all generic points of $D$.  

Assume first that $\cF$ is isoclinic along $D$, and let $R$ be its conductor relatively to $X$. 
We say \eqref{clean7} that $\cF$ is {\em clean} along $D$ if the following conditions are satisfied~:
\begin{itemize}
\item[{\rm (i)}] the ramification of $\cF$ along $D$ is bounded by $R+$; 

\item[{\rm (ii)}] there exists a log-smooth morphism of snc-pairs $f\colon (X',D')\rightarrow (X,D)$ over $k$
such that the morphism $X'\rightarrow X$ is faithfully flat, that $R'=f^*(R)$ has integral coefficients, 
and if we put $U'=X'-D'$ and $\cF'=\cF|U'$, 
that the $R'$-specialization $\nu'_{R'}(\cH(\cF'), X')$ of $\cH(\cF')$ in the sense of \eqref{intro10c}
relatively to $(X',D')$, is additive and non-degenerate.
\end{itemize}
Note that we may replace (ii) by the stronger condition that it holds for any morphism $f$ satisfying the same
assumptions (cf. \ref{clean71}).

This notion can be extended to general sheaves as follows. 
Let $\ox$ be a geometric point of $X$.  We say that $\cF$ is {\em clean} at $\ox$ 
if there exists an étale neighborhood $X'$ of $\ox$ in $X$ such that, 
if we put $U'=U\times_XX'$ and denote by $D'$ the pull-back of $D$ over $X'$, 
there exists a finite decomposition 
\begin{equation}
\cF|U'=\oplus_{1\leq i\leq n} \cF'_i
\end{equation}
of $\cF|U'$ into a direct sum of locally constant constructible sheaves of $\Lambda$-modules 
$\cF'_i$ $(1\leq i\leq n)$ on $U'$ which are isoclinic and clean along $D'$ in the previous sense. 
We say that $\cF$ is {\em clean} along $D$ if it is clean at all geometric points of $X$ (cf. \ref{clean8}). 
Note that for isoclinic sheaves, the two definitions are equivalent \eqref{clean11}.

The notion of cleanliness was first introduced by Kato for rank $1$ sheaves in \cite{kato1}. Our definition extends his. 
It was extended to isoclinic sheaves by the second author (T. S.) in (\cite{saito1} § 3.2).

Roughly speaking, if $\cF$ is clean along $D$, then its ramification along $D$ is controlled by 
its ramification at the generic points of $D$. This is the main idea behind the following definition of the 
{\em characteristic cycle} of $\cF$. 

\subsection{}\label{cc}
We assume that $X$ is connected and denote by $d$ the dimension of $X$, 
by $\rT^*_X(\log D)=\bV(\Omega^1_{X/k}(\log D))$ 
the logarithmic cotangent bundle of $X$ and by $\xi_1,\dots,\xi_n$ the generic points of $D$.
For each $1\leq i\leq n$, we denote by $F_i$ the residue
field of $X$ at $\xi_i$, by $S_i=\Spec(\co_{K_i})$ the henselization of $X$ at $\xi_i$ and by $\eta_i=\Spec(K_i)$ 
the generic point of $S_i$. We fix a separable closure $\oK_i$ of $K_i$ and 
denote by $\cG_i$ the Galois group of $\oK_i/K_i$. 
 
Let $\cF$ be a locally constant constructible sheaf of free $\Lambda$-modules on $U$ which is {\em clean} along $D$.
We denote by $M_i$ the $\Lambda[\cG_i]$-module corresponding to $\cF|\eta_i$. Let 
\begin{equation}\label{cca}
M_i=\oplus_{r\in \mQ_{\geq 0}}M_i^{(r)}
\end{equation}
be its slope decomposition and, for each rational number $r>0$,
\begin{equation}\label{ccb}
M_i^{(r)}=\oplus_{\chi} M_{i,\chi}^{(r)}
\end{equation}
the central character decomposition of $M_i^{(r)}$.  
Note that $M_{i,\chi}^{(r)}$ is a free $\Lambda$-module of finite type for all $r>0$ and all $\chi$. 
By enlarging $\Lambda$, we may assume that for all rational numbers $r>0$ and all central 
characters $\chi$ of $M_i^{(r)}$ (i.e., all characters $\chi\colon \Gr_{\log}^r\cG_i\rightarrow \Lambda_\chi^\times$
that appear in the decomposition \eqref{ccb}), we have $\Lambda_\chi=\Lambda$.
Since $\Gr_{\log}^r\cG_i$ is abelian and killed by $p$ \eqref{rrt12}, $\chi$ factors uniquely as 
$\Gr_{\log}^r\cG_i\rightarrow \mF_p\stackrel{\psi}{\rightarrow} \Lambda^\times$, where $\psi$ is the non-trivial 
additive character fixed in \ref{intro10}. 
We denote  also by $\chi\colon \Gr_{\log}^r\cG_i\rightarrow \mF_p$ the induced character and by  
\begin{equation}
\rsw(\chi)\colon \fm^r_{\oK_i}/\fm^{r+}_{\oK_i}\rightarrow \Omega^1_{F_i}(\log)\otimes\oF_i
\end{equation}
its refined Swan conductor \eqref{intro11d} (where the notation are defined as in \ref{intro11} with $K=K_i$). 
Let $F_\chi$ be the field of definition of $\rsw(\chi)$, which is a finite
extension of $F_i$ contained in $\oF_i$.
The refined Swan conductor $\rsw(\chi)$ defines a line $L_\chi$ in $\rT^*_X(\log D)\otimes_X F_\chi$.
Let $\oL_\chi$ be the closure of the image of $L_\chi$ in $\rT^*_X(\log D)$. For each $1\leq i\leq n$, we put
\begin{equation}
CC_i(\cF)=\sum_{r\in \mQ_{>0}}\sum_{\chi}\frac{r\cdot \rk_{\Lambda}(M_{i,\chi}^{(r)})}{[F_\chi:F_i]}[\oL_\chi],
\end{equation}
which is a $d$-cycle on $\rT^*_X(\log D)\times_XD_i$. It follows from the proof of (\cite{saito1} 1.26) 
that the coefficient of $[\oL_\chi]$ is an element of $\mZ[\frac 1 p]$, and hence gives an element of $\Lambda$.  

Let $\sigma\colon X\rightarrow \rT^*_X(\log D)$ be the zero-section of $\rT^*_X(\log D)$ over $X$. 
We define the {\em characteristic cycle} of $\cF$ and denote by $CC(\cF)$, 
the $d$-cycle on $\rT^*_X(\log D)$ defined by
\begin{equation}
CC(\cF)=\rk_\Lambda(\cF)[\sigma]-\sum_{1\leq i\leq n}CC_i(\cF).
\end{equation}

Recall (\cite{as3} 2.1.1) that we associated to $j_!\cF$ a {\em characteristic class}, denoted by $C(j_!\cF)$,
which is a section of $\rH^0(X,\cK_X)$, where $\cK_X=f^!\Lambda$ and $f\colon X\rightarrow \Spec(k)$
is the structural morphism.

\begin{conj}\label{cc1}
Under the assumptions of \eqref{cc}, we have in $\rH^0(X,\cK_X)$
\begin{equation}
C(j_!\cF)=(CC(\cF),[\sigma]),
\end{equation} 
where the right hand side is the intersection pairing relatively to $\rT^*_X(\log D)$.
\end{conj}

Kato defined the characteristic cycle of a clean sheaf of rank $1$ in \cite{kato3}.  
The second author (T. S.) extended the definition to isoclinic and clean sheaves in (\cite{saito1} 3.6) 
and proved conjecture \ref{cc1} for these sheaves in (loc. cit. 3.7). 

\subsection{}
We may optimistically expect that for any locally constant constructible sheaf $\cF$ of $\Lambda$-modules on $U$, 
there exists an snc-pair $(X',D')$ over $k$ and a proper morphism of snc-pairs $(X',D')\rightarrow (X,D)$
inducing an isomorphism $X'-D'\stackrel{\sim}{\rightarrow} U$ such that  $\cF$ is clean along $D'$.  
Kato proved this property for rank $1$ sheaves on surfaces (\cite{kato3} 4.1).

\subsection{} 
We introduce in § \ref{not} the general notation and conventions for this article
and prove some preliminary results. Section \ref{add} is devoted to studying additive sheaves on vector bundles.
We recall in § \ref{tub} the classical notion of dilatation. The first part of  section \ref{fram} 
contains a detailed review of the notion of frame in logarithmic geometry and some representability 
results following \cite{ks1}. Its second part is devoted to the study of snc-pairs over $k$~;
we introduce the framed products and extend the notion of dilatation to rational divisors. 
Section \ref{rrt} is a review of ramification theory of local fields with imperfect residue fields. 
The last two sections, § \ref{ram-cov} and  § \ref{ram}, are the heart of this article. The former is devoted to 
studying the ramification of galois torsors and the latter to studying the ramification of $\ell$-adic sheaves.

\section*{Acknowledgment}

The first author would like to thank the Institut des Hautes \'Etudes Scientifiques 
and the University of Tokyo for their hospitalities. Both authors are grateful
to the Centre \'Emile Borel at the Institut Henri Poincaré for its hospitality during the Galois trimester. 
This research is partially supported by Grants-in-Aid for scientific research (B) 18340002 and
JSPS Core-to-Core Program 18005 {\em New Developments of Arithmetic Geometry,
Motive, Galois Theory, and Their Practical Applications}.

\section{Notation and preliminaries} \label{not}

\subsection{}\label{not1}
In this article, we fix a prime number $p$, 
a perfect field $k$ of characteristic $p$ and an algebraic closure $\ok$ of $k$. 
All $k$-schemes are assumed to be separated of finite type over $k$. 
We fix also a prime number $\ell$ different from $p$, 
a finite local $\mZ_\ell$-algebra $\Lambda$ 
and a non-trivial additive character $\psi\colon \mF_p\rightarrow \Lambda^\times$. 

\subsection{}\label{not15}
Let $X$ be a scheme and $\cE$ a locally free $\co_X$-module of finite type. We call 
the spectrum of the quasi-coherent $\co_X$-algebra $\Sym_{\co_X}(\cE)$ the vector bundle over 
$X$ defined by $\cE$ and denote it by $\bV(\cE)$. 

\subsection{}
Let $X$ be a locally noetherian scheme. In this article, a {\em Galois torsor} over $X$ of group  $G$
stands for a torsor over $X$ for the étale topology under a finite constant group $G$,
that is, a principal covering of $X$ of Galois group $G$ in the sense of (\cite{sga1}  V 2.8).

\subsection{}\label{NS15}
Let $X$ be a normal and locally noetherian scheme, $U$ a dense open subscheme of $X$, 
$V$ a Galois torsor over $U$ of group $G$, and $Y$ the integral closure of $X$ in $V$. 
Then $G$ acts on $Y$ and we have $X=Y/G$. 
Let $y$ be a point of $Y$ and $\oy$ a geometric point of $Y$ above $y$. 
Recall that the inertia group $I_y$ of $y$ is the subgroup of elements $\sigma\in G$ such that 
$g(y)=y$ and that $g$ acts trivially on $\kappa(y)$. 
It is convenient to denote $I_{y}$ also by $I_{\oy}$ and to call it also the inertia group of $\oy$. 
Assume that $X$ is {\em universally Japanese}, which means that every point 
of $X$ has an affine open neighborhood whose ring is  universally Japanese (\cite{ega4} 0.23.1.1). 
Let $\ox$ be the image of $\oy$ in $X$. Then we have a canonical isomorphism 
\begin{equation}\label{NS15a}
X_{(\ox)}\times_XV\simeq \bigsqcup_{\oz\in Y\otimes_X\kappa(\ox)}Y_{(\oz)}\times_YV,
\end{equation}
where $Y_{(\oz)}$ is the strict localization of $Y$ at $\oz$. 
Since $Y$ is normal, $Y_{(\oz)}\times_YV$ is integral. Therefore, $I_{\oy}$ is the stabilizer of 
$Y_{(\oy)}\times_YV$ in $G$.

\subsection{}
Recall that a scheme locally of finite type over a universally Japanese scheme is universally Japanese,
and that the ring $\mZ$  (resp. any field) is universally Japanese  (\cite{ega4} 7.7.4).

\subsection{}\label{not2}
Let $X$ be a $k$-scheme.  
We denote by $\pr_1,\pr_2\colon X\times_kX\rightarrow X$ the canonical projections. 
If $Y$ is an $X$-scheme and $Z$ is an $(X\times_kX)$-scheme, we denote by $Y\times_XZ$ (resp. $Z\times_XY$)
the fibered product of $Y$ and $Z$ over $X$, 
where $Z$ is considered as an $X$-scheme by $\pr_1$ (resp. $\pr_2$). 
In particular, in $Z\times_XZ$, the first factor 
is considered as an $X$-scheme by $\pr_2$ while the second factor is considered as an $X$-scheme by $\pr_1$. 
Let $\cF$ be an étale sheaf of $\Lambda$-modules on $X$. We denote by $\cH(\cF)$ the sheaf on $X\times_kX$
defined by 
\begin{equation}\label{not2a}
\cH(\cF)=\cHom(\pr_2^*\cF,\pr_1^*\cF).
\end{equation} 
If $f\colon Y\rightarrow X$ is a morphism of schemes, we denote (abusively) the pull-back $f^*(\cF)$ also by $\cF|Y$. 

\begin{lem}\label{rc2}
Consider a commutative diagram of finite morphisms of locally noetherian schemes
\begin{equation}
\xymatrix{
Z'\ar[r]^h\ar[d]_{i'}&Z\ar[d]^i\\
Y'\ar[r]^{g'}\ar[d]_{f'}&Y\ar[d]^f\\
X'\ar[r]^g&X}
\end{equation}
and let $X_0$ be a dense open subscheme of $X$, $z'\in Z'$, $y'=i'(z')$, $x'=f'(y')$, $z=h(z')$, $y=i(z)=g'(y')$
and $x=f(y)=g(x')$. 
We denote by the index $_0$ the base change of schemes or morphisms by the canonical injection $X_0\rightarrow X$.
Assume that $X,X',Y$ and $Y'$ are normal, that $Y_0$ is dense in $Y$, that $f_0$ is étale,
that $Y'_0\simeq Y_0\times_XX'$, that $Y'_0$ is dense in $Y'$ and that 
$f\circ i$ and $f'\circ i'$ are closed immersions. 

{\rm (i)}\ If $f$ is étale at $y$,  then $f'$ is étale at $y'$. 

{\rm (ii)}\ Assume moreover that the irreducible component of $X'$ containing $x'$ dominates the irreducible 
component of $X$ containing $x$, that  $Z'_0\simeq Z_0\times_XX'$ and that $Z'_0$ is schematically dense in $Z'$. 
Then $f$ is étale at $y$ if and only if $f'$ is étale at $y'$. 
\end{lem}

(i) We denote by $V$ (resp. $V'$) the maximal open subscheme of $Y$ (resp. $Y'$) where $f$ (resp. $f'$)
is étale. Since $V\times_XX'$ is étale over $X'$, it is normal. 
But $Y'$ is the integral closure of $Y\times_XX'$ in $Y'_0$.
Therefore, $V\times_XX'$ is isomorphic to $g'^{-1}(V)$, and $g'^{-1}(V)\subset V'$, which implies the proposition. 

(ii) Observe first that $i$ and $i'$ are closed immersions.
By (i), it is enough to prove that if $f'$ is étale at $y'$, then $f$ is étale at $y$. 
We may replace $X$ (resp. $X'$) by its strict henselization 
at a geometric point above $x$ (resp. $x'$) and $Y$ and $Z$ (resp. $Y'$ and $Z'$)
by their pull-back~; so we may assume $X,X', Z$ and $Z'$ strictly local. 
Then we may replace $Y$ by its localization $Y_y$ and $Y'$ by $g'^{-1}(Y_y)$. 
Let $Y^\dagger$ be a connected component of $Y'$.  By assumption $Y^\dagger_0$ is dense in $Y^\dagger$. 
Since the restriction $Y^\dagger_0\rightarrow X'_0$ of $f'$ is finite and étale, it is surjective~;
hence $f'(Y^\dagger)=X'$. If $F$ is a reduced closed subscheme of $Y$ such that $f(F)=X$,
then $F=Y$.  We deduce that $g'(Y^\dagger)=Y$. Since $Z'_0$ is dense in $Z'$, 
it is not empty.  Then $g'^{-1}(g'(Z'_0))\cap Y^\dagger =Z'_0\cap Y^\dagger\not=\emptyset$, and hence $Z'\subset Y^\dagger$.
Therefore $Y'$ is connected and $f'$ is an isomorphism (as it is étale).
Let $\xi$ (resp. $\eta$) be the generic point of $X$ (resp. $Y$). 
It follows that $f$ induces an isomorphism $\eta\simeq \xi$.
Since $Y$ is the integral closure of $X$ in $\eta$, $f$ is an isomorphism.

\begin{prop}\label{ram-cov5}
Let $X$ be a regular, locally noetherian and universally Japanese scheme, $U$ an open dense subscheme of $X$, 
$V$ a finite étale covering of $U$, $Y$ the integral closure of $X$ in $V$, 
$V'$ the maximal open subscheme of $Y$
which is étale over $X$, and $T$ a closed subscheme of $Y$. Assume the following conditions satisfied~:

{\rm (i)}\  All codimension one points of $T$ are contained in $V'$.

{\rm (ii)}\   There exists a Galois torsor $W$ over $U$, with nilpotent group $G$,
and a subgroup $H$ of $G$, such that $V$ is $U$-isomorphic to the quotient of $W$ by $H$. 

Then $T\subset V'$.
\end{prop}

Let $Z$ be the integral closure of $X$ in $W$. 
For every geometric point $\oz$ of $Z$, we denote by $I_{\oz}\subset G$ the inertia group of $\oz$.
By condition (i), if $\oz$ is above a codimension one point of $T$, then $I_{\oz}\subset H$.
We proceed by induction on $[G:H]$. Let $n\geq 1$. 
We assume that the proposition holds true if $[G\colon H]<n$ and prove it if $[G\colon H]=n$. 
The proposition is obvious if $G=H$; so we may assume that $n>1$. 
There exists a normal subgroup $G'$ of $G$ 
containing $H$ and different from $G$ such that $G/G'$ is abelian (\cite{alg} I §6.3 prop.~8). 
Observe that $G'$ is nilpotent. We denote by $U'$ the 
quotient of $W$ by $G'$, and by $X'$ the integral closure of $X$ in $U'$. 
We denote by $h\colon Y\rightarrow X'$ and $g\colon X'\rightarrow X$ the canonical morphisms, 
and put $f=g\circ h$. Let $N$ be the maximal open subscheme of $X$ over which $g$ is étale (observe that $g$ is finite). 
It follows from the assumption that $g$ is étale at all points $h(t)\in X'$, where $t$ is a codimension one point of $T$.
Since $g$ is Galois, we conclude that for any codimension one point $t$ of $T$, $f(t)\in N$.
Let $S$ be the reduced closed subscheme of $X$ with support $X-N$. 
By the Zariski-Nagata purity theorem (\cite{sga2} X 3.4), $S$ is a Cartier divisor on $X$. 
If $f(T)\cap S\not=\emptyset$, then there exists a codimension one point $t$ of $T$
such that $f(t)\in S$, which is a contradiction. Hence $f(T)\subset N$. We may replace $X$ by $N$
and $U,V$ and $W$ by their pull-backs. Then $X'$ is étale above $X$; in particular, it is regular. 
The induction assumption implies that $h$ is étale over an open neighbourhood of $T$ in $Y$. 
Then $T\subset V'$. 

\begin{rema}
Under the assumptions of \eqref{ram-cov5}, if moreover $V$ is connected, 
then condition (ii) is equivalent to the following condition~:

(ii') $V$ is dominated by a connected finite étale Galois covering $W$ of $U$ with nilpotent Galois group $G$ 
(i.e. there exists a dominating $U$-morphism $W\rightarrow V$).
\end{rema}

\begin{lem}\label{NS3}
Let $A$ be a strictly henselian valuation ring, with fraction field $K$ and residue field of characteristic $p$, 
$\oK$ a separable closure of $K$ and $G$ the Galois group of $\oK$ over $K$. 
Then $G$ has a normal $p$-Sylow subgroup $P$, and the quotient $I^t=G/P$ is abelian. 
In particular, $G$ is solvable.
\end{lem} 

Let $L$ be a finite Galois extension of $K$ and $G_L$ the Galois group of $L/K$. 
The integral closure $B$ of $A$ in $L$ is a strictly henselian valuation ring (\cite{nagata} Theorem 9). 
Let $P_L$ be the large valuation group  of $B$ defined in (\cite{zs} VI § 12 page 75), which is a normal subgroup of $G_L$. 
Then $P_L$ is a $p$-group (loc. cit., Theorem 24 page 77). The quotient $I_L^t=G_L/P_L$ 
is abelian (by construction),  and its order is prime to $p$ (loc. cit., (23) page 76). 
Hence $P_L$ is a normal $p$-Sylow subgroup of $G_L$. 
The lemma follows by passing to the limit over finite Galois extensions of $K$ contained in $\oK$.

\begin{prop}\label{NS31}
Let $A$ be a local ring of maximal ideal $\fm$, $\fp$ a prime ideal of $A$, $\kappa(\fp)$
the residue field of $A$ at $\fp$, and $\nu\colon A\rightarrow A_\fp$ the canonical homomorphism. 
Assume that $\fp A_\fp\subset \nu(A)$, and consider the following conditions~:

{\rm (i)}\ $A$ is henselian. 

{\rm (ii)}\ $A_\fp$ and $A/\fp$ are henselian. 

Then we have {\rm (i)}$\Rightarrow${\rm (ii)}. If moreover $\nu$ is injective, the two conditions are equivalent.  
\end{prop}

Observe first that we may assume that $\nu$ is injective.

(i)$\Rightarrow$(ii). We need only to prove that $A_\fp$ is henselian. Let $B'$ be a finite free $A_\fp$-algebra.  
We need to prove that $B'$ is decomposed. By (\cite{ega4} 8.8.2 and 8.10.5), there exists $f\in A-\fp$ and a finite 
$A_f$-algebra of finite presentation $B''$ such that $B'\simeq B''\otimes_{A_f}A_\fp$. Then by Zariski's main theorem
(\cite{ega4} 8.12.6), there exists a finite $A$-algebra $B$ that fits into a commutative digram 
\[
\xymatrix{
{\Spec(B'')}\ar[r]^j\ar[d]&{\Spec(B)}\ar[d]\\
{\Spec(A_f)}\ar[r]&{\Spec(A)}}
\]
where $j$ is an open immersion. The induced morphism $\Spec(B'')\rightarrow \Spec(B_f)$ being an open and 
closed immersion, we may replace $B'$ by $B_\fp$. Replacing $B$ by its canonical image in $B_\fp$, we may
assume that $B\subset B_\fp$. We have isomorphisms of $A$-modules 
\[
B_\fp/B\simeq B\otimes_A(A_\fp/A)\simeq B\otimes_A(\kappa(\fp)/(A/\fp)),
\]
where the second follows from the assumption  $\fp A_\fp\subset A\subset A_\fp$. 
We deduce that $\fp B_\fp\subset B$. 

We put $C=B\otimes_A\kappa(\fp)=B_\fp/\fp B_\fp$. Since $C$ is an artinian ring, we have 
\[
C\simeq \prod_{\fq \in Q} C_\fq,
\]
where $Q$ is the set of prime ideals of $B$ above $\fp$.  For each $\fq\in Q$, we denote by $\tC_\fq$ 
the canonical image of $B$ in $C_\fq$. We put $\tC=\prod_{\fq \in Q} \tC_\fq$ and denote by $\tB$ its inverse image  
by the canonical morphism $B_\fp\rightarrow B_\fp/\fp B_\fp=C$.
We have an exact sequence of $A$-modules
\begin{equation}
0\rightarrow B/\fp B_\fp\rightarrow \tC  \rightarrow \tB/B\rightarrow  0.
\end{equation}
We deduce that $\tB$ is finite over $A$ and that $B_\fp\simeq \tB_\fp$. 
On the other hand, since $\fp A_\fp\subset A\subset A_\fp$,
we have $\fp=\ker(A\rightarrow \kappa(\fp))= \fp A_\fp$, and hence $\fp \tB= \fp A_\fp \tB=\fp \tB_\fp$.  
Therefore, the canonical morphism $\tB/\fp \tB\rightarrow \tC$ is an isomorphism. 

Consider the following commutative diagram 
\[
\xymatrix{
{\idemp(B_\fp)}\ar[r]^-(0.5)u&{\idemp(C)}\\
{\idemp(\tB)}\ar[r]^-(0.5)v\ar[u]\ar[dr]_{\beta}&{\idemp(\tC)}\ar[u]_\alpha\ar[d]^{\gamma}\\
&{\idemp(\tB/\fm \tB)}}
\]
where $\idemp(-)$ denotes the set of idempotents. Since $A$ is henselian, $\beta$ and $\gamma$ are bijective;
then so is $v$. 
For each $q\in Q$, since $C_\fq$ is an artinian ring, $\tC_\fq$ is a local ring and $\alpha$ is a bijection. 
We deduce that $u$ is surjective and hence that $B_\fp$ is decomposed. Note that $u$ is always injective. 

(ii)$\Rightarrow$(i) Let $B$ be a finite free $A$-algebra.  
We need to prove that $B$ is decomposed. Consider the following commutative diagram.
\[
\xymatrix{
{\idemp(B)}\ar[r]^-(0.5)u\ar[d]&{\idemp(B/\fp B)}\ar[r]^-(0.5)w\ar[d]&{\idemp(B/\fm B)}\\
{\idemp(B_\fp)}\ar[r]^-(0.5)v&{\idemp(B_\fp/\fp B_\fp)}&}
\]
By assumption, $v$ and $w$ are bijections. On the other hand, it follows from the assumption 
$\fp A_\fp\subset A\subset A_\fp$ that the canonical diagram 
\[
\xymatrix{
B\ar[r]\ar[d]&{B/\fp B}\ar[d]\\
{B_\fp}\ar[r]&{B_\fp/\fp B_\fp}}
\]
is cartesian with injective vertical arrows. We deduce that $u$ is surjective and hence that $B$ is decomposed.

\begin{defi}\label{NS1}
Let $X$ be a locally noetherian and normal scheme, $U$ a dense open subscheme of $X$, 
$V$ a Galois torsor over $U$ of group $G$, $Y$ the integral closure of $X$ in $V$, and
$\ox$ a geometric point of $X$. We say that $V/U$ has the property (NpS)
at $\ox$ if for every geometric point $\oy$ of $Y$ above $\ox$, the inertia group 
$I_{\oy}$ of $\oy$ has a normal  $p$-Sylow subgroup (or equivalently, $I_{\oy}$ is a semi-direct product
of a group of order prime to $p$ by a $p$-group (\cite{serre} theorem 4.10)).  
\end{defi}

\begin{lem}\label{NS19}
Let $X$ be a normal, locally noetherian and universally Japanese scheme, $U$ a dense open subscheme of $X$, 
$V$ a Galois torsor over $U$ of group $G$, $Y$ the integral closure of $X$ in $V$, 
$\oy$ and $\oy'$ geometric points of $Y$, and $I_\oy$ and $I_{\oy'}$ the inertia 
groups of $\oy$ and $\oy'$, respectively. 
If $\oy$ is a specialization of $\oy'$, then $I_{\oy'}\subset I_{\oy}$.
\end{lem}

Let $Y_{(\oy)}$ and $Y_{(\oy')}$ be the strict localizations of $Y$ at $\oy$ and $\oy'$, respectively, and
let $v\colon Y_{(\oy')}\rightarrow Y_{(\oy)}$ be a specialization map. 
Let $\ox$ and $\ox'$ be the images of $\oy$ and $\oy'$ in $X$, respectively, and let $X_{(\ox)}$ and $X_{(\ox')}$ be 
the corresponding strict localizations of $X$. There exists a specialization map 
$u\colon X_{(\ox')}\rightarrow X_{(\ox)}$ such that the diagram 
\begin{equation}
\xymatrix{
{Y_{(\oy')}}\ar[r]^v\ar[d]&{Y_{(\oy)}}\ar[d]\\
{X_{(\ox')}}\ar[r]^u&{X_{(\ox)}}}
\end{equation}
where the vertical arrows are the canonical morphisms, is commutative (\cite{sga4} VIII 7.4). The morphism 
\begin{equation}
w=u\times_XV\colon X_{(\ox')}\times_XV\rightarrow X_{(\ox)}\times_XV
\end{equation}
is $G$-equivariant. If we identify $Y_{(\oy)}\times_YV$ with a connected component of $X_{(\ox)}\times_XV$
and $Y_{(\oy')}\times_YV$ with a connected component of $X_{(\ox')}\times_XV$ \eqref{NS15}, then we have 
$w(Y_{(\oy')}\times_YV)\subset Y_{(\oy)}\times_YV$. We deduce that $I_{\oy'}\subset I_{\oy}$ \eqref{NS15}.

\begin{cor}\label{NS2}
Let $X$ be a normal, locally noetherian and universally Japanese scheme, $U$ a dense open subscheme of $X$, 
$V$ a Galois torsor over $U$ of group $G$, and $\ox$ a geometric point of $X$.  
Assume that $V/U$ has the property {\rm (NpS)} at $\ox$. 
Then, there exists an open neighborhood $X_0$ of $\ox$ in $X$ such that
$V/U$ has the property {\rm (NpS)} at every geometric point of $X_0$.
\end{cor}

This follows from \ref{NS19}, (\cite{sga4} VIII 7.5) and the fact that if a finite group has a normal 
$p$-Sylow subgroup, then so is any subgroup (\cite{alg} I § 6.6 cor.~3 of theo.~3).

\begin{lem}\label{NS25}
Let $X,X'$ be normal, locally noetherian and universally Japanese schemes, $U$ a dense open subscheme of $X$, 
$V$ a Galois torsor over $U$ of group $G$, $Y$ the integral closure of $X$ in $V$, 
$f\colon X'\rightarrow X$ a morphism, $U'=f^{-1}(U)$, $V'=U'\times_UV$, 
$Y'$ the integral closure of $X'$ in $V'$, $g\colon Y'\rightarrow Y$ the canonical morphism, 
$\oy'$ a geometric point of $Y'$, $\oy=g(\oy')$, $\ox'$ the image of $\oy'$ in $X'$, $\ox=f(\ox')$,
and $I_\oy$ and $I_{\oy'}$ the inertia groups of $\oy$ and $\oy'$, respectively. Then $I_{\oy'}\subset I_{\oy}$.
In particular, if $V/U$ has the property {\rm (NpS)} at $\ox$, $V'/U'$ has the property {\rm (NpS)} at $\ox'$.
\end{lem}

Replacing $X'$ by the schematic closure of $U'$ in $X'$, we may assume that $U'$ is dense in $X'$. 
Let $Y_{(\oy)}$ (resp. $Y'_{(\oy')}$) be the strict localization of $Y$ at $\oy$ (resp. $Y'$ at $\oy'$) and let
$X_{(\ox)}$ (resp. $X'_{(\ox')}$) be the strict localization of $X$ at $\ox$ (resp. $X'$ at $\ox'$). The morphism 
\begin{equation}
h=f\times_XV\colon X'_{(\ox')}\times_{X'}V'\rightarrow X_{(\ox)}\times_XV
\end{equation}
is $G$-equivariant. If we identify $Y_{(\oy)}\times_YV$ with a connected component of $X_{(\ox)}\times_XV$
and $Y'_{(\oy')}\times_{Y'}V'$ with a connected component of $X'_{(\ox')}\times_{X'}V'$ \eqref{NS15}, then we have 
$h(Y'_{(\oy')}\times_{Y'}V')\subset Y_{(\oy)}\times_YV$. We deduce that $I_{\oy'}\subset I_{\oy}$ \eqref{NS15}.
The second assertion follows immediately from the first one. 

\begin{prop}\label{AB1}
Let $A$ be a ring, $t\in A$, $X=\Spec(A)$, $U=\Spec(A_t)$, $B$ a finite sub-$A$-algebra of $A_t$ and 
$Y=\Spec(B)$. Assume that $t$ is not a zero divisor in $A$. 
Then the canonical morphism $Y\rightarrow X$ is a $U$-admissible blow-up. 
\end{prop}

We refer to (\cite{rg} § 5.1 and \cite{egr1} §1.13) for generalities on admissible blow-ups. 
Let $f_i$ $(1\leq i\leq n)$ be generators of the $A$-algebra $B$, $a_i$ $(1\leq i\leq n)$ elements 
of $A$, and $r$ an integer $\geq 1$ such that $a_i=t^rf_i \in A_t$. We put $a_0= t^r$, $I=(a_0,a_1,\dots,a_n)$
and let $\varphi\colon X'\rightarrow X$ be the blow-up of $I$ in $X$. For $0\leq i\leq n$, we put
\begin{eqnarray*}
A'_i&=&A\left[\frac{a_0}{a_i},\dots, \frac{a_n}{a_i}\right],\\
A_i&=&A'_i/J_i,
\end{eqnarray*}
where $J_i$ is the ideal of $a_i$-torsion in $A'_i$ (i.e. the ideal of $x\in A'_i$ such that $a_i^mx=0$ for some $m\geq 1$). 
We see easily that the $\Spec(A_i)$'s $(0\leq i\leq n)$ form an open covering of $X'$; 
for every $0\leq i\leq n$, $\Spec(A_i)$ is the maximal open subscheme of $X'$ where $\varphi^*(a_i)$ 
generates the ideal $I\co_{X'}$ (cf. \cite{egr1} 3.1.6 and 3.1.7). It is clear that $B=A_0$; so
$Y$ is canonically identified with an open subscheme of $X'$. 
Since $Y$ is finite over $X$, the open immersion $Y\rightarrow X'$ is also closed. 
On the other hand, as $U$ is schematically dense in $X$ (\cite{egr1} 1.8.30.2), $\varphi^{-1}(U)$
is schematically dense in $X'$ (\cite{egr1} 1.13.3(i)). But $\varphi^{-1}(U)\subset Y$, so $Y=X'$.

\begin{cor}\label{AB2}
Let $X$ be a quasi-compact and quasi-separated scheme, $D$ an effective Cartier divisor on $X$, 
$U=X-D$ and $f\colon Y\rightarrow X$ a finite morphism 
inducing an isomorphism above $U$ such that $f^{-1}(U)$ is schematically dense in $Y$. 
Then, there exists a $U$-admissible blow-up $\varphi\colon X'\rightarrow X$ 
and an $X$-morphism $g\colon X'\rightarrow Y$. 
\end{cor}

Let $X_i=\Spec(A_i)$ $(1\leq i\leq n)$ be a finite affine open covering of $X$ such that, for each $i$, 
$D$ is defined over $X_i$ by one equation in $A_i$.
For each $1\leq i\leq n$, we put $Y_i=X_i\times_XY$ and let $f_i\colon Y_i\rightarrow X_i$
be the restriction of $f$. By \ref{AB1}, each $f_i$ is a $(U\cap X_i)$-admissible blow-up. 
By (\cite{rg} 5.3.1), there exists a $U$-admissible blow-up $\varphi_i\colon X'_i\rightarrow X$
extending $f_i$. Assume that $\varphi_i$ is the blow-up of an ideal of finite type $\cA_i$ of $\co_X$
such that $\cA_i|U=\co_X|U$. Let $\varphi\colon X'\rightarrow X$ be the blow-up of $\prod_{i=1}^n\cA_i$. 
By the universal property of blow-ups, for each $1\leq i\leq n$, there exists an $X$-morphism 
$h_i\colon X'\rightarrow X'_i$. Its restriction above $X_i$ is a morphism $g_i\colon X'\times_XX_i\rightarrow Y_i$.
For each $1\leq i,j\leq n$, the restrictions of $g_i$ and $g_j$ above $X_i\cap X_j$ are canonically identified. 
By gluing the $g_i$'s, we get an $X$-morphism $g\colon X'\rightarrow Y$.

\subsection{}\label{ZR1}
Let $X$ be a coherent scheme (i.e. a quasi-compact and quasi-separated scheme) and
$U$ an open subscheme of $X$.  We denote by $\Sch_{/X}$ the category of $X$-schemes and 
by $\cB$ the full subcategory of $\Sch_{/X}$ of objects $(X',\varphi)$,
where $\varphi\colon X'\rightarrow X$ is a $U$-admissible blow-up (cf. \cite{rg} § 5.1 and \cite{egr1} §1.13).
The {\em Zariski-Riemann space} of the pair $(X,U)$ is the topological space defined by
\begin{equation}\label{ZR1a}
X_{\ZR}=\underset{\underset{(X',\varphi)\in \cB}{\longleftarrow}}{\lim}\ |X'|,
\end{equation}
where $|X'|$ denotes the topological space underlying to $X'$. For every $\xi\in X_\ZR$, we put 
\begin{equation}\label{ZR1c}
\co_{X_\ZR,\xi}=\underset{\underset{(X',\varphi)\in \cB^\circ}{\longrightarrow}}{\lim}\  \co_{X',\xi_\varphi},
\end{equation}
where $\xi_\varphi$ is the image of $\xi$ in $X'$. 
By (\cite{zs} VI §17), $X_\ZR$ is quasi-compact. 
If $U$ is schematically dense in $X$, then the canonical map
$X_\ZR\rightarrow |X|$ is surjective. 

\subsection{} \label{ZR2}
Let $X$ be a coherent scheme, $D$ an effective Cartier divisor on $X$ and
$U=X-D$. We keep the notation of \eqref{ZR1} and  
denote by $\cC$ the full subcategory of $\Sch_{/X}$ of objects $(X'',\psi)$, where $\psi$ is composed of
two morphisms $X''\stackrel{\rho}{\rightarrow}X'\stackrel{\varphi}{\rightarrow}X$ satisfying 
the following conditions~:

(a) $\varphi$ is a $U$-admissible blow-up; 

(b) $\rho$ is a finite morphism inducing an isomorphism above $U$;

(c) $\psi^{-1}(U)$ is schematically dense in $X''$.

Then every object of $\cB$ is an object of $\cC$ (\cite{egr1} 1.13.3(i)). We denote by  
\begin{equation}\label{ZR2a}
\iota\colon \cB\rightarrow \cC
\end{equation}
the canonical injection functor. Then $\iota^\circ$ is cofinal (cf. \cite{sga4} I 8.1.1) and $\cC$ is cofiltered.
Indeed, since $\cB$ is cofiltered, it is enough to prove that $\iota^\circ$ satisfies conditions F$1)$ and F$2)$ 
of (\cite{sga4} I 8.1.3). Condition F$1)$ follows from \ref{AB2} and (\cite{rg} 5.1.4), 
and condition F$2)$ is an immediate consequence of condition (c) above. 
We deduce that the canonical morphism 
\begin{equation}\label{ZR2b}
\underset{\underset{(X'',\psi)\in \cC}{\longleftarrow}}{\lim}\ |X''| \rightarrow
\underset{\underset{(X',\varphi)\in \cB}{\longleftarrow}}{\lim}\ |X'|=X_\ZR
\end{equation}
is an isomorphism. For every $\xi\in X_\ZR$, we have a canonical isomorphism 
\begin{equation}\label{ZR2c}
\co_{X_\ZR,\xi} \stackrel{\sim}{\rightarrow} \underset{\underset{(X'',\psi)\in \cC^\circ}{\longrightarrow}}{\lim}\  
\co_{X'',\xi_\psi},
\end{equation}
where $\xi_\psi$ is the image of $\xi$ in $X''$ by the map $X_\ZR\rightarrow |X''|$ induced by \eqref{ZR2b}.
Note that for any object $(X'',\psi)$ of $\cC$, the map $X_\ZR\rightarrow |X''|$ is surjective (cf. \ref{AB2} and \ref{ZR1}).

\begin{lem}\label{ZR3} 
Let $X$ be a coherent scheme, $D$ a closed subscheme of finite presentation of $X$, 
$U=X-D$, $x\in D$, $\xi$ a point of $X_\ZR$ above $x$, and
$J$ the ideal of $\co_{X,x}$ defined by $D$. Assume that $U$ is schematically dense in $X$,
and put $\co_\xi=\co_{X_\ZR,\xi}$ and $\fp=\bigcap_nJ^n \co_{\xi}$. Then~:

{\rm (i)}\ $\co_{\xi}$ equipped with the $J$-adic topology, is a prevaluative  ring, 
which means that it is local and that every open ideal of finite type is invertible {\rm (\cite{egr1}  1.9.1)}.
Let $t\in \co_{\xi}$ be a generator of $J \co_{\xi}$. 

{\rm (ii)}\ $\co_{\xi}[\frac 1 t]$ is a local ring. 

{\rm (iii)}\ $\co_{\xi}/\fp$ is a valuation ring with fraction field the residue field of $\co_{\xi}[\frac 1 t]$. 
In particular, $\co_\xi[\frac 1 t]$ is the localization of $\co_\xi$ at $\fp$.

{\rm (iv)}\ The ideal $\fp\co_\xi[\frac 1 t]$ is contained in the image of the canonical homomorphism 
$\co_\xi\rightarrow \co_\xi[\frac 1 t]$.
\end{lem}

Since the transition homomorphisms of the inductive limit \eqref{ZR1c}
are local, the ring $\co_{\xi}$ is local and $J$ is contained in the maximal ideal of $\co_{\xi}$.
Replacing $X$ by its blow-up along $D$, we may assume that $J$ is invertible, generated by $t\in \co_{X,x}$.  
For every object $(X',\varphi)$ of $\cB$ \eqref{ZR1},  $\varphi^{-1}(U)$ is schematically dense in $X'$ (\cite{egr1} 1.13.3(i)). 
It follows that $t$ is not a zero divisor in $\co_{\xi}$.  
Let $I$ be an open ideal of finite type of $\co_{\xi}$. Then $I$ is induced by an ideal of finite type 
$\cI$ of $\co_{X'}$ for an object $(X',\varphi)$ of $\cB$ such that $\cI|U=\co_U$. By blowing-up $\cI$ in $X'$, 
we obtain an object $(X_1,\phi)$ of $\cB$ (\cite{rg} 5.1.4) such that $\cI\co_{X_1}$ is invertible. 
Therefore, the ideal $I$ is monogenic. Since $I$ is open and since $t$ is not a zero divisor in $\co_{\xi}$, 
$I$ is invertible, which proves the proposition (i). 

Propositions (ii) and (iii) follow from (\cite{egr1} 1.9.4). Proposition (iv) is obvious.

\begin{lem}\label{ZR4}
Let $X$ be a normal, locally noetherian and universally Japanese scheme, 
$D$ an effective Cartier divisor on $X$, $U=X-D$ and
$V$ a Galois torsor over $U$ of group $G$. Then with the notation of \eqref{ZR1} and \eqref{ZR2}, 
for every $\xi\in X_\ZR$, there exists an object $(X'',\psi)$ of $\cC$ satisfying 
the following properties~:

{\rm (i)}\ $X''$ is normal. 

{\rm (ii)}\ If $\pi_\psi\colon X_\ZR\rightarrow |X''|$ is the morphism induced by the isomorphism \eqref{ZR2b}
and $\xi_\psi=\pi_\psi(\xi)$, then $V/U$ has property {\rm (NpS)} at every geometric point $\oxi_\psi$ of $X''$ above $\xi_\psi$.
\end{lem}

We put $\co_\xi=\co_{X_{\ZR},\xi}$ and $S=\Spec(\co_\xi)$.
Let $\os$ be a geometric closed point of $S$ and $\oS=\Spec(A)$ the corresponding strictly local scheme. 
For each object $(X'',\psi)$ of $\cC$, we denote by $\pi_\psi\colon X_\ZR\rightarrow |X''|$ the morphism induced by 
the isomorphism \eqref{ZR2b}, and put $\xi_\psi=\pi_\psi(\xi)$, 
$\co_{\xi_{\psi}}=\co_{X'',\xi_{\psi}}$ and $S_\psi=\Spec(\co_{\xi_{\psi}})$. 
We have a canonical isomorphism \eqref{ZR2c}
\begin{equation}
S\stackrel{\sim}{\rightarrow}\underset{\underset{(X'',\psi)\in \cC}{\longleftarrow}}{\lim}\ S_{\psi}.
\end{equation} 
Let $(X'',\psi)$ be an object of $\cC$. 
Since the canonical homomorphism $\co_{\xi_{\psi}}\rightarrow \co_\xi$ is local, $\os$ determines a geometric 
closed point of $S_\psi$ (also denoted by $\os$). 
We denote by $\oxi_\psi$ the geometric point of $X''$ above $\xi_\psi$ corresponding to $\os$,
and by $\oS_\psi$ the strict localization of $S_\psi$ at $\os$. 
Then the canonical morphisms $\oS\rightarrow \oS_\psi$ induce an isomorphism 
\begin{equation}\label{ZR4a}
\oS\stackrel{\sim}{\rightarrow}\underset{\underset{(X'',\psi)\in \cC}{\longleftarrow}}{\lim}\ \oS_\psi.
\end{equation}
Indeed, the projective limit above is a strictly local scheme, with the same residue field as $\oS$. 

We may assume that $x=\xi_\id\in D$. 
Let $J$ be the ideal of  $\co_{X,x}$ defined by $D$. We put $\fp=\bigcap_n J^n\co_\xi$. 
It follows from \ref{ZR3} that the ideal $J\co_\xi$ is invertible, 
generated by $t\in \co_\xi$, that $\co_\xi/\fp$ is a valuation ring 
and that  $\co_\xi[\frac 1 t]$ is the localization of $\co_\xi$ at $\fp$. 
The scheme $\oT=\Spec(A/\fp A)$ is the strict localization of $T=\Spec(\co_\xi/\fp)$ at $\os$.   
Since $\oS$ has only one point above $\fp\in S$ (\cite{sga4} VIII 7.6), namely $\fp A$, 
we have $A_{\fp A}=A\otimes_{\co_\xi}(\co_\xi)_\fp=A[\frac 1 t]$. Then it follows from \ref{ZR3}(iv) 
that $\fp A_{\fp A}$ is contained in the image of the canonical homomorphism $A\rightarrow A_{\fp A}$. 
Therefore, $A[\frac 1 t]$ is a henselian local ring by \ref{NS31}, and hence the canonical map
\begin{equation}\label{ZR4b}
\pi_0(\oT\times_XV)\rightarrow \pi_0(\oS\times_XV)
\end{equation}
is bijective. 

Since $\co_\xi/\fp$ is a valuation ring, $A/\fp A$ is a strictly henselian valuation ring. 
This follows from (\cite{blr-nm} § 2.4, prop.~11) and (\cite{nagata} Fundamental lemma on 
the extensions of valuations, page 50).
Let $K$ be the fraction field of $A/\fp A$, $\oK$ a separable closure of $K$ and $\cG$ the Galois group
of $\oK$ over $K$. 
By \eqref{ZR4a}, \eqref{ZR4b} and (\cite{ega4} 8.4.1), there exists an object $(X'',\psi)$ of $\cC$ such that 
the canonical map
\begin{equation}\label{ZR4c}
\pi_0(\oT\times_XV)\rightarrow \pi_0(\oS_\psi\times_XV)
\end{equation}
is injective. We may assume that $X''$ is normal. 
Since the map \eqref{ZR4c} is $G$-equivariant and since $G$ acts transitively on the source and on the target, 
it is bijective. Let $Y''$ be the normalization of $X''$ in $V$. 
The set $\pi_0(\oS_\psi\times_XV)$ is isomorphic to $Y''\otimes_{X''}\kappa(\oxi_\psi)$ (cf. \eqref{NS15a}).
For $\oy''\in Y''\otimes_{X''}\kappa(\oxi_\psi)$, we denote by $I_{\oy''}\subset G$ the inertia group of $y''$.
As $I_{\oy''}$ is the stabilizer in $G$ of the connected component of $\oS_\psi\times_XV$ corresponding 
to $\oy''$ \eqref{NS15}, it follows from the bijection \eqref{ZR4c} that $I_{\oy''}$ is isomorphic to a quotient of $\cG$. 
Therefore,  by \ref{NS3}, $V/U$ has property (NpS) at $\oxi_\psi$.

\begin{prop}\label{NS4}
Let $X$ be a normal, locally noetherian and universally Japanese scheme, 
$U$ a dense open subscheme, and $V$ a Galois torsor over $U$. 
Then there exists a $U$-admissible blow-up $\varphi\colon X'\rightarrow X$, such that if we denote by 
$X''$ the normalization of $X'$, $V/U$ has the property {\rm (NpS)} at every geometric point of $X''$. 
\end{prop}

Replacing $X$ by a  $U$-admissible blow-up, we may assume that there exists an effective Cartier divisor
$D$ on $X$ such that $U=X-D$; we take again the notation of \eqref{ZR1} and \eqref{ZR2}.
Let $\cD$ be the full subcategory of $\cC$ of objects $(X'',\psi)$ such that $X''$ is normal. 
It follows from (\cite{sga4} I 8.1.3(c)) that $\cD^\circ$ is cofinal in $\cC^\circ$. 
For each object $(X'',\psi)$ of $\cD$, we denote by $\pi_\psi\colon X_\ZR\rightarrow |X''|$ the morphism induced by 
\eqref{ZR2b}, and by $X''_0$ the maximal open subscheme of $X''$ such that $V/U$ has the property (NpS) 
at every geometric point of $X''_0$ \eqref{NS2}. By \ref{ZR4}, 
we have 
\begin{equation}
X_\ZR=\bigcup_{(X'',\psi)\in \cD}\pi_\psi^{-1}(X''_0).
\end{equation}
Since $X_\ZR$ is quasi-compact \eqref{ZR1}, there exists an object $(X'',\psi)$ of $\cD$ such that 
$X_\ZR=\pi_\psi^{-1}(X''_0)$. As $\pi_\psi$ is surjective, we deduce that $V/U$ has the property (NpS) 
at every geometric point of $X''$.

\begin{lem}\label{ram-cov9}
Let $G,G'$ be smooth connected group schemes over $\ok$, and let $f\colon G\rightarrow G'$ be an étale 
morphism of $\ok$-group schemes. Then~:

{\rm (i)}\ $f$ is finite and surjective. 

{\rm (ii)}\ If $G'$ is commutative, then so is $G$; if $G'$ is isomorphic to $\mA_{\ok}^n$ for some integer $n\geq 1$, 
then so is $G$ and the kernel of $f$ is a finite dimensional $\mF_p$-vector space.
\end{lem}

(i) The proposition follows from (\cite{sga3} VI$_{\rm B}$ 1.3.2 and 1.4.1). 

(ii) Assume first that $G'$ is commutative. Then the derived group $(G,G)$ is contained in the kernel of $f$ and hence is the unit group. Therefore $G$ is commutative. Assume next that $G'\simeq \mA_{\ok}^n$. 
Any maximal torus of $G$ is contained in the kernel of $f$, and hence is the unit group.
Therefore, $G$ is unipotent (\cite{sga3} XVII 4.1.1). Since $pG$ is contained in the kernel of $f$ and $G$
is connected, we deduce that $pG=0$. Therefore, $G$ is isomorphic to $\mA_{\ok}^n$.

\subsection{}\label{ram1-2}
Let $f\colon X\rightarrow Y$ be a morphism of schemes, $\cF$ a locally constant and constructible 
sheaf of $\Lambda$-modules on $Y$, and $\cG$ a sheaf of $\Lambda$-modules on $Y$.
Then the canonical morphism
\begin{equation}\label{ram1-2a}
f^*(\cHom(\cF,\cG))\rightarrow \cHom(f^*(\cF),f^*(\cG))
\end{equation}
is an isomorphism. 
Indeed, the statement is obviously true if $f$ is étale (even is $\cF$ is not locally constant and constructible). 
Hence, by replacing $Y$ by an étale covering $Y'$ and $X$ by $X\times_YY'$, we may assume
that $\cF$ is constant on $Y$ of value a finite $\Lambda$-module $M$.  Since $\Lambda$ is noetherian, 
we have an exact sequence $\Lambda^m\rightarrow \Lambda^n\rightarrow M\rightarrow 0$. 
We deduce the following commutative diagram with exact lines.
\begin{equation}
\xymatrix{
0\ar[r]&{f^*(\cHom(\cF,\cG))}\ar[r]\ar[d]^\alpha&{f^*(\cG^n)}\ar[r]\ar[d]^\beta&{f^*(\cG^m)}\ar[d]^\gamma\\
0\ar[r]&{\cHom(f^*(\cF),f^*(\cG))}\ar[r]&{f^*(\cG)^n}\ar[r]&{f^*(\cG)^m}}
\end{equation}
Since $\beta$ and $\gamma$ are clearly isomorphisms, $\alpha$ is an isomorphism. 

\begin{prop}\label{ram00}
Let $X$ be a normal scheme, $U$ a dense open subscheme of $X$, 
$\cF$ a locally constant and constructible sheaf of $\Lambda$-modules on $U$, and
$f\colon Y\rightarrow X$ a morphism. 
We put $V=f^{-1}(U)$ and denote by $i\colon U\rightarrow X$ and  $j\colon V\rightarrow Y$
the canonical injections. We assume that $V$ is schematically dense in $Y$. 
Then the base change morphism 
\begin{equation}\label{ram00b}
\alpha\colon i_*(\cF)|Y\rightarrow j_*(\cF|V)
\end{equation}
relatively to $f$ is injective.
\end{prop}

Let $\oy$ be a geometric point of $Y$ and $\ox=f(\oy)$. It is enough to prove that the stalk $\alpha_{\oy}$ of 
$\alpha$ at $\oy$ is injective. We may replace $X$ and $Y$ by their strict henselizations at $\ox$ and $\oy$. 
Since $V$ is dense in $Y$, it is not empty. Let $\oeta$ be a geometric point of $V$, $\oeta\rightarrow \oy$ 
a specialization map and $\oxi=f(\oeta)$.  Then we have a commutative diagram 
\begin{equation}
\xymatrix{
{(i_*(\cF))_{\ox}}\ar[r]^{u}\ar[d]_{\alpha_\oy}&{\cF_\oxi}\ar[d]^{\alpha_{\oeta}}\\
{(j_*(\cF|V))_{\oy}}\ar[r]^v&{(\cF|V)_{\oeta}}}
\end{equation}
where $u$ and $v$ are the specialization homomorphisms. Since $\alpha_\oeta$ is an isomorphism, 
it is enough to prove that $u$ is injective. 
Since $X$ is normal and strictly local, $U$ is connected and 
we have $(i_*(\cF))_{\ox}=\Gamma(X,i_*(\cF))=\Gamma(U,\cF)$. 
There exists a connected Galois torsor $U'$ over $U$ 
that trivializes $\cF$. Then $u$ is identified with the canonical morphism
$\Gamma(U,\cF)\rightarrow \Gamma(U',\cF)$, which is obviously injective. This concludes the proof. 

\begin{lem}\label{ram0}
Let $X$ be a normal, locally noetherian and universally Japanese scheme,  
$U$ a dense open subscheme of $X$, 
$j\colon U\rightarrow X$ the canonical injection, $V$ a Galois torsor over $U$ of group $G$, 
and $M$ a $\Lambda[G]$-module. The constant étale sheaf $M_V$ on $V$ defines by Galois 
descent a locally constant and constructible sheaf $\cF$ of $\Lambda$-modules on $U$. 
Let $s\in M$, $H$ the stabilizer of $s$ in $G$, $U'$ the quotient of 
$V$ by $H$, $X'$ the integral closure of $X$ in $U'$, and $j'\colon U'\rightarrow X'$ the canonical injection.  
\begin{equation}\label{ram0a}
\xymatrix{
U'\ar[r]^{j'}\ar[d]&X'\ar[d]\\
U\ar[r]^j&X}
\end{equation}
We put $\cF'=\cF|U'$ and consider $s$ as a section of $j'_*(\cF')(X')=\cF'(U')=\cF(U')=\cF(V)^H$.
Let $\ox'$ be a geometric point of $X'$ and $\ox$ its image in $X$. 
Then the base change morphism
\begin{equation}\label{ram0c}
\alpha\colon j_*(\cF)|X'\rightarrow j'_*(\cF')
\end{equation}
relatively to the Cartesian diagram \eqref{ram0a} is injective. 
Moreover, the following conditions are equivalent~:

{\rm (i)}\  The stalk
\begin{equation}\label{ram0b}
\alpha_{\ox'}\colon (j_*\cF)_{\ox}\rightarrow (j'_*\cF')_{\ox'}
\end{equation} 
of the morphism $\alpha$ at $\ox'$ is an isomorphism. 

{\rm (ii)}\  The image of $s$ in $(j'_*\cF')_{\ox'}$ is in the image of the morphism \eqref{ram0b}. 

{\rm (iii)}\ The morphism $X'\rightarrow X$ is étale at $\ox'$.  
\end{lem} 

Observe first that the implications (iii)$\Rightarrow$(i)$\Rightarrow$(ii) are obvious. 
Let $Y$ be the integral closure of $X$ in $V$, $Y\rightarrow X'$ the canonical morphism,
$\oy$ a geometric point of $Y$ above $\ox'$, and
$I\subset G$ the inertia group of $\oy$.  
Then the morphism \eqref{ram0b} can be canonically identified with the canonical injection
\begin{equation}\label{ram0d}
M^I\rightarrow M^{I\cap H},
\end{equation}
which proves the first assertion. 
It follows from \eqref{ram0d} that 
condition (ii) is equivalent to $I\subset H$, which is also equivalent to each of the conditions (i) and (iii).

\begin{lem}\label{ram-loc4}
Let $X$ be a scheme, $U$ an open subscheme of $X$,
$\cF$ a locally constant constructible sheaf of $\Lambda$-modules on $X$, 
$\ox$ a geometric point of $X$,
$X_{(\ox)}$ the corresponding strictly local scheme,  and
$V_1$ a Galois torsor over $U_1=X_{(\ox)}\times_XU$ trivializing $\cF|U_1$. 
Then, there exists an étale morphism $f\colon X'\rightarrow X$, a geometric point $\ox'$
above $\ox$ and a Galois torsor $V'$ over $U'=f^{-1}(U)$ trivializing 
$\cF|U'$ such that if we identify the strictly local schemes $X'_{(\ox')}$ and $X_{(\ox)}$ by $f$, 
there exists a $U_1$-isomorphism $V'\times_{X'}X'_{(\ox')}\simeq V_1$. 
\end{lem}
It follows from (\cite{ega4} 8.8.2 and 10.8.5)  (cf. \cite{as4} 6.2).

\section{Additive sheaves on vector bundles} \label{add}

\begin{defi}\label{add1}
Let $X$ be a scheme with residual characteristics different from $\ell$, 
$\pi\colon E\rightarrow X$ a vector bundle, and $\cF$ a 
constructible sheaf of $\Lambda$-modules on $E$. 
We say that $\cF$ is {\em additive}  if for every geometric point $\xi$ of $X$ and for every $e\in E(\xi)$, 
denoting by  $\tau_e$ the translation by $e$ on $E_\xi=E\times_X\xi$,  
$\tau_e^*(\cF|E_\xi)$ is isomorphic to $\cF|E_\xi$. 
\end{defi}

We can make the following remarks~:

\vspace{2mm}

(i)\ We may restrict to the geometric points $\xi$ of $X$ with algebraically closed residue fields. 

(ii)\ If $\cF$ is additive, then for any $X$-scheme $X'$, 
denoting by $E'$ the vector bundle $E\times_XX'$ over $X'$, $\cF|E'$ is additive.

(iii)\ $\cF$ is additive if and only if for every geometric point $\xi$ of $X$ with algebraically closed residue field, 
$\cF|E_\xi$ is additive.

\begin{prop}\label{add15}
Let $X$ be a scheme with residue characteristics different from $\ell$, $f\colon E'\rightarrow E$ 
a morphism of vector bundles over $X$, and
$\cF$ (resp. $\cF'$) a constructible sheaf of $\Lambda$-modules on $E$ (resp. $E'$). Then~:

{\rm (i)}\ If $\cF$ is additive, $f^*(\cF)$ is additive.

{\rm (ii)}\  If $\cF'$ is additive and if $f$ surjective, $\rR^nf_!(\cF')$ is additive for all $n\geq 0$. 

{\rm (iii)}\  Assume $f$ surjective. Then $\cF$ is additive if and only if $f^*(\cF)$ is additive. 
\end{prop}

Propositions (i) and (ii) follow immediately from the definition \eqref{add1}. To prove (iii), 
it remains to show that if $f^*(\cF)$ is additive then so is $\cF$. 
The problem being local on $X$, we may assume that there exists a section $\sigma\colon E\rightarrow E'$
of $f$. Then the required property follows from (i).

\subsection{}\label{add25}
Let $\cL_\psi$ be the Artin-Schreier sheaf of $\Lambda$-modules of rank $1$ 
over the additive group $\mA^1_{\mF_p}$ over $\mF_p$,
associated to the character $\psi$ fixed in \eqref{not1} (\cite{laumon} 1.1.3).
Then $\cL_\psi$ is additive. Indeed, if $\mu\colon  \mA^1_{\mF_p}\times_{\mF_p}\mA^1_{\mF_p}
\rightarrow \mA^1_{\mF_p}$ denotes the addition, we have an isomorphism
\begin{equation}
\mu^*\cL_\psi\simeq \pr_1^*\cL_\psi\otimes \pr_2^*\cL_\psi.
\end{equation}
We will show that to a certain extent, all additive sheaves in characteristic $p$ come from $\cL_\psi$. 
If $f\colon X\rightarrow \mA^1_{\mF_p}$ is a morphism of schemes, we put $\cL_{\psi}(f)=f^*\cL_{\psi}$. 

\subsection{}\label{add3}
Let $X$ be a $k$-scheme, $\pi\colon E\rightarrow X$ a vector bundle of constant rank $d$, and
$\chpi\colon \chE\rightarrow X$ the dual vector bundle. We denote by 
$\langle \ ,\ \rangle \colon E\times_X\chE\rightarrow 
\mA^1_{\mF_p}$ the canonical pairing, by $\pr_1\colon E\times_X\chE\rightarrow E$ and
$\pr_2\colon E\times_X\chE\rightarrow \chE$ the canonical projections and by
\begin{equation}\label{add3a}
\fF_{\psi}\colon \bD_c^b(E,\Lambda)\rightarrow \bD_c^b(\chE,\Lambda)
\end{equation}
the Fourier-Deligne transform defined by 
\begin{equation}\label{add3b}
\fF_\psi(K)=\rR\pr_{2!}(\pr_1^*K\otimes \cL_\psi(\langle \ ,\ \rangle )). 
\end{equation}
We recall some properties of this transform that will be used later.

Let $\pi^{\flat}\colon  E^{\flat}\rightarrow X$ be the bidual vector bundle of $\pi\colon  E\rightarrow X$,
$a\colon E\stackrel{\sim}{\rightarrow} E^{\flat}$ the anti-canonical isomorphism defined by $a(x)=-\langle x,\ \rangle$, 
and $\fF^\vee_\psi$ the Fourier-Deligne transform for $\chpi\colon \chE\rightarrow X$.
For every object $K$ of $\bD_c^b(E,\Lambda)$, we have a canonical isomorphism (\cite{laumon} 1.2.2.1) 
\begin{equation}\label{add3c}
\fF_\psi^\vee\circ \fF_\psi(K)\simeq a^*(K)(-d)[-2d].
\end{equation}

Let $\pi'\colon E'\rightarrow X$ be a vector bundle of constant rank $d'$, 
$\fF'_\psi$ its Fourier-Deligne transform,  
$f\colon E\rightarrow E'$ a morphism of vector bundles, and
$\chf\colon \chE'\rightarrow \chE$ its dual. 
For every object $K$ of $\bD_c^b(E,\Lambda)$, we have canonical isomorphisms 
\begin{eqnarray}
\fF'_\psi \circ \rR f_!(K)&\stackrel{\sim}{\rightarrow}& \chf^*\circ \fF_\psi(K),\label{add3d}\\
\fF'_\psi \circ \rR f_*(K)(d')[2d']&\stackrel{\sim}{\rightarrow} &\chf^!\circ \fF_\psi(K)(d)[2d],\label{add3e}
\end{eqnarray}
and for every object $K'$ of $\bD_c^b(E',\Lambda)$, we have canonical isomorphisms
\begin{eqnarray}
\rR \chf_!\circ \fF'_\psi(K')(d')[2d'] &\stackrel{\sim}{\rightarrow}& \fF_\psi \circ f^*(K')(d)[2d],\label{add3f}\\
\rR \chf_*\circ \fF'_\psi(K')&\stackrel{\sim}{\rightarrow}& \fF_\psi \circ f^!(K').\label{add3g}
\end{eqnarray}
Indeed, isomorphism \eqref{add3d} is proved in (\cite{laumon} 1.2.2.4). It implies isomorphism \eqref{add3f} by \eqref{add3c},
and isomorphism \eqref{add3e} by duality and \eqref{add3c}. Finally, isomorphism \eqref{add3g} is obtained from  
\eqref{add3e} by \eqref{add3c}.

For any section $e\in E(X)$ and any object $K$ of $\bD_c^b(E,\Lambda)$, if we denote by $\tau_e\colon E\stackrel{\sim}{\rightarrow} E$ the translation by $e$, we have a canonical isomorphism (\cite{laumon} 1.2.3.2)
\begin{equation}\label{add3h}
\fF_\psi(\tau_{e*}K)\stackrel{\sim}{\rightarrow} \fF_\psi(K)\otimes \cL_\psi(\langle e,\ \rangle).
\end{equation}

\subsection{}\label{supp}
Let $X$ be a scheme and $K$ an object of $\bD_c^b(X,\Lambda)$.  
The {\em support} of $K$ is the subset of points of $X$ where the stalks of the cohomology sheaves of $K$ 
are not all zero. It is constructible in $X$. This definition is in general different from 
the one introduced in (\cite{sga4} IV 8.5.2).

\begin{prop}\label{add4}
Let $X$ be a $k$-scheme, $\pi\colon E\rightarrow X$ a vector bundle of constant rank, 
$\chpi\colon \chE\rightarrow X$ the dual vector bundle, $\cF$ a 
constructible sheaf of $\Lambda$-modules on $E$, and $S\subset \chE$ the support of $\fF_\psi(\cF)$. 
Then $\cF$ is additive if and only if for every $x\in X$,  the set $S\cap \chE_x$ is finite. 
\end{prop}

By the proper base change theorem, we may assume $X=\Spec(k)$ and $k$ algebraically closed.
Then we are reduced to the following~:

\begin{prop}\label{add41}
Assume $k$ is algebraically closed and let $E$ be a vector bundle over $k$, 
$\chpi\colon \chE\rightarrow X$ the dual vector bundle and
$\cF$ a constructible sheaf of $\Lambda$-modules over $E$. The following conditions 
are equivalent~:

{\rm (i)}\ $\cF$ is additive.

{\rm (ii)}\ The support of $\fF_\psi(\cF)$ is finite. 

{\rm (iii)}\ $\cF$ is isomorphic to a finite direct sum of sheaves of the form $M\otimes \cL_\psi(f)$, 
where $M$ is a $\Lambda$-module of finite type and $f\colon E\rightarrow \mA^1_k$ is a linear form. 

{\rm (iv)}\  $\cF$ is locally constant and all its Jordan-Hölder subquotients 
are of the form $\cL_\psi(f)\otimes_\Lambda\oLambda$,
where $\oLambda$ is the residue field of $\Lambda$ and $f\colon E\rightarrow \mA^1_k$ is a linear form. 

{\rm (v)}\  $\cF$ is locally constant and all its Jordan-Hölder subquotients are additive. 
\end{prop}

First we prove (i)$\Rightarrow$(ii). For every point $e\in E(k)$, we have \eqref{add3h}
\begin{equation}\label{add4c}
\fF_\psi(\cF)\simeq \fF_\psi(\cF) \otimes \cL_\psi(\langle e, \ \rangle).
\end{equation}
If $\cG$ is a cohomology sheaf of $\fF_\psi(\cF)$, we have $\cG\simeq \cG \otimes \cL_\psi(\langle e, \ \rangle)$
for every point $e\in E(k)$. Let $U$ be an integral locally closed subscheme of $\chE$
such that $\cG|U$ is locally constant and not zero. It is enough to prove that $U$ is a closed point of $\chE$.
We may assume that $\Lambda$ is a field. Let $\pi\colon V\rightarrow U$ be a finite étale connected covering 
such that $\pi^*(\cG|U)$ is constant. Then for every point $e\in E(k)$,
the sheaf $\pi^*(\cL_\psi(\langle e, \ \rangle)|U)$ is constant;
equivalently, for every linear form $f\colon \chE\rightarrow \mA^1_k$, the sheaf
$\cL_\psi(f)|V$ is constant. So the equation $T^p-T=f$ has a solution in the function field $k(V)$ of $V$. 
If $U$ is not a closed point of $\chE$, there exists a linear form $f\colon \chE\rightarrow \mA^1_k$ 
such that $f|U\colon U\rightarrow \mA^1_k=\Spec(k[t])$ is dominant. 
For all $c\in k^\times$, the equation $T^p-T=ct$ has a solution in $k(V)$. 
We obtain infinitely many linearly disjoint extensions of degree $p$ of $k(t)$ contained in $k(V)$, which is not possible.
So $U$ is a closed point. 

Next we prove (ii)$\Rightarrow$(iii). By  \eqref{add3c}, 
we may assume that the support of $\fF_\psi(\cF)$ is a point $i\colon \Spec(k)\rightarrow \chE$;
then it is enough to observe that for any $\Lambda$-module of finite type $M$, 
we have 
\begin{equation}\label{add4b}
\fF^\vee_\psi(i_*M)\simeq M \otimes \cL_\psi(f),
\end{equation}
which follows from \eqref{add3h} and \eqref{add3d}.

It is clear that we have (iii)$\Rightarrow$(i), (iii)$\Rightarrow$(iv)$\Rightarrow$(ii) and (vi)$\Rightarrow$(v).
Finally, since conditions (i) and (ii) are equivalent and that the latter is stable by extensions, 
we have (v)$\Rightarrow$(i).

\begin{defi}\label{add5}
Let $X$ be a $k$-scheme, $\pi\colon E\rightarrow X$ a vector bundle of constant rank, 
$\chpi\colon \chE\rightarrow X$ the dual vector bundle,
and $\cF$ an additive constructible sheaf of $\Lambda$-modules on $E$. 
We call the {\em Fourier dual support} of $\cF$ the support of $\fF_\psi(\cF)$ in $\chE$. 
We say that $\cF$ is {\em non-degenerate} if the {\em closure} of its Fourier dual support 
does not meet the zero section of $\chE$. 
\end{defi}

We can make the following remarks~:

\vspace{2mm}

(i) If we replace $\psi$ by $a\psi$ for an element $a\in \mF_p^\times$, then the Fourier dual support of $\cF$
will be replaced by its inverse image by the multiplication by $a$ on $\chE$. In particular,  
the notion of being non-degenerate does not depend on $\psi$. 

(ii)\ Let $X'$ be an $X$-scheme and $E'$ the vector bundle $E\times_XX'$ over $X'$.
Then the Fourier dual support of $\cF|E'$ is the inverse image of the Fourier dual 
support of $\cF$ (\cite{laumon} 1.2.2.9).

(iii)\ Let $f\colon E\rightarrow \mA^1_X$ be a linear form, $i\colon X\rightarrow \chE$
the associated section, $M$ a non zero $\Lambda$-module of finite type, and $d$ the rank of $E$. 
Then the Fourier dual support of 
$M\otimes \cL_\psi(-f)$ is $i(X)$. Indeed, by \eqref{add3c}, \eqref{add3d} and \eqref{add3h}, we have
\begin{equation}\label{add5a}
\fF_\psi(M \otimes \cL_\psi(-f))\simeq i_*M (-d)[-2d].
\end{equation}

(iv) Assume $X=\Spec(k)$ and $k$ algebraically closed. Then $\cF$ is locally constant, 
and its Fourier dual support is the union of the Fourier dual supports of its Jordan-Hölder subquotients \eqref{add41}.

\begin{lem}\label{add9}
Let $f\colon X'\rightarrow X$ be a finite morphism of $k$-schemes, 
$\pi\colon E\rightarrow X$ a vector bundle of constant rank,
$E'=E\times_XX'$, $\pi'\colon E'\rightarrow X'$ and $f_E\colon E'\rightarrow E$ the canonical projections, 
and $\cF'$ an additive constructible sheaf of $\Lambda$-modules on $E'$. Then $f_{E*}(\cF')$ 
is additive and its Fourier dual support is the image of the Fourier dual support of $\cF'$. 
\end{lem}

By the proper base change theorem, we may assume $X=\Spec(k)$ and $k$ is algebraically closed.
Then we may reduce the proof to the case where $X'$ is a finite disjoint sum of copies of $X$, where the assertion is obvious.

\begin{lem}\label{add7}
Let $X$ be a $k$-scheme, $\pi\colon E\rightarrow X$ a vector bundle of constant rank, 
and $\cF$ an additive constructible sheaf of $\Lambda$-modules on $E$. 
If $\cF$ is non-degenerate then $\rR \pi_*\cF=\rR \pi_!\cF=0$. 
\end{lem}

It follows from \eqref{add3c}, \eqref{add3f}
and \eqref{add3g} (applied to $f$ the zero section of the dual vector bundle $\chE$ of $E$ 
and $K'=\fF_\psi(\cF)$).

\begin{lem}\label{add8}
Let $X$ be a $k$-scheme, $\pi\colon E\rightarrow X$ a vector bundle of constant rank,
$\cG$ an additive constructible sheaf of $\Lambda$-modules on $E$, 
$\cF$ a constructible sheaf of $\Lambda$-modules on $E$, and
$u\colon \cG\rightarrow \cF$ a surjective morphism (resp. $v\colon \cF\rightarrow \cG$ 
an injective morphism). Assume that $\cF$ is locally constant on all geometric fibers of $\pi$.
Then $\cF$ is additive and its Fourier dual support is contained in the Fourier dual support of $\cG$.
\end{lem}
We may assume $X=\Spec(k)$ and $k$ algebraically closed.
Then $\cG$ is locally constant \eqref{add41}, and the assertion follows from \ref{add41} and \ref{add5}(iv). 

\begin{lem}\label{add16}
Let $X$ be a $k$-scheme, $\pi\colon E\rightarrow X$ a vector bundle of constant rank,
$G$ a group scheme over $X$, 
$\rho\colon G\rightarrow E$ an étale surjective morphism of 
group schemes over $X$, and $\cF$ a constructible sheaf of $\Lambda$-modules on $E$.
Assume that for every geometric point $\ox$ of $X$, $\rho^*(\cF)|G_\ox$ is constant; then $\cF$ is additive. 
\end{lem}

We prove first that $\rho_*(\rho^*\cF)$ is additive.
Let $\ox$ be a geometric point of $X$ and $a\in E(\ox)$.
There exists $g\in G(\ox)$ such that $\rho(g)=a$. If we denote by $\tau_a$ (resp. $\tau_g$)
the translation by $a$ on $E$ (resp. $g$ on $G$), then  
\[
\tau_a^*((\rho_*(\rho^*\cF))|E_{\ox})\simeq \rho_*(\tau_g^*((\rho^*\cF)|G_{\ox}))\simeq  (\rho_*(\rho^*\cF))|E_\ox.
\]
On the one hand, the adjunction morphism $\cF\rightarrow \rho_*(\rho^*\cF)$ is injective
and $\cF$ is locally constant on the geometric fibers of $\pi$. So $\cF$ is additive by \ref{add8}.

\begin{lem}[\cite{saito1} 2.7]\label{add10}
Let $X$ be a normal $k$-scheme, $\pi\colon E\rightarrow X$ a vector bundle of constant rank,  
$\chpi\colon \chE\rightarrow X$ the dual vector bundle, 
$U$ a dense open subscheme of $X$, $\cF$ a constructible sheaf of $\Lambda$-modules on $E$, 
and $S\subset \chE$ the support of $\fF_\psi(\cF)$. 
We put $E_U=\pi^{-1}(U)$, $\chE_U=\chpi^{-1}(U)$, $S_U=S\cap \chE_U$,
and denote by $j\colon E_U\rightarrow E$ the canonical injection. 
Assume that the following conditions are satisfied~:

{\rm (i)}\ The adjunction morphism $u\colon \cF\rightarrow j_*j^*(\cF)$ is injective. 

{\rm (ii)}\ $j^*(\cF)$ is locally constant and additive.  

{\rm (iii)}\ $\cF$ is locally constant on all fibers of $\pi$. 

Then $\cF$ is additive and its Fourier dual support $S$ is contained in the Zariski 
closure $\overline{S_U}$ of $S_U$ in $\chE$. Moreover, if $\cF$ is locally constant, then $S=\overline{S_U}$.
\end{lem}

It is enough to prove that for every $x\in X$, the set $S_x=S\cap \chE_x$ is finite and is contained in $\overline{S_U}$ 
\eqref{add41}, moreover, if $\cF$ is locally constant then $S_x=\overline{S_U}\cap \chE_x$. 
We may shrink $U$. The assertion is obvious for the generic points of $X$;
so we assume that $x$ is not a generic point of $X$. 
Let $f\colon X'\rightarrow X$ be a proper surjective morphism such that $X'$ is normal and $U'=f^{-1}(U)$
is dense in $X'$. Let $E'=E\times_XX'$, $\chE'=\chE\times_XX'$, 
$f_E\colon E'\rightarrow E$ and $f_\chE\colon \chE'\rightarrow \chE$ the canonical projections, 
$\cF'=f_E^*(\cF)$, and $S'\subset \chE'$ the support of $\fF'_\psi(\cF')$. 
We put $E'_U=f^{-1}_E(E_U)$, $\chE'_U=f^{-1}_\chE(\chE_U)$, $S'_U=S'\cap \chE'_U$
and denote by $j'\colon E'_U\rightarrow E'$ the canonical injection. 
We have $S'=f_{\chE}^{-1}(S)$, $S=f_{\chE}(S')$ and $f_\chE(\overline{S'_U})=\overline{S_U}$. 
On the other hand, the adjunction morphism $u'\colon \cF'\rightarrow j'_*j'^*(\cF')$ is composed of 
\begin{equation}\label{add10a}
\xymatrix{
{\cF'}\ar[r]^-(0.5){f_E^*(u)}&{f_E^*j_*j^*(\cF)}\ar[r]^v&{j'_*j'^*(\cF')}},
\end{equation}
where $v$ is the base change morphism relatively to $f_E$. Since $v$ is injective by \ref{ram00}, $u'$ is injective.  
Hence, it is enough to prove the assertions after replacing $X$ by $X'$ and $x$ by a point $x'$ of $X'$ above it. 
Taking for $X'$ the normalization of the blow-up of $X$ along the Zariski closure 
of $x$ in $X$, we are reduced to the case where $\co_{X,x}$ is a discrete valuation ring. 

We may assume $X$ integral. Let $\eta$ be the generic point of $X$ and $K=k(\eta)$
the residue field of $\eta$. Replacing $X$ by its normalization in a finite extension of $K$,
we may assume that $S_\eta=S\cap \chE_\eta$ is a finite set of $K$-rational points. 
After shrinking $U$, we may assume that every $f\in S_\eta$ extends to a linear form on $E_U$.  
Then there exist constant sheaves $(\cG_f)_{f\in S_\eta}$ on $E_U$ such that 
$j^*(\cF)=\oplus_{f\in S_\eta} \cL_\psi(f)\otimes \cG_f$ by \ref{add41} and (\cite{sga4} IX 2.14.1).
Let $S_{\eta,x}\subset S_\eta$ be the subset of elements $f\in S_\eta$ which are regular at $x$.
For $f\in S_{\eta,x}$, we denote by $\of\in \chE_x$ its reduction.   We claim that 
\begin{equation}\label{add10b}
j_*(\cL_\psi(f)\otimes \cG_f)|E_x=\left\{
\begin{array}{clcr}
\cL_\psi(\of)\otimes \cG_f &{\rm if} \ f\in S_{\eta,x},\\
0&{\rm if} \ f\not\in S_{\eta,x}.
\end{array}
\right.
\end{equation}
Indeed, let $f\in S_\eta$, $\ot$ be a geometric generic point of $E_x$.
It follows from (\cite{saito1} 2.8) that $\cL_\psi(f)$ is ramified at $\ot$ if and only if $f\not\in S_{\eta,x}$. 
If $f\in S_{\eta,x}$, then $j_*(\cL_\psi(f)\otimes \cG_f)|E_x=\cL_\psi(\of)\otimes \cG_f$ by  (\cite{sga4} IX 2.14.1).
On the other hand, for every geometric point $\oy$ of $E_x$ and
for every specialization map $\ot\rightarrow \oy$, the specialization homomorphism 
$j_*(\cL_\psi(f)\otimes \cG_f)_{\oy}\rightarrow j_*(\cL_\psi(f)\otimes \cG_f)_{\ot}$ is injective
(cf. the proof of \ref{ram00}). Hence  $j_*(\cL_\psi(f)\otimes \cG_f)|E_x=0$ if $f\not\in S_{\eta,x}$.

Since $\cF|E_x$ is locally constant
and the adjunction morphism $u\colon \cF\rightarrow j_*j^*(\cF)$ is injective, 
we deduce from \eqref{add10b} that the Jordan-Hölder subquotients of $\cF|E_x$ are of the form $\cL_\psi(\of)\otimes \cG_f$,
for some $f\in S_{\eta,x}$. Hence, $\cF|E_x$ is additive and $S_x\subset \{\of; f\in S_{\eta,x}\}=\overline{S_U}\cap \chE_x$ 
by  \ref{add41} and \ref{add5}(iv). 
Assume that $\cF$ is locally constant. Then $u$
is an isomorphism (\cite{sga4} IX 2.14.1), and $S_{\eta,x}=S_\eta$ as the rank of the stalks 
of $j_*j^*(\cF)$ is constant. Therefore $S_x=\{\of| f\in S_{\eta,x}\}$, which conclude the proof of the required assertion.

\begin{cor}\label{add11}
Let $X$ be a $k$-scheme, $\pi\colon E\rightarrow X$ a vector bundle of constant rank,  
$\chpi\colon \chE\rightarrow X$ the dual vector bundle, and
$\cF$ a locally constant, constructible and additive sheaf of $\Lambda$-modules on $E$. 
Then the Fourier dual support $S$ of $\cF$ is the underlying space of a closed subscheme of $\chE$
which is finite over $X$. 
\end{cor}

Let $f\colon X'\rightarrow X$ be a proper surjective morphism, $E'=E\times_XX'$, $\chE'=\chE\times_XX'$,
$f_E\colon E'\rightarrow E$ and $f_{\chE}\colon \chE'\rightarrow \chE$ the canonical projections, and $S'$ 
the Fourier dual support of $f^*(\cF)$. We have $S'=f_{\chE}^{-1}(S)$ and $S=f_{\chE}(S')$. 
First, we take for $X'$ the normalization of $X$. Since $S'$ is closed in $\chE'$ by \ref{add10}, 
$S$ is closed in $\chE$. We denote also by $S$ the reduced closed subscheme of $\chE$ 
with support $S$. We prove that $S$ is finite over $X$. We may assume 
that $X$ is normal and integral. Let $\eta$ be the generic point of $X$ and $K=k(\eta)$. 
Replacing $X$ by its normalization in a finite extension of $K$, we may assume that 
$S_\eta=S\cap \chE_\eta$ is a finite set of $K$-rational points. 

We know \eqref{add4} that $S$ is quasi-finite over $X$. 
It is enough to prove that it is proper over $X$. Let $R$ be a discrete valuation ring, $Y=\Spec(R)$,
and $y$ (resp. $\kappa$) the closed (resp. generic) point of $Y$. Consider a commutative diagram 
\begin{equation}
\xymatrix{
\kappa\ar[r]^{\beta}\ar[d]_\rho&S\ar[d]\\
Y\ar[r]^{\alpha}&X}
\end{equation}
where $\rho$ is the canonical injection. It is enough to prove that there exists an $X$-morphism
$\gamma\colon Y\rightarrow S$ such that $\beta=\gamma\circ \rho$. We may replace $X$ by any normal scheme $X'$
such that $\alpha$ factors as $Y\rightarrow X'\stackrel{f}{\rightarrow}X$, 
where $f$ is proper and surjective. 
Replacing $X$ successively by the normalization of the blow-up of $X$ along the Zariski closure of $x=\alpha(y)$ in $X$, 
we may assume that $\co_{X,x}$ is a discrete valuation ring. Then the assertion follows from the proof of \ref{add10}.
Indeed, with the notation of loc. cit., we have $S_{\eta,x}=S_{\eta}$ because $u\colon \cF\rightarrow j_*j^*(\cF)$ 
is an isomorphism.

\section{Dilatations}\label{tub}

\subsection{}\label{tub1}
Let $X$ be a scheme, $u\colon P\rightarrow X$ a morphism, 
$Y$ a closed subscheme of $P$
defined by a quasi-coherent ideal $\cI_Y$ of $\co_P$,
$R$ a closed subscheme of $X$ defined by a quasi-coherent ideal $\cJ$ of $\co_X$ and $R_Y=R\times_XY$.
\[
\xymatrix{
{R_Y}\ar[r]\ar[rd]&{Y}\ar[r]\ar[rd]&{P}\ar[d]^u\\
&R\ar[r]&X}
\]
Let $\cI$ be the ideal of $\co_P$ associated to the closed immersion   
$R_Y\rightarrow P$; we have $\cI=\cI_Y+\cJ\co_P$. 
We denote by $X'$ the blow-up of $X$ along $R$ and by $P'$ the blow-up of $P$ along $R_Y$. 
We call the maximal open subscheme $W$ of $P'$ where we have $\cJ\co_W=\cI\co_W$ 
the {\em dilatation of $P$ along $Y$ of thickening $R$} and denote it by $P^{(R)}$
(also called the dilatation of $R_Y$ in $P$ relatively to $R$ in \cite{as4} 2.4; see loc. cit. 2.8). 
There is a unique $X$-morphism $P^{(R)}\rightarrow X'$.
\[
\xymatrix{
{P^{(R)}}\ar[r]\ar[d]&{P'}\ar[r]&P\ar[d]^u\\
{X'}\ar[rr]&&X}
\]  

\begin{lem}\label{tub2}
We keep the assumptions of \eqref{tub1} and assume moreover that $R$ is a Cartier divisor on $X$
and that $\cJ_Y$ is of finite type. 
We put $U=X-R$ and denote by $j_P\colon P_U\rightarrow P$ the canonical injection. Then
$P^{(R)}$ is affine over $P$ and it corresponds to the quasi-coherent sub-$\co_P$-algebra
of $j_{P*}(\co_{P_U})$ generated by the image of the canonical morphism 
$u^*(\co_X(R))\otimes_{\co_P}\cI_Y\rightarrow j_{P*}(\co_{P_U})$. 
\end{lem}

We may assume that $X=\Spec(A)$ and $P=\Spec(B)$ are affine, that $R$ is defined in $X$ by 
an equation $t\in A$ and that $Y$ is defined by an ideal of finite type $J_Y$ of $B$. Let $I$ be the ideal
of $B$ generated by $J_Y$ and $t$, and let $P'$ be the blow-up of $P$ along $I$. Then $P^{(R)}$
is the maximal open subscheme of $P'$ where the exceptional divisor $I\co_{P'}$ is generated by $t$. 
Let $a_1,\dots,a_n\in B$ be generators of $J_Y$. We put
\begin{eqnarray*}
C'&=&B[\frac{a_1}{t},\dots,\frac{a_n}{t}]=\frac{B[\xi_1,\dots,\xi_n]}{(a_1-t\xi_1,\dots,a_n-t\xi_n)},\\
C&=&C'/C'_{t{\rm \text -tor}},
\end{eqnarray*}
where $C'_{t{\rm \text -tor}}$ is the ideal of $C'$ of elements annihilated by a power of $t$. 
Then we have $P^{(R)}=\Spec(C)$, which implies the assertion. 

\subsection{}\label{tub3} 
Let $X$ be a scheme, $u\colon P\rightarrow X$ and  $g\colon Q\rightarrow P$ morphisms, 
$i\colon Y\rightarrow P$ and $j\colon Z\rightarrow Q$ closed immersions and $h\colon Z\rightarrow Y$
a morphism such that $g\circ j=i\circ h$; so the diagram 
\[
\xymatrix{
Z\ar[r]^j\ar[d]_{h}&Q\ar[d]^g\\
Y\ar[r]^i&P}
\]
is commutative. Let $R$ be a closed subscheme of $X$ and let
$P^{(R)}$ (resp. $Q^{(R)}$) be the dilatation of $P$ (resp. $Q$) along $Y$ (resp. $Z$) of thickening $R$.
By the functorial property of dilatations (\cite{as4} 2.6), there is a canonical morphism  
\begin{equation}\label{tub3a}
g^{(R)}\colon Q^{(R)}\rightarrow P^{(R)}
\end{equation}
lifting $g$. 

\begin{lem}\label{tub31} 
We keep the assumptions of \eqref{tub3} and assume moreover that 
$Z\simeq Y\times_PQ$. We put $R_Y=R\times_XY$ and $R_Z=R\times_XZ$ and denote by $P'$ (resp. $Q'$) 
the blow-up of $P$ along $R_Y$ (resp. $Q$ along $R_Z$). Then~:

{\rm (i)}\ There exists a unique morphism $g'\colon Q'\rightarrow P'$ lifting $g$; we have 
$Q^{(R)}=g'^{-1}(P^{(R)})$ and $g^{(R)}$ is the restriction of $g'$. 

{\rm (ii)}\  If $g$ is flat, the morphism $Q^{(R)}\rightarrow P^{(R)}\times_PQ$ induced 
by $g^{(R)}$ is an isomorphism.
\end{lem}

Since $R_Z\simeq R_Y\times_PQ$, there exists 
a unique morphism $g'\colon Q'\rightarrow P'$ lifting $g$. We know (\cite{as4} 2.6) that $g'(Q^{(R)})\subset P^{(R)}$
and that $g^{(R)}$ is the restriction of $g'$.
On the other hand, $g'^{-1}(P^{(R)})\subset Q^{(R)}$ by (\cite{as4} 2.7). Therefore, $Q^{(R)}=g'^{-1}(P^{(R)})$,
which proves assertion (i).  Assertion (ii) is an immediate consequence of (i).

\subsection{}\label{tub25} 
Let $X$ be a scheme, $u\colon P\rightarrow X$ a separated morphism, 
$\sigma\colon X\rightarrow P$ a section of $u$, $R$ an effective Cartier divisor on $X$, $U=X-R$ and
$P^{(R)}$ the dilatation of $P$ along $\sigma$ of thickening $R$. Then we have a canonical isomorphism 
\begin{equation}\label{tub25a}  
P^{(R)}\times_XU\simeq P_U. 
\end{equation}
By the universal property of dilatations (\cite{as4} 2.7), 
there exists a unique $X$-morphism 
\begin{equation}\label{tub25c}  
\sigma^{(R)}\colon X\rightarrow P^{(R)}
\end{equation}
lifting $\sigma$.

\begin{lem}\label{tub26} 
We keep the assumptions of \eqref{tub25} and assume moreover that $X$ is locally noetherian 
and that $u$ is smooth. 
Then $P^{(R)}$ is smooth over $X$, and we have a canonical $R$-isomorphism
\begin{equation}\label{tub25b}
P^{(R)}\times_XR\stackrel{\sim}{\rightarrow}\bV(\sigma^*(\Omega^1_{P/X})\otimes_{\co_X}\co_X(R))\times_XR.  
\end{equation} 
\end{lem}

The isomorphism \eqref{tub25b} follows from (\cite{as4} 3.5). 
Since $R$ is a Cartier divisor on $X$, the isomorphisms \eqref{tub25a} and \eqref{tub25b} imply that 
$P^{(R)}$ is flat over $X$ (\cite{egr1} 1.12.9). Since all fibers of $P^{(R)}$ over $X$ are smooth, 
$P^{(R)}$ is smooth over $X$.

\begin{lem}\label{tub6}
Let $X$ be a locally noetherian scheme, $R$ a Cartier divisor on $X$, $u\colon P\rightarrow X$ and
$v\colon Q\rightarrow X$ separated morphisms of finite type and $\sigma\colon X\rightarrow P$ and
$\tau\colon X\rightarrow Q$ sections of $u$ and $v$, respectively. 
We denote by $P^{(R)}$ (resp. $Q^{(R)}$, resp. $(P\times_XQ)^{(R)}$) the dilatation of $P$ (resp. $Q$,
resp. $P\times_XQ$) along $\sigma$ (resp. $\tau$, resp. $(\sigma,\tau)$) of thickening $R$. 
If $P$ or $Q$ is smooth over $X$, then the canonical morphism 
\begin{equation}\label{tub6a}
w\colon (P\times_XQ)^{(R)}\rightarrow P^{(R)}\times_XQ^{(R)}
\end{equation}
is an isomorphism.
\end{lem}

We denote by $\cJ_P$ (resp. $\cJ_Q$, resp. $\cJ_{P\times_XQ}$) the ideal of $\co_P$ (resp. $\co_Q$,
resp. $\co_{P\times_XQ}$) defined by $\sigma$ (resp. $\tau$, resp. $(\sigma,\tau)$). 
Since $\cJ_{P\times_XQ}=\cJ_P\co_{P\times_XQ}+\cJ_Q\co_{P\times_XQ}$, 
it follows from \ref{tub2} that the canonical morphism $w$ \eqref{tub6a} is a closed immersion. 
By construction, $P_U$ (resp. $Q_U$, resp. $P_U\times_UQ_U$) is schematically dense in $P^{(R)}$ (resp. $Q^{(R)}$,
resp. $(P\times_XQ)^{(R)}$). Since $P^{(R)}$ or $Q^{(R)}$ is smooth over $X$ \eqref{tub26}, $P_U\times_UQ_U$
is schematically dense in $P^{(R)}\times_XQ^{(R)}$. Therefore $w$ is an isomorphism.

\subsection{}\label{tub4}
Consider a commutative diagram of morphisms of finite type of locally noetherian schemes
\begin{equation}\label{tub4a}
\xymatrix{
Y\ar[r]^\tau\ar[d]_f&Q\ar[r]^v\ar[d]^g&Y\ar[d]^f\\
X\ar[r]^\sigma&P\ar[r]^u&X}
\end{equation}
such that $u$ and $v$ are smooth and separated, $u\circ \sigma=\id_X$ and $v\circ \tau=\id_Y$.
Let $R$ be a Cartier divisor on $X$  such that $R_Y=R\times_XY$ is a Cartier divisor on $Y$.
We denote by $P^{(R)}$ (resp. $Q^{(R)})$ the dilatation of $P$ (resp. $Q$) along $\sigma$ (resp. $\tau$) 
of thickening $R$ and by 
\begin{equation}\label{tub4b}
g^{(R)}\colon Q^{(R)}\rightarrow P^{(R)}
\end{equation}
the morphism induced by $g$ \eqref{tub3a}. Note that $Q^{(R)}$ is also the dilatation of $Y$ in $Q$ of thickening $R_Y$.
Let $\cN_{X/P}$ and $\cN_{Y/Q}$ be the conormal bundles of $X$ in $P$ and $Y$ in $Q$, respectively.
Then the morphism $g^{(R)}\times_XR\colon Q^{(R)}\times_XR\rightarrow P^{(R)}\times_XR$
can be identified with the morphism 
\begin{equation}\label{tub4c}
\bV(\cN_{Y/Q}\otimes_{\co_Y}\co_Y(R_Y))\times_YR_Y\rightarrow 
\bV(\cN_{X/P}\otimes_{\co_X}\co_X(R))\times_XR
\end{equation}
induced by the canonical morphism $f^*(\cN_{X/P})\rightarrow \cN_{Y/Q}$ (\cite{as4} 3.4).

\begin{lem}\label{tub5}
We keep the assumptions of \eqref{tub4}, and assume moreover that $g$ is smooth and that $f$ is an isomorphism.
Then $g^{(R)}\colon Q^{(R)}\rightarrow P^{(R)}$ is smooth. 
\end{lem}
Observe first that $P^{(R)}$ and $Q^{(R)}$ are smooth over $X$ \eqref{tub26}. We put $U=X-D$. 
Then $g^{(R)}\times_XU=g_U$, which is smooth by assumption. 
On the other hand, $g^{(R)}\times_XR$ is the morphism of vector bundles over $R$ induced 
by the canonical morphism $\cN_{X/P}\rightarrow \cN_{X/Q}$ \eqref{tub4c}. Since 
the latter is locally left invertible, $g^{(R)}\times_XR$ is smooth. 
Then the assertion follows from (\cite{ega4} 17.8.2).

\section{Frames and strict normal crossing pairs}\label{fram}

\subsection{}\label{fram1}
In this article, a {\em monoid} stands for a commutative monoid. 
If $M$ is a monoid, we denote by $M^\gp$ the associated group,
by $M^\times$ the group of units in $M$ and by $\oM$ the orbit space $M/M^\times$ (which is also the quotient 
of $M$ by $M^\times$ in the category of monoids). 
We say that a monoid $M$ is {\em integral} if the canonical homomorphism $M\rightarrow M^\gp$ is injective, 
that $M$ is {\em fine} if it is finitely generated and integral and that $M$ is {\em saturated} if it is integral
and equal to its saturation in $M^\gp$ (i.e., equal to $\{m\in M^\gp; m^n\in M \ {\rm for \ some}\ n\geq 1\}$).
If a monoid $M$ is integral, $\oM$ is integral.
We say that a homomorphism of monoids $u\colon M\rightarrow N$ is {\em strict} if the induced morphism 
$\ou\colon \oM\rightarrow \oN$ is an isomorphism.  
We denote by $\Mon$ the category of monoids and by $\Mon_\fs$ the full subcategory of fine
and saturated monoids (usually called fs-monoids for short).

\subsection{}\label{fram2}
A {\em pre-logarithmic structure} on a scheme $X$ is a pair $(\cM,\alpha)$ 
where $\cM$ is a sheaf of abelian monoids on the étale site of $X$
and $\alpha$ is a homomorphism from $\cM$ to the multiplicative monoid $\co_X$. 
A pre-logarithmic structure $(\cM,\alpha)$ is called a {\em logarithmic structure} if  $\alpha$
induces an isomorphism $\alpha^{-1}(\co_X^\times)\stackrel{\sim}{\rightarrow} \co_X^\times$.
Pre-logarithmic structures on $X$ form naturally a category, containing the full subcategory of 
logarithmic structures on $X$. 
The canonical injection from the category of logarithmic structures on $X$
to the category of pre-logarithmic structures on $X$ has a left adjoint. It associates to 
a pre-logarithmic structure $(\cP,\beta)$ the logarithmic structure $(\cM,\alpha)$, where $\cM$
is defined by the following co-cartesian diagram.
\begin{equation}
\xymatrix{
{\beta^{-1}(\co_X^\times)}\ar[r]\ar[d]& \cP\ar[d]\\
{\co_X^\times}\ar[r]& \cM}
\end{equation}
We say that $(\cM,\alpha)$ is the logarithmic structure {\em associated} to $(\cP,\beta)$. 
 
If $f\colon X\rightarrow Y$ is a morphism of schemes and $(\cM,\alpha)$ is a pre-logarithmic structure on $Y$, 
the sheaf of monoids $f^{-1}(\cM)$ equipped with the composed homomorphism
$f^{-1}(\cM)\rightarrow f^{-1}(\co_Y)\rightarrow \co_X$ is a pre-logarithmic structure on $X$ 
called the inverse image of $(\cM,\alpha)$ and denoted by $f^{-1}(\cM,\alpha)$.

\subsection{}\label{fram3}
A {\em logarithmic scheme} is a triple $(X,\cM_X,\alpha_X)$, usually simply
denoted by $(X,\cM_X)$ or even by $X$, 
consisting of a scheme $X$ and a logarithmic structure $(\cM_X,\alpha_X)$ on $X$.
Logarithmic schemes form a category; we refer to \cite{kato2} for more details. 
If $(X,\cM_X,\alpha_X)$ is a logarithmic scheme, we denote by $\cM^\times_X$ the sheaf of units in $\cM_X$, 
by $\cM^\gp_X$ the sheaf associated to the presheaf $U\mapsto \Gamma(U, \cM_X)^\gp$ 
and by $\ocM_X$ the sheaf associated to the presheaf $U\mapsto \Gamma(U, \cM_X)/\Gamma(U,\cM_X)^\times$
(which is the quotient of $\cM_X$ by $\cM^\times_X$ in the category of sheaves of monoids). 
Observe that $\alpha_X$ identifies $\cM^\times_X$ with $\co_X^\times$. 

We say that a morphism of logarithmic schemes $f\colon (X,\cM_X,\alpha_X)\rightarrow (Y,\cM_Y,\alpha_Y)$
is {\em strict} if $(\cM_X,\alpha_X)$ is the logarithmic structure associated to the pre-logarithmic structure
$f^{-1}(\cM_Y,\alpha_Y)$ on $X$, or equivalently if the canonical morphism $f^{-1}(\ocM_Y)\rightarrow \ocM_X$ is 
an isomorphism.
 
We say that a logarithmic scheme $(X,\cM_X,\alpha_X)$ is {\em integral} 
(resp. {\em fine}, resp. {\em saturated}) if for every $x\in X$, 
there exists an étale neighbourhood $U$ of $x$ in $X$ such that $(\cM_X|U,\alpha_X|U)$ is associated 
to a pre-logarithmic structure $(P_U,\beta)$ on $U$, where $P_U$ is a constant sheaf of monoids 
on $U$ of value an integral (resp. a fine, resp. a saturated) monoid $P$. 
If $(X,\cM_X,\alpha_X)$ is integral (resp. saturated), the monoid $\Gamma(X,\cM_X)$ is integral (resp. saturated). 
If $(X,\cM_X,\alpha_X)$ is integral (resp. fine), for every geometric point $\ox$ of $X$, the monoid 
$\ocM_{X,\ox}$ is integral (resp. fine).

\subsection{}\label{fram4}
Let $M$ be a monoid and $X$ a logarithmic scheme.
We denote by $B[M]$ the scheme $\Spec(\mZ[M])$ equipped with the logarithmic structure induced by 
the pre-logarithmic structure $M\rightarrow \mZ[M]$ (denoted by ${\bf S}[M]$ in \cite{ks1} §4.1), 
and by $M_X$ the constant sheaf of monoids on $X$ of value $M$.
Then the following data are equivalent (and will be identified in what follows)~:

(i) A homomorphism $M\rightarrow \Gamma(X,\cM_X)$; 

(ii) A homomorphism $M_X\rightarrow \cM_X$;

(iii) A morphism of logarithmic schemes $X\rightarrow B[M]$. 

Moreover, the following conditions are equivalent~:

(a) $\cM_X$ is associated to the pre-logarithmic 
structure induced on $M_X$. 

(b) The morphism $X\rightarrow B[M]$ is strict. 

We say then that $M$ is a {\em chart} for $X$. 

\subsection{}\label{fram5}
We denote by $\LS$ the category of fine and saturated logarithmic schemes (usually called fs-logarithmic schemes for short)
and by $\hLS$ the category of presheaves of sets over $\LS$. 
Since the canonical functor $\LS\rightarrow \hLS$ is fully faithful, 
we will identify the objects of $\LS$ with their canonical images in $\hLS$. 
Fibred products are representable in the category $\LS$. 
For morphisms $X\rightarrow S$ and $Y\rightarrow S$ of $\LS$, 
we will denote by $X\timesl_SY$ the fibered product in the category $\LS$
and reserve the notation $X\times_SY$ for the fibered product of the underlying schemes. 
To avoid any risk of confusion, we will usually use the same notation for fibered products in $\hLS$ 
(but not for products as there is no risk of confusion).

\subsection{}\label{fram51}
We have a functor 
\begin{equation}\label{fram51a}
\Mon^\circ_\fs\rightarrow\LS, \ \ \ M\mapsto B[M].
\end{equation} 
Let  $g\colon N\rightarrow M$ and $g'\colon N\rightarrow M'$ be two morphisms of $\Mon_\fs$. 
We denote by $M^\gp\oplus_{N^\gp}M'^\gp$ the cokernel of the homomorphism $g^\gp-g'^\gp\colon N^\gp
\rightarrow M^\gp\oplus M'^\gp$ and by $M\oplus_{N}^\sat M'$ the saturation of the image of the canonical 
homomorphism $M\times M'\rightarrow  M^\gp \oplus_{N^\gp}M'^\gp$ ($M\oplus_{N}^\sat M'$ is the amalgamated 
sum of $g$ and $g'$ in $\Mon_\fs$).  Then we have a canonical isomorphism in $\LS$
\begin{equation}
B[M\oplus_{N}^\sat M'] \stackrel{\sim}{\rightarrow} B[M]\timesl_{B[N]}B[M']. 
\end{equation}
Consider a commutative diagram of $\LS$
\begin{equation}
\xymatrix{
X\ar[r]\ar[d]&Y\ar[d]&X'\ar[l]\ar[d]\\
{B[M]}\ar[r]&{B[N]}&{B[M']}\ar[l]}
\end{equation}
where the vertical arrows are strict and the lower horizontal morphisms are induced by $g$ and $g'$.
Then we have a canonical isomorphism of the underlying schemes
\begin{equation}
X\timesl_YX'\stackrel{\sim}{\rightarrow}(X\times_YX')\times_{B[M\times M']}B[M\oplus_{N}^\sat M']
\end{equation} 
and $X\timesl_YX'$ is strict over $B[M\oplus_{N}^\sat M']$.

\subsection{}\label{fram6}
Following  (\cite{ks1} 4.1.1), we denote by 
\begin{equation}\label{fram6a}
\Mon_\fs^\circ \rightarrow \hLS, \ \ \ M\mapsto [M],
\end{equation}
the functor defined, for a fine and saturated monoid $M$ and $X\in \ob(\LS)$, by 
\begin{equation}\label{fram6b}
[M](X)=\Hom(M,\Gamma(X,\ocM_X)).
\end{equation}
We denote by $\MLS$ the following category. Objects of $\MLS$ are triples $(X,M,u)$,
where $X\in \ob(\LS)$, $M$ is a  fine and saturated monoid and $u\colon X\rightarrow [M]$
is a morphism of $\hLS$. Let $(X,M,u)$ and $(Y,N,v)$ be two objects of $\MLS$. A morphism from 
$(X,M,u)$ to $(Y,N,v)$ is a pair $(f,g)$ made of a morphism $f\colon X\rightarrow Y$ of $\LS$
and a homomorphism of monoids $g\colon N\rightarrow M$ such that the diagram 
\begin{equation}\label{fram6c}
\xymatrix{
X\ar[r]^-(0.5)u\ar[d]_f&{[M]}\ar[d]^{[g]}\\
Y\ar[r]^-(0.5)v&{[N]}}
\end{equation}
is commutative. An object $(X,M,u)$ of $\MLS$ is called a {\em framed logarithmic scheme}
(and $(M,u)$ is called a {\em frame} on $X$) if for every geometric point $\ox$ of $X$, there exists 
an étale neighbourhood $U$ of $\ox$ in $X$ such that the morphism $U\rightarrow [M]$ induced by $u$ 
factors as $U\stackrel{v}{\rightarrow}B[M]\stackrel{w}{\rightarrow} [M]$, where $v$ is a strict morphism 
and $w$ is the canonical morphism (\cite{ks1} 4.1.2); we say also that $u$ is {\em strict}.

\begin{prop}[\cite{ks1} 4.2.1]\label{fram7}
Let $g\colon N\rightarrow M$ be a morphism of $\Mon_\fs$ such that $g^\gp\colon N^\gp\rightarrow M^\gp$
is surjective. Then~: 

{\rm (i)}\ The morphism $[g]\colon [M]\rightarrow [N]$ is representable, log étale and affine, i.e., 
for every $X\in \ob(\LS)$ and every morphism $u\colon X\rightarrow [N]$, 
the fibre product $X\timesl_{[N]}[M]$ is representable
by an object of $\LS$ which is log-étale and affine over $X$. 

{\rm (ii)}\ Let $\tM$ be the inverse image of $M$ by $g^\gp\colon N^\gp\rightarrow M^\gp$. 
Then the canonical morphism $[M]\rightarrow [\tM]$ is an isomorphism, 
and for every morphism $u\colon X\rightarrow B[N]$, the canonical morphism 
\begin{equation}\label{fram7a}
X\timesl_{B[N]}B[\tM] \rightarrow X\timesl_{[N]}[\tM]
\end{equation}
is an isomorphism. 
\end{prop}

\subsection{}\label{fram8}
Let $X,Y,S$ be objects of $\LS$, $M$ a finitely generated and saturated monoid and
$X\rightarrow S\times [M]$ and $Y\rightarrow S\times [M]$ two morphisms of $\hLS$.  
We will denote such a diagram by $X,Y\rightrightarrows S\times [M]$ and its projective limit by $X\timesl_{S\times [M]}Y$.  
Let $\mu\colon M\times M\rightarrow M$ be the multiplication and $i_1,i_2\colon M\rightarrow M\times M$
the homomorphisms defined by $i_1(m)=(m,1)$ and $i_2(m)=(1,m)$. Since the diagram
\begin{equation}\label{fram8b}
\xymatrix{
M\ar@<0.5ex>[r]^-(0.5){i_1}\ar@<-0.5ex>[r]_-(0.5){i_2}&{M\times M}\ar[r]^-(0.5)\mu&M}
\end{equation}
is co-exact in the category of monoids (i.e., $\mu$ is the cokernel of the pair of morphisms 
$i_1$ and $i_2$), the diagram 
\begin{equation}\label{fram8c}
\xymatrix{
{X\timesl_{S\times [M]}Y}\ar[r]\ar[d]&{[M]}\ar[d]^{[\mu]}\\
{X\timesl_SY}\ar[r]^-(0.5){u\times v}&{[M\times M]}}
\end{equation}
where $u\times v$ is induced by $u$ and $v$, is cartesian in $\hLS$. 

\begin{cor}[\cite{ks1} 4.2.3]\label{fram9}
For every diagram $X,Y\rightrightarrows S\times [M]$ of $\hLS$, $X\timesl_{S\times [M]}Y$  is representable
by an object of $\LS$, which is log-étale and affine over $X\timesl_SY$. 
\end{cor}

\begin{prop}[\cite{ks1} 4.2.5] \label{fram11}
Let $g\colon N\rightarrow M$ be a homomorphism of finitely generated and saturated monoids, 
$X$, $Y$ and $S$ objects of $\LS$ and
\begin{equation}
\xymatrix{
X\ar[r]\ar[d]&S\ar[d]&Y\ar[l]\ar[d]\\
{[M]}\ar[r]^{[g]}&{[N]}&{[M]}\ar[l]_{[g]}}
\end{equation}
a commutative diagram. 
Assume that $X\rightarrow [M]$ and $S\rightarrow [N]$ are strict. Then the canonical projection $X\timesl_{S\times [M]}Y\rightarrow Y$ is strict.
\end{prop}

\begin{cor} \label{fram12}
Under the assumptions of \eqref{fram11}, if moreover $X\rightarrow S$ is log-smooth, 
then the canonical projection $X\timesl_{S\times [M]}Y\rightarrow Y$ is strict and smooth. 
\end{cor}
It follows from \ref{fram9} and \ref{fram11} as $X\timesl_{S\times [M]}Y\rightarrow Y$ 
is composed of $X\timesl_{S\times [M]}Y\rightarrow X\timesl_SY\rightarrow Y$.

\subsection{}\label{fram13}
Let $(X,M,u)$ be an object of $\MLS$ and $Y\rightarrow X\rightarrow S$ two morphisms
of $\LS$. Then we can form a diagram $X,Y\rightrightarrows S\times [M]$,
and the canonical morphism 
\begin{equation}\label{fram13a}
X\timesl_{S\times [M]}Y\rightarrow (X\timesl_{S\times [M]}X)\timesl_XY
\end{equation}
is an isomorphism. 
In particular, if the morphism $X\rightarrow S$ can be extended to a morphism of framed logarithmic 
schemes $(X,M,u)\rightarrow (S,N,v)$, then the canonical morphism of the underlying schemes  
\begin{equation}\label{fram13c}
X\timesl_{S\times [M]}Y\rightarrow (X\timesl_{S\times [M]}X)\times_XY
\end{equation}
is an isomorphism, and $X\timesl_{S\times [M]}Y$ is strict over $Y$ \eqref{fram11}.

\subsection{}\label{fram10}
Let $(X,M,u)$ and $(Y,N,v)$ be two objects of $\MLS$, $X\rightarrow S$ and $Y\rightarrow S$ two morphisms
of $\LS$ and $g\colon M\rightarrow N$ a homomorphism; so we can form a diagram $X,Y\rightrightarrows S\times [M]$. 
We denote by $\theta\colon M\times N\rightarrow N$ the homomorphism defined by $\theta(m,n)=g(m)\cdot n$. 
Since the diagram 
\begin{equation}\label{fram10a}
\xymatrix{
{M\times M}\ar[r]^-(0.5)\mu\ar[d]_{\id\times g}&M\ar[d]^g\\
{M\times N}\ar[r]^-(0.5)\theta&N}
\end{equation} 
is co-cartesian in the category of monoids, the diagram
\begin{equation}\label{fram10b}
\xymatrix{
{X\timesl_{S\times [M]}Y}\ar[r]^-(0.5)w\ar[d]&{[N]}\ar[d]^{[\theta]}\\
{X\timesl_SY}\ar[r]^-(0.5){u\times v}&{[M\times N]}}
\end{equation}
where $w$ is the composed morphism $X\timesl_{S\times [M]}Y\rightarrow Y\stackrel{v}{\rightarrow}[N]$,
is cartesian in $\hLS$. 

\subsection{}\label{fram14}
Let $(f,g)\colon (Y,N,v)\rightarrow (X,M,u)$ be a morphism of $\MLS$ and $h\colon X\rightarrow S$ a morphism of $\LS$; 
so we have the following commutative diagram.
\begin{equation}\label{fram14a}
\xymatrix{
Y\ar[r]^v\ar[d]_f&{[N]}\ar[d]^{[g]}\\
X\ar[d]_h\ar[r]^u&{[M]}\\
S&}
\end{equation}
Then the canonical projection $\pr_2\colon X\timesl_{S\times [M]}Y\rightarrow Y$ induces an isomorphism
\begin{equation}\label{fram14b}
Y\timesl_{X\times [N]}(X\timesl_{S\times [M]}Y) \stackrel{\sim}{\rightarrow} Y\timesl_{S\times [N]}Y.
\end{equation}
Indeed, the commutative diagram 
\begin{equation}
\xymatrix{
{Y\timesl_{S\times[N]}Y}\ar[r]\ar[d]_{\pr_1}&{X\timesl_{S\times[M]}Y}\ar[d]^{\id\times v}\\
Y\ar[r]^{f\times v}&{X\times [N]}}
\end{equation}
defines the inverse.

\subsection{}\label{fram15}
Let $X$ be a regular noetherian scheme, $D$ a normal crossing divisor on $X$, 
$U=X-D$ and $j\colon U\rightarrow X$ the canonical injection.  There is a canonical fine and saturated 
logarithmic structure $(\cM_D,\alpha_D)$ on $X$ defined by $\cM_D=\co_X\times_{j_*(\co_U)}j_*(\co_U^\times)$.
We denote $(X,\cM_D,\alpha_D)$ by $X_{\log D}$.
The sheaf $\ocM_D$ is canonically isomorphic to the sheaf $\uGamma_D(\cDiv_X^+)$ 
of effective Cartier divisors on $X$ with support in $D$. Assume that $D$ has simple normal 
crossings, and let $D_1\dots,D_m$ be the irreducible components of $D$. We denote by 
$M_D$ the free abelian monoid generated by $D_1,\dots,D_m$. Then the canonical 
morphism $u_D\colon X_{\log D}\rightarrow [M_D]$ defines a frame on $X$. 

\subsection{}\label{log1}
A {\em  strict normal crossing pair} over $k$ (or an {\em snc-pair} over $k$ for short)
stands for a pair $(X,D)$ where $X$ is a smooth $k$-scheme and $D$ 
is a simple normal crossing divisor on $X$. 
Let $(X,D)$ and $(Y,E)$ be two snc-pairs over $k$. 
A morphism $f\colon  (Y,E)\rightarrow (X,D)$ is a $k$-morphism 
$f\colon Y\rightarrow X$ such that the support of $f^{-1}(D)$ is contained in $E$. We denote by $\SNCP_k$ 
the category of snc-pairs over $k$. We have a canonical functor 
\begin{equation}
\SNCP_k\rightarrow \MLS,\ \ \ (X,D)\mapsto (X_{\log D}, M_D, u_D)
\end{equation}
where $X_{\log D}$, $M_D$ and $u_D$ are defined in \eqref{fram15}.
We say that a morphism $f\colon (Y,E)\rightarrow (X,D)$ of snc-pairs over $k$  
is {\em log-smooth} (resp. {\em log-étale}) if the associated morphism $Y_{\log E}\rightarrow X_{\log D}$ is log-smooth
(resp. log-étale). 

\begin{lem}\label{rc1}
Let $f\colon (X',D')\rightarrow (X,D)$ be a log-smooth morphism of snc-pairs over $k$ such that the morphism
of the underlying schemes $X'\rightarrow X$ is flat, $x\in f(X')\subset X$ and $D_1,\dots,D_n$ the irreducible components 
of $D$ containing $x$. 
Then there exists $x'\in X'$ contained in exactly $n$ irreducible components $D'_1,\dots,D'_n$ of $D'$
such that $f(x')=x$ and $D'_i$ dominates $D_i$ for all $1\leq i\leq n$.
\end{lem}

We may shrink $X$, so we may assume that the irreducible components of $D$ are $D_1,\dots,D_n$. 
Let $y'\in X'$ such that $f(y')=x$. Since $f$ is flat, for each $1\leq i\leq n$, there exists an irreducible component 
$D'_i$ of $D'$ containing $y'$ and dominating $D_i$. We put $Y=\bigcap_{1\leq i\leq n}D_i$ and 
$Y'=\bigcap_{1\leq i\leq n}D'_i$ and denote by $g\colon Y'\rightarrow Y$ the restriction of $f$. 
We equip $Y'$ with the strictly normal crossing divisor $E'$ defined by the irreducible components of $D'$
different from  $D'_1,\dots,D'_n$. We have the following canonical commutative diagram of $\co_{Y'}$-modules 
with exact lines.
\begin{equation}
\xymatrix{
0\ar[r]&{g^*(\Omega^1_{Y/k})}\ar[r]\ar[d]&{g^*(\Omega^1_{X/k}(\log D)\otimes_{\co_X}\co_Y)}\ar[r]^-(0.5){g^*(\res)}\ar[d]&
{\co_{Y'}^n}\ar[r]\ar@{=}[d]&0\\
0\ar[r]&{\Omega^1_{Y'/k}(\log E')}\ar[r]&{\Omega^1_{X'/k}(\log D')\otimes_{\co_{X'}}\co_{Y'}}\ar[r]^-(0.5){\res}&
{\co_{Y'}^n}\ar[r]&0}
\end{equation}
Therefore, the morphism $g^*(\Omega^1_{Y/k})\rightarrow \Omega^1_{Y'/k}(\log E')$ is injective
and its cokernel is locally free. Then the morphism of snc-pairs $(Y',E')\rightarrow (Y,\emptyset)$ induced by $g$ 
is log-smooth (\cite{kato2} 3.12), which implies that $g$ is smooth (in the usual sense) 
and that $E'$ is a strict normal crossing divisor on $Y'$ relatively to $Y$ (\cite{kato2} 3.5). 
In particular, $Y'_x$ is smooth over $k(x)$ and $E'_x$ is a divisor on $Y'_x$.
Hence, there exists $x'\in Y'_x-E'_x$. 

\begin{lem}\label{clean72}
Let $f\colon (X^\dagger,D^\dagger)\rightarrow (X,D)$ and $g\colon (X',D')\rightarrow (X,D)$ be two morphisms 
of snc-pairs over $k$ and $x^\dagger\in X^\dagger$. Assume that $g$ is log-smooth, that the morphism $X'\rightarrow X$ is flat
and that $f(x^\dagger)\in g(X')$. Then there exists a commutative diagram of snc-pairs over $k$ 
\begin{equation}\label{clean72a}
\xymatrix{
{(X^\ddagger,D^\ddagger)}\ar[r]^{f'}\ar[d]_{g^\dagger}&{(X',D')}\ar[d]^g\\
{(X^\dagger,D^\dagger)}\ar[r]^f&{(X,D)}}
\end{equation}
such that $g^\dagger$ is log-smooth, that the morphism $X^\ddagger\rightarrow X^\dagger$ 
is flat and that $x^\dagger\in g^\dagger(X^\ddagger)$. Moreover, if $f$ is 
log-smooth, then we can choose $(X^\ddagger,D^\ddagger)$ such that $f'$ is log-smooth.
\end{lem}
We denote by $\cM,\cM'$ and  $\cM^\dagger$ the sheaves of monoids over $X,X'$ and $X^\dagger$ 
defined by the divisors $D,D'$ and $D^\dagger$ respectively; we use the notation of § \ref{fram}. 
Let $\ox^\dagger$ be a geometric point of $X^\dagger$ above $x^\dagger$, $\ox=f(\ox^\dagger)$
and $\ox'$ a geometric point of $X'$ above $\ox$. We put $M=\ocM_{\ox}$, $M'=\ocM'_{\ox'}$ and
$M^\dagger=\ocM^\dagger_{\ox^\dagger}$. By \ref{rc1}, we may assume that the following 
condition is satisfied~: 

(C$_1$)\ The canonical homomorphism $u\colon M\rightarrow M'$ induces an isomorphism  
$M^\gp\otimes_{\mZ}\mQ \stackrel{\sim}{\rightarrow} M'^\gp\otimes_{\mZ}\mQ$ 
and $M$ is the inverse image of $M'$ by the morphism $u^\gp\colon M^\gp\rightarrow M'^\gp$. 

Replacing $X'$ by an étale neighborhood of $\ox'$, we may assume that the following 
condition is satisfied~: 

(C$_2$)\ The divisors $D'$ and $g^*(D)$ have the same support. 

Let $N$ be an integer annihilating the cokernel of $u^\gp$. We put $Q=(M^{\dagger})^\gp\times M^\dagger$
and denote by $q\colon M^\dagger\rightarrow Q$ the morphism defined by $t\mapsto (t,t^N)$. 
Replacing $X^\dagger$ by an étale neighborhood of $\ox^\dagger$, we may assume that there exists 
a chart $X^\dagger\rightarrow B[M^\dagger]$. We put $X_1=X^\dagger\timesl_{B[M^\dagger]}B[Q]$.
Then $X_1\rightarrow X^\dagger$ is log-smooth and the underlying morphism of schemes is faithfully flat.
Hence, by replacing $X^\dagger$ by $X_1$, we may further assume that the following condition is satisfied~:

(C$_3$)\ The canonical homomorphism $u^\dagger\colon M\rightarrow M^\dagger$ factors as 
$M\stackrel{u}{\rightarrow}M'\stackrel{v}{\rightarrow}M^\dagger$. 
 
Under assumptions (C$_1$), (C$_2$) and (C$_3$), we put $X^\ddagger=X^\dagger\timesl_XX'$. 
The saturation $M^\dagger\oplus_M^\sat M'$ of the amalgamated sum of $M^\dagger$ and $M'$ over $M$,
is equal to $M^\dagger$ \eqref{fram51}.
Hence, the canonical projection $g^\dagger\colon X^\ddagger\rightarrow X^\dagger$ is strict. 
Since $g^\dagger$ is log-smooth, the underlying morphism of schemes is smooth. 
Let $D^\ddagger$ be the inverse image of $D^\dagger$ by $g^\dagger$. Then the logarithmic structure on $X^\ddagger$
is induced by $D^\ddagger$. On the other hand, $X^\ddagger$ is the normalization of $X^\dagger\times_XX'$.
Hence $X^\ddagger\rightarrow X^\dagger\times_XX'$ is surjective. The first assertion is proved, and the second one is 
obvious from the definition of $X^\ddagger$.

\subsection{}\label{fram16}
Let $(X,D)$ and $(Y,E)$ be two snc-pairs over $k$ and $g\colon M_D\rightarrow M_E$ a homomorphism
of monoids \eqref{log1}; so we can form the diagram $X_{\log D},Y_{\log E}\rightrightarrows \Spec(k)\times [M_D]$.
We call the {\em $g$-framed product} of $(X,D)$ and $(Y,E)$ over $k$ and denote by $X\Asterisk_{k,g}Y$ the logarithmic
scheme
\begin{equation}\label{fram16a}
X\Asterisk_{k,g}Y=X_{\log D}\timesl_{\Spec(k)\times [M_D]}Y_{\log E}.
\end{equation} 
We know that the canonical morphism $X\Asterisk_{k,g}Y\rightarrow X_{\log D}\timesl_kY_{\log E}$ 
is log-étale \eqref{fram9} and that the second projection $X\Asterisk_{k,g}Y\rightarrow Y$ is strict and smooth \eqref{fram12}.
Observe that $X_{\log D}\timesl_kY_{\log E}$ is the logarithmic scheme associated to the snc-pair 
$(X\times_kY,\pr_1^*(D)+\pr_2^*(E))$. 

\subsection{}\label{fram17}
Let $f\colon (Y,E)\rightarrow (X,D)$ be a morphism of snc-pairs over $k$. Then $f$ induces a homomorphism
of monoids $g\colon M_D\rightarrow M_E$. We call the {\em $f$-framed product} of $(X,D)$ and $(Y,E)$ over $k$
and denote by $X\Asterisk_{k,f}Y$ the logarithmic scheme $X\Asterisk_{k,g}Y$ \eqref{fram16a};
we omit $f$ from the terminology and the notation if there is no risk of confusion.  In particular, 
we call $X\Asterisk_{k,\id}X$ the {\em framed self-product} of $(X,D)$ over $k$ and denote it simply by 
$X\Asterisk_{k}X$.

There is a canonical morphism 
\begin{equation}
\gamma_f\colon Y_{\log E}\rightarrow X\Asterisk_{k,f}Y
\end{equation}
called the {\em framed graph} of $f$. The framed graph of the identity $\delta\colon X\rightarrow X\Asterisk_kX$ is called the {\em framed diagonal} of $(X,D)$.

The formation of $X\Asterisk_{k,f}Y$ is functorial in $f$. In particular, we have canonical morphisms 
\begin{equation}\label{fram17a}
Y\Asterisk_kY\stackrel{f_1}{\longrightarrow} X\Asterisk_kY \stackrel{f_2}{\longrightarrow}  X\Asterisk_kX.
\end{equation}
We put $\tf=f_2\circ f_1$. By \eqref{fram13c}, the canonical morphism of the underlying schemes
\begin{equation}\label{fram17b}
X\Asterisk_kY\rightarrow (X\Asterisk_kX)\times_XY
\end{equation}
is an isomorphism.

\begin{prop}\label{log6}
Let $(X,D)$ be an snc-pair over $k$ and $D_1,\dots,D_m$ the irreducible components of $D$. 
For $1\leq i\leq m$, let $\cI_i$ be the ideal of the closed subscheme $D_i\times_kD_i$ of $X\times_kX$, 
and let $(X\times_kX)'$ be the blow-up of $X\times_kX$ along the ideal $\prod_{1\leq i\leq m}\cI_i$. 
Then $X\Asterisk_kX$ is canonically isomorphic to 
the open subscheme $Z$ of $(X\times_kX)'$ complementary to the strict transforms
of $D\times_kX$ and $X\times_kD$. 
\end{prop}
Let $V=\Spec(A)$ and $W=\Spec(B)$ be affine open subscheme of $X$ such that for all $1\leq i\leq m$, 
$D_i|V$ (resp. $D_i|W$) is defined by an equation $t_i\in A$ (resp. $s_i\in B$). 
The inverse image of $V\times_kW$ in $Z$ is the affine scheme of the ring 
\begin{equation}\label{log6a}
\frac{A\otimes_kB[u_1^{\pm1},\dots,u_m^{\pm1}]}{(t_i\otimes 1-u_i\cdot 1\otimes s_i \  (1\leq i\leq m))}.
\end{equation} 
Hence, $Z$ equipped with the exceptional divisor $E$ is an snc-pair over $k$. 
By construction, for all $1\leq i\leq m$,
the ideals $\pr_1^*(\co_X(-D_i))\co_Z$ and $\pr_2^*(\co_X(-D_i))\co_Z$ are equal and invertible.
Therefore, the canonical morphism $Z\rightarrow X\times_kX$ lifts to a morphism 
$Z\rightarrow X\Asterisk_kX$, which is an isomorphism as it can be easily checked from the local description \eqref{log6a}. 

\begin{remas}\label{log3}
We keep the assumptions of \eqref{log6}. 

(i) The universal property of $X\Asterisk_kX$ can be restated without logarithmic geometry as follows. 
Let $\varphi\colon Y\rightarrow X\times_kX$ be a morphism such that for all $1\leq i\leq m$,
the ideals $\pr_1^*(\co_X(-D_i))\co_Y$ and $\pr_2^*(\co_X(-D_i))\co_Y$ are equal and invertible.
Then there exists a unique morphism $\psi \colon Y\rightarrow X\Asterisk_kX$ lifting $\varphi$. 

(ii) Let $U=X-D$. Then, with the conventions of \ref{not2}, we have canonical isomorphisms
\begin{equation}\label{log3a}
U\times_X(X\Asterisk_kX)\simeq U\times_kU\simeq (X\Asterisk_kX)\times_XU.
\end{equation}

(iii) The canonical projections $X\Asterisk_kX \rightarrow X$ will be denoted also by $\pr_1$ and $\pr_2 $; they
are smooth. The framed diagonal $\delta\colon X\rightarrow X\Asterisk_kX$ is a regular closed immersion with conormal
bundle canonically isomorphic to the sheaf of logarithmic differentials $\Omega^1_{X}(\log D)$ (\cite{ks1} 4.2.8). 
\end{remas}

\begin{prop}\label{fram18}
Let $f\colon (Y,E)\rightarrow (X,D)$ be a log-smooth (resp. log-étale) morphism of snc-pairs over $k$. 
Then the canonical morphism $Y\Asterisk_kY\rightarrow X\Asterisk_kY$ is smooth (resp. étale). 
\end{prop}

By \eqref{fram14b}, the canonical morphism 
\begin{equation}\label{fram18a}
Y\Asterisk_kY\rightarrow Y_{\log E}\timesl_{X_{\log D}\times [M_E]}(X\Asterisk_kY)
\end{equation}
is an isomorphism. Therefore $Y\Asterisk_kY$ is log-étale over $Y_{\log E}\timesl_{X_{\log D}}(X\Asterisk_kY)$ 
\eqref{fram9}. Hence $Y\Asterisk_kY\rightarrow X\Asterisk_kY$ is log-smooth (resp. log-étale). 
Since the second projections $Y\Asterisk_kY\rightarrow Y$ and $X\Asterisk_kY\rightarrow Y$
are strict, $Y\Asterisk_kY\rightarrow X\Asterisk_kY$ is smooth (resp. étale).

\subsection{}\label{log2}
Let  $(X,D)$ be a snc-pair over $k$ and $D_1,\dots,D_m$ the irreducible components of $D$.
A {\em rational divisor on $X$ with support in $D$} is an element $R=\sum_{i=1}^mr_iD_i$ 
of the $\mQ$-vector space generated by $D_1,\dots,D_m$. 
We say that $R$ is {\em effective} if $r_i\geq 0$ for all $i$, and that $R$ has {\em integral coefficients} if $r_i$ 
is integral for all $i$.
We call {\em generic points} of $R$ the generic points of the $D_i$'s such that $r_i\not=0$. 
For every integer $n\geq 0$, we denote by $\lfloor nR\rfloor$ the  
divisor $\sum_{i=1}^m\lfloor nr_i\rfloor D_i$ on $X$, where $\lfloor nr_i\rfloor$ is the integral part of $nr_i$.
If $R$ and $R'$ are two rational divisors on $X$ with support in $D$, we say that $R'$ is {\em bigger} than $R$
and use the notation $R'\geq R$ if $R'-R$ is effective. 
If $f\colon (Y,E)\rightarrow (X,D)$ is a morphism of snc-pairs over $k$ and
$R$ is an effective rational divisor on $X$ with support in $D$, we can define the pull-back $f^*(R)$
as a rational divisor on $Y$ with support in $E$.

\subsection{} \label{log4}
Let  $(X,D)$ be a snc-pair over $k$, $R$ an effective rational divisor on $X$ with support in $D$,
$u\colon P\rightarrow X$ a smooth separated morphism of finite type and
$s\colon X\rightarrow P$ a section of $u$. 
We put $U=X-D$, and denote by $j\colon U\rightarrow X$ and $j_P\colon P_U\rightarrow P$ 
the canonical injections and by $\cI_X$ the ideal of $X$ in $P$. 
We call {\em dilatation of $P$ along $s$ of thickening $R$} and 
denote by $P^{(R)}$ the affine scheme over $P$ 
defined by the quasi-coherent sub-$\co_P$-algebra of $j_{P*}(\co_{P_U})$
\begin{equation}\label{log4a}
\sum_{n\geq 0}u^*(\co_X(\lfloor nR\rfloor))\cdot \cI_X^{n}.
\end{equation}
This notion extends the one introduced in \ref{tub1} if $R$ has integral coefficients (cf. \ref{tub2}). 
We have a canonical isomorphism 
\begin{equation}\label{log4b}
P^{(R)}\times_XU\simeq P_U.
\end{equation}
The image of the algebra \eqref{log4a} by the surjective homomorphism $j_{P*}(\co_{P_U})\rightarrow s_*j_*(\co_U)$
is canonically isomorphic to $s_*(\co_X)$. Hence we have a canonical section 
\begin{equation}\label{log4c}
s^{(R)}\colon X\rightarrow P^{(R)}
\end{equation}
lifting $s$. 

Let $R'$ be another rational divisor on $X$ with support in $D$ such that $R'\geq R$. 
Then for every $n\geq 0$, there is a canonical injection $\co_X(\lfloor nR\rfloor)\rightarrow \co_X(\lfloor nR'\rfloor)$.
We deduce a canonical $P$-morphism 
\begin{equation}\label{log4d}
P^{(R')}\rightarrow P^{(R)}
\end{equation}
that fits into the following commutative diagram. 
\begin{equation}\label{log4e}
\xymatrix{
&{P^{(R')}}\ar[d]&\\
X\ar[ru]^{s^{(R')}}\ar[r]^{s^{(R)}}&{P^{(R)}}&{P_U}\ar[lu]\ar[l]}
\end{equation}

\begin{prop}\label{log5}
We keep the assumptions of \eqref{log4}, moreover, let $\lambda$ be an integer $\geq 1$ such that $\lambda R$
has integral coefficients, $X^{(\lambda)}$ the $\lambda$-th infinitesimal neighbourhood of $s\colon X\rightarrow P$,
and $P^\dagger$ the dilatation of $P$ along $X^{(\lambda)}$ of thickening $\lambda R$ in the sense of \eqref{tub1}.  
Then $P^{(R)}$ is canonically isomorphic to the integral closure of $P^\dagger$ in $P_U$.
\end{prop}

Recall \eqref{tub2} that $P^\dagger$ is the affine scheme over $P$ defined by the quasi-coherent sub-$\co_P$-algebra 
of $j_{P*}(\co_{P_U})$ generated by the image of the canonical morphism 
$u^*(\co_{X}(\lambda R))\cdot \cI_X^{\lambda}\rightarrow j_{P*}(\co_{P_U})$. 
Therefore, there is a canonical integral $P$-morphism
$P^{(R)}\rightarrow P^\dagger$ extending the identity of $P_U$. 
To prove the proposition, it is enough to show that $P^{(R)}$ is normal. 
We put $R'=\lambda R$. 
The canonical morphism $\rho\colon P^{(R)}\rightarrow P$ is an isomorphism above $u^{-1}(X-R')$, and 
$\rho(P^{(R)}\times_XR')\subset s(R')$. Then $P^{(R)}$ is normal at all points above $P-s(R')$. 
Hence, it is enough to prove that $P^{(R)}$ is normal at $s(x)$ for $x\in R'$. 

By the Jacobian criterion of smoothness,  
there exists an open neighborhood $V$ of $s(x)$ in $P$ and sections 
$g_1,\dots,g_m\in \Gamma(V,\cI_X)$ such that $g_1,\dots,g_m$ generate $\cI_X$ at $s(x)$ and 
$dg_1,\dots,dg_m$ generate $\Omega^1_{P/X}$ at $s(x)$.   
Let $g\colon V\rightarrow \mA_X^m$ be the $X$-morphism defined by $g_1,\dots,g_m$. 
Then $g$ is étale at $s(x)$ by (\cite{ega4} 17.11.1). After shrinking $V$, we may assume that $s(X)\cap V$ coincides with 
the inverse image by $g$ of the zero section of $\mA^m_X$. Hence, we are reduced to the case where $P=\mA^m_X$ 
and $s$ is the zero section. By shrinking $X$, we may further assume that there exists a smooth morphism 
$X\rightarrow \mA^d_k$ such that $D$ is the pull-back of the union of the coordinates hyperplanes of $\mA^d_k$. 
Hence, we are further reduced to the case where $X=\mA^d_k$ and $D$ is the union of the coordinates hyperplanes. 

We denote by $M$ the free monoid $\mN^{d}$ with basis $T_1,\dots,T_{d}$, 
and by $N$ the free monoid $\mN^m$ with basis $T_{d+1},\dots,T_{d+m}$. 
We identify $X$ with $\Spec(k[M])$ and let $f=\prod_{i=1}^dT_i^{a_i}$ be 
an equation defining the divisor $\lambda R$ on $X$. 
Let $H$ be the submonoid of $M^\gp\times N$ generated by $T_i$ for $1\leq i\leq d$ and 
$T_j^\lambda/f$ for $d+1\leq j\leq d+m$. 
We denote by $H^\sat$ the saturation of $H$ and by $\sigma\colon N\rightarrow \mN$ the homomorphism sending 
$T_{d+1},\dots,T_{d+m}$ to $1$. Then $H^\sat$ is the submonoid of $M^\gp\times N$ of elements 
$(\alpha,\beta)$ such that if we write $\alpha=(\alpha_1,\dots,\alpha_d)\in \mZ^d=M^\gp$, we have for all $1\leq i\leq d$, 
\begin{equation}
\lambda \alpha_i+a_i\sigma(\beta)\geq 0.
\end{equation}
Therefore, we have 
\begin{equation}
H^\sat=\coprod_{n \geq 0} \{ (\alpha,\beta)\in M^\gp\times N \ ; \ \alpha \geq \lfloor n R \rfloor \ {\rm and}\ 
\sigma(\beta)=n\},
\end{equation}
where $R$ is considered as an element of $M^\gp\otimes_{\mZ}\mQ$.
Hence, we have a canonical isomorphism $P^{(R)}\simeq \Spec(k[H^\sat])$, and $P^{(R)}$ is normal
as $H^\sat$ is saturated.

\subsection{}\label{log51}
Let $f\colon (Y,E)\rightarrow (X,D)$ be a morphism of snc-pairs over $k$, $U=X-D$, $V=Y-E$, 
$R$ a rational effective divisor on $X$ with support in $D$ and $R_Y=f^*(R)$. 
Let $u\colon P\rightarrow X$ and $v\colon Q\rightarrow Y$ be smooth separated morphisms of finite type, 
$s\colon X\rightarrow P$ a section of $u$, $t\colon Y\rightarrow Q$ a section of $v$
and $g\colon Q\rightarrow P$ a morphism such that the diagram 
\begin{equation}\label{log51a}
\xymatrix{
Y\ar[r]^t\ar[d]_f&Q\ar[d]^g\ar[r]^v&Y\ar[d]^f\\
X\ar[r]^s&P\ar[r]^u&X}
\end{equation}
is commutative. We denote by $P^{(R)}$ (resp. $Q^{(R_Y)}$) the dilatation of 
$P$ (resp. $Q$) along $s$ (resp. $t$) of thickening $R$ (resp. $R_Y$) and by $s^{(R)}\colon X\rightarrow P^{(R)}$
(resp. $t^{(R_Y)}\colon Y\rightarrow Q^{(R_Y)}$) the canonical lifting of $s$ (resp. $t$). 
Let $\cI_X$ be the ideal of $X$ in $P$, $\cI_Y$ the ideal of $Y$ in $Q$ and
$j_P\colon P_U\rightarrow P$ and $j_Q\colon Q_V\rightarrow Q$ the canonical injections. 
The morphism $g$ induces a homomorphism of $\co_Q$-algebras
\begin{equation}\label{log51b}
g^*j_{P*}(\co_{P_U})\rightarrow j_{Q*}(\co_{Q_V}).
\end{equation}
We have $g^*(\cI_X)\co_Q\subset \cI_Y$ and 
$f^*(\lfloor nR\rfloor)\leq \lfloor nR_Y\rfloor$ for every $n\geq 0$. Therefore, \eqref{log51b} 
induces a homomorphism of $\co_Q$-algebras $g^*(\co_{P^{(R)}})\rightarrow \co_{Q^{(R_Y)}}$, and hence a morphism 
\begin{equation}\label{log51c}
g^{(R)}\colon Q^{(R_Y)}\rightarrow P^{(R)}
\end{equation}
lifting $g$. We clearly have 
\begin{equation}\label{log51d}
g^{(R)}\circ t^{(R_Y)}=s^{(R)}\circ f.
\end{equation}
If $R$ has integral coefficients, $g^{(R)}$ is the morphism defined in \eqref{tub4b}.  

\begin{prop}\label{log52}
We keep the assumptions of \eqref{log51} and assume moreover that $Q=P\times_XY$ 
and $g$ and $v$ are the canonical projections. Then $Q^{(R_Y)}$ is the integral closure
of $P^{(R)}\times_PQ$ in $Q_V$. 
\end{prop}
Let $n$ be an integer $\geq 1$ such that $nR$ has integral coefficients, 
$X^{(n)}$ (resp. $Y^{(n)}$) the $n$-th infinitesimal neighbourhood of $s$ (resp. $t$), and
$P^\dagger$ (resp. $Q^\dagger$) the dilatation of $P$ along $X^{(n)}$ (resp. $Q$ along $Y^{(n)}$)
of thickening $nR$ in the sense of \eqref{tub1}. We denote by $F$ (resp. $G$)
the inverse image of $nR$ in $X^{(n)}$  (resp. $Y^{(n)}$) and by $P'$ (resp. $Q'$) the blow-up of $P$
(resp. $Q$) along $F$ (resp. $G$). Since $Y^{(n)}\simeq X^{(n)}\times_PQ$, $g$ lifts to a morphism 
$g'\colon Q'\rightarrow P'$ and we have $Q^{\dagger}=g'^{-1}(P^{\dagger})$ \eqref{tub31}.
Let $g^\dagger\colon Q^\dagger\rightarrow P^\dagger$ be the restriction of $g'$. 
By \ref{log5}, $P^{(R)}$ (resp. $Q^{(R_Y)}$) is the integral closure of $P^\dagger$ (resp. $Q^\dagger$)
in $P_U$ (resp. $Q_V$). It is clear that the canonical morphism $g^{(R)}\colon Q^{(R_Y)}\rightarrow P^{(R)}$ 
\eqref{log51c} is induced by $g^\dagger$. 
The morphism $Q'\rightarrow P'\times_PQ$ induced by $g'$ is a closed immersion (\cite{as4} 2.5).
Therefore, the morphism  $Q^\dagger\rightarrow P^\dagger\times_PQ$ induced by $g^\dagger$ is a closed immersion,
and the integral closures in $Q_V$ of $Q^\dagger$ and $P^\dagger\times_PQ$ are isomorphic.
The proposition follows because the integral closures in $Q_V$ of $P^{(R)}\times_PQ$ and $P^\dagger\times_PQ$ are isomorphic.

\begin{prop}\label{log53}
We keep the assumptions of \eqref{log51} and assume moreover that $g$ is étale 
and that $f$ is an isomorphism. Then $g^{(R)}\colon Q^{(R)}\rightarrow P^{(R)}$ is étale.
\end{prop}
Let $n$ be an integer $\geq 1$ such that $nR$ has integral coefficients, $X^{(n)}$ (resp. 
$Y^{(n)}$) the $n$-th infinitesimal neighborhood of $s$ (resp. $t$), and
$P^\dagger$ (resp. $Q^\dagger$) be the dilatation of $P$ along 
$X^{(n)}$ (resp. $Q$ along $Y^{(n)}$) of thickening $nR$ in the sense of \eqref{tub1}. We prove first that the morphism 
$g^\dagger\colon Q^\dagger\rightarrow P^\dagger$ 
induced by $g$ is étale \eqref{tub3a}. 
Since $g$ is étale, $Y^{(n)}$ is an open and closed subscheme of $Z=X^{(n)}\times_PQ$. 
Let $Q^\ddagger$ be the dilatation of $Q$ along $Z$ of thickening $nR$. Then the canonical
morphism $Q^\dagger\rightarrow Q^\ddagger$ is an open immersion.  Since the morphism
$g^\ddagger\colon Q^\ddagger\rightarrow P^\dagger$ induced by $g$ is étale by \ref{tub31}(iii), 
$g^\dagger$ is étale. On the other hand, by \ref{log5}, $P^{(R)}$ (resp. $Q^{(R)}$) is the integral closure 
of $P^\dagger$ in $P_U$ (resp. $Q^\dagger$ in $Q_V$), and $g^{(R)}$ is induced by $g^\dagger$. 
We deduce that the morphism $Q^{(R)}\rightarrow P^{(R)}\times_{P^\dagger}Q^\dagger$ induced by $g^{(R)}$
is an isomorphism, which implies that $g^{(R)}$ is étale.

\subsection{}\label{log7}
Let $(X,D)$ be an snc-pair over $k$, $U=X-D$ and $R$ an effective rational divisor on $X$ with support in $D$. 
We consider $X\Asterisk_kX$ as an $X$-scheme by $\pr_2$,
and denote by $\delta\colon X\rightarrow X\Asterisk_kX$ the framed diagonal of $(X,D)$ and
by $(X\Asterisk_kX)^{(R)}$ the dilatation of $X\Asterisk_kX$ along $\delta$ of thickening $R$ \eqref{log4}.
We can make the following remarks~:

{\rm (i)}\ If we consider $X\Asterisk_kX$ as an $X$-scheme by $\pr_1$ instead of $\pr_2$, 
then the dilatation of $X\Asterisk_kX$ along $\delta$ of thickening $R$
is equal to  $(X\Asterisk_kX)^{(R)}$. In particular, the automorphism of $X\Asterisk_kX$ 
switching the  factors  induces an isomorphism 
\begin{equation}\label{log7d}
\sigma\colon (X\Asterisk_kX)^{(R)}\stackrel{\sim}{\rightarrow} (X\Asterisk_kX)^{(R)}.
\end{equation}

{\rm (ii)}\ There is a canonical morphism 
\begin{equation}\label{log7b}
\delta^{(R)}\colon X\rightarrow (X\Asterisk_kX)^{(R)}
\end{equation}
lifting $\delta$, and a canonical open immersion
\begin{equation}\label{log7a}
j^{(R)}\colon U\times_kU\rightarrow (X\Asterisk_kX)^{(R)}.
\end{equation}

{\rm (iii)}\ If $R$ has integral coefficients, then the canonical projections  
$(X\Asterisk_kX)^{(R)}\rightarrow X$ are smooth \eqref{tub26} and we have canonical $R$-isomorphisms \eqref{tub25b}
\begin{equation}\label{log7c}
R\times_X(X\Asterisk_kX)^{(R)}\stackrel{\sim}{\rightarrow} \bV(\Omega^1_{X/k}(\log D)\otimes_{\co_X}\co_X(R))\times_XR
\stackrel{\sim}{\rightarrow} (X\Asterisk_kX)^{(R)}\times_XR.
\end{equation}

\subsection{}\label{log71}
Let $f\colon (Y,E)\rightarrow (X,D)$ be a morphism of snc-pairs over $k$, $U=X-D$, $V=Y-E$, 
$R$ a rational effective divisor on $X$ with support in $D$ and $R_Y=f^*(R)$. We have a commutative 
diagram 
\begin{equation}\label{log71a}
\xymatrix{
Y\ar@{=}[r]\ar[d]_{\delta_Y}&Y\ar[r]^f\ar[d]^{\gamma_f}&X\ar[d]^{\delta_X}\\
{Y\Asterisk_kY}\ar[r]^{f_1}&{X\Asterisk_kY}\ar[r]^{f_2}&{X\Asterisk_kX}}
\end{equation}
where $f_1$ and $f_2$ are the morphisms defined in \eqref{fram17a},
$\delta_X$ and $\delta_Y$ are the framed diagonals and $\gamma_f$ is the framed graph of $f$. 
We put $\tf=f_2\circ f_1$. We consider $Y\Asterisk_kY$ and $X\Asterisk_kY$ as $Y$-schemes and 
$X\Asterisk_kX$ as an $X$-scheme by the second projections. We denote by 
$(Y\Asterisk_kY)^{(R_Y)}$ (resp. $(X\Asterisk_kY)^{(R_Y)}$) the dilatation of $Y\Asterisk_kY$
(resp. $X\Asterisk_kY$) along $\delta_Y$ (resp. $\gamma_f$) of thickening $R_Y$, 
and by $(X\Asterisk_kX)^{(R)}$ the dilatation of  
$X\Asterisk_kX$ along $\delta_X$ of thickening $R$. Then we have canonical morphisms 
\begin{equation}\label{log71b}
\xymatrix{
Y\ar@{=}[r]\ar[d]_{\delta_Y^{(R_Y)}}&Y\ar[r]^f\ar[d]^{\gamma_f^{(R_Y)}}&X\ar[d]^{\delta^{(R)}_X}\\
{(Y\Asterisk_kY)^{(R_Y)}}\ar[r]^{f_1^{(R_Y)}}&{(X\Asterisk_kY)^{(R_Y)}}\ar[r]^{f_2^{(R)}}& 
{(X\Asterisk_kX)^{(R)}}}
\end{equation}
where $f_1^{(R_Y)}$ and $f_2^{(R)}$ are defined in \eqref{log51c} and the vertical morphisms are
the canonical liftings of the vertical morphisms in \eqref{log71a}; we have $\tf^{(R)}=f_2^{(R)}\circ f_1^{(R_Y)}$ . 

Since the canonical morphism $X\Asterisk_kY\rightarrow (X\Asterisk_kX)\times_XY$ is an isomorphism \eqref{fram17b}, 
the morphism $f_2^{(R)}$ is described by the proposition \ref{log52}.

\begin{prop}\label{log73}
We keep the assumptions of \eqref{log71}. 

{\rm (i)}\ If $f$ is log-étale, then the morphism 
\begin{equation}\label{log73a}
f_1^{(R_Y)}\colon (Y\Asterisk_kY)^{(R_Y)}\rightarrow (X\Asterisk_kY)^{(R_Y)}
\end{equation}
is étale. 

{\rm (ii)}\ If $f$ is log-smooth and if $R_Y=f^*(R)$ has integral coefficients, 
then the morphism $f_1^{(R_Y)}$ is smooth.
\end{prop}

(i) Indeed, $f_1$ is étale by \ref{fram18}, and hence $f_1^{(R_Y)}$ is étale by \ref{log53}. 

(ii) Indeed, $f_1$ is smooth by \ref{fram18}, and hence $f_1^{(R_Y)}$ is smooth by \ref{tub5}.

\subsection{}\label{log8} 
Let $(X,D)$ be an snc-pair over $k$. We put 
\begin{equation}\label{log8a} 
X\Asterisk_kX\Asterisk_kX=X_{\log D}\timesl_{\Spec(k)\times [M_D]}X_{\log D}\timesl_{\Spec(k)\times[M_D]}X_{\log D}. 
\end{equation}
We denote by $\pr_i\colon X\Asterisk_kX\Asterisk_kX\rightarrow X$ $(1\leq i\leq 3)$ and
$\pr_{ij}\colon X\Asterisk_kX\Asterisk_kX\rightarrow X\Asterisk_kX$ $(1\leq i<j\leq 3)$ the canonical projections,
by $\Delta\colon X\rightarrow X\Asterisk_kX\Asterisk_kX$ the unique morphism 
such that $\pr_i\circ \Delta=\id_X$ for all $1\leq i\leq 3$ and by $\delta\colon X\rightarrow X\Asterisk_kX$
the framed diagonal of $(X,D)$. It follows immediately from the definition that the diagram 
\begin{equation}\label{log8b} 
\xymatrix{
{X\Asterisk_kX\Asterisk_kX}\ar[r]^-(0.5){\pr_{23}}\ar[d]_{\pr_{12}}&{X\Asterisk_kX}\ar[d]^{\pr_1}\\
{X\Asterisk_kX}\ar[r]^-(0.5){\pr_2}&X}
\end{equation}
is cartesian. The projection $\pr_{23}$ is strict and smooth by \ref{fram12}. Then so is $\pr_{13}$ by symmetry. 

Let $R$ be an effective divisor on $X$ with support in $D$. We consider 
$X\Asterisk_kX\Asterisk_kX$ (resp. $X\Asterisk_kX$) as an $X$-scheme by $\pr_3$ (resp. $\pr_2$) and denote by 
$(X\Asterisk_kX\Asterisk_kX)^{(R)}$ (resp. $(X\Asterisk_kX)^{(R)}$) the dilatation 
of $X\Asterisk_kX\Asterisk_kX$ (resp. $X\Asterisk_kX$) along $\Delta$ (resp. $\delta$) of thickening $R$. 
Observe that if we consider $X\Asterisk_kX\Asterisk_kX$ 
as an $X$-scheme by any of the projections $\pr_i$ $(1\leq i\leq 3)$, 
the dilatation of $X\Asterisk_kX\Asterisk_kX$ along $\Delta$ of thickening $R$ does not change. 
We deduce by \ref{tub6} that we have a canonical isomorphism 
\begin{equation}\label{log8c}
(X\Asterisk_kX\Asterisk_kX)^{(R)}\stackrel{\sim}{\rightarrow} (X\Asterisk_kX)^{(R)}\times_X(X\Asterisk_kX)^{(R)}.
\end{equation}
By the universal property of dilatations \eqref{tub3}, $\pr_{13}$ induces a morphism 
\begin{equation}\label{log8d}
\mu\colon (X\Asterisk_kX)^{(R)}\times_X(X\Asterisk_kX)^{(R)}\rightarrow (X\Asterisk_kX)^{(R)}
\end{equation}
that fits into the following commutative diagram.
\begin{equation}\label{log8e}
\xymatrix{
{(X\Asterisk_kX)^{(R)}\times_X(X\Asterisk_kX)^{(R)}}\ar[rr]^-(0.5)\mu\ar[d]&&{(X\Asterisk_kX)^{(R)}}\ar[d]\\
{(X\Asterisk_kX)\times_X(X\Asterisk_kX)}\ar@{=}[r]&
{X\Asterisk_kX\Asterisk_kX}\ar[r]^-(0.5){\pr_{13}}&{X\Asterisk_kX}}
\end{equation}

\begin{prop}[\cite{saito1} 2.24]\label{log9}
Under the assumptions of \eqref{log8}, $\mu$ is smooth and $\mu\times_XR$ is the addition of the 
vector bundle $E=(X\Asterisk_kX)^{(R)}\times_XR$ over $R$ \eqref{log7c}.
\end{prop} 

Since $\pr_{13}$ is smooth, $\mu$ is smooth \eqref{tub5}. The closed subscheme $R\times_X(X\Asterisk_kX)^{(R)}$
of $(X\Asterisk_kX)^{(R)}$ is equal to $E$, and the canonical projections  $(X\Asterisk_kX)^{(R)}\rightrightarrows
X$ induce the same morphism $E\rightarrow R$. 
On the other hand, $\alpha=\mu\times_XR$ is a linear morphism of vector bundles 
$E\times_RE\rightarrow E$ \eqref{tub4}. Let $i_1,i_2\colon E\rightrightarrows E\times_RE$
be the homomorphisms defined by $i_1(x)=(x,0)$ and $i_2(x)=(0,x)$. 
To prove that $\alpha$ is the addition of $E$,
it is enough to prove that $\alpha\circ i_1=\alpha\circ i_2=\id_{E}$. 
Consider the morphism 
$\iota_1=\id\times_X\delta^{(R)}\colon (X\Asterisk_kX)^{(R)}\rightarrow (X\Asterisk_kX\Asterisk_kX)^{(R)}$. 
We have $i_1=\iota_1\times_XR\colon E\rightarrow E\times_RE$. 
Since $\mu\circ \iota_1=\id_{(X\Asterisk_kX)^{(R)}}$, 
then $\alpha\circ i_1=\id_{E}$. The same argument shows that $\alpha\circ i_2=\id_{E}$.

\section{Review of ramification theory of local fields with imperfect residue fields}\label{rrt}

\subsection{}\label{rrt1}
In this section,  $K$ denotes a discrete valuation field, $\co_K$ the valuation ring of $K$, 
$\fm_K$ the maximal ideal of $\co_K$,   
$F$ the residue field of $\co_K$, $\oK$ a separable closure of $K$, $\cG$ the Galois group of $\oK/K$
and $\ord$ the valuation of $\oK$ normalized by $\ord(K^\times)=\mZ$. 
We assume that $\co_K$ is {\em henselian} and that $F$ has characteristic $p$. 
In (\cite{as1} 3.12), we defined a decreasing filtration $\cG_{\log}^r$ $(r\in \mQ_{>0})$ of $\cG$
by closed normal subgroups, called the {\em logarithmic ramification filtration}. 
Unlike the convention in loc. cit., it is more convenient to extend it by letting $\cG_{\log}^0$ 
be the inertia subgroup of $\cG$.  For a rational number $r\geq 0$, 
we put 
\begin{eqnarray}
\cG_{\log}^{r+}&=&\overline{\bigcup_{s>r}\cG_{\log}^{s}},\label{rrt1a}\\
\Gr^r_{\log}(\cG)&=&\cG_{\log}^r/\cG_{\log}^{r+}.\label{rrt1b}
\end{eqnarray}
Then $\cP=\cG_{\log}^{0+}$ is the wild inertia subgroup of $\cG$, i.e., the $p$-Sylow subgroup of $\cG_{\log}^0$
(\cite{as1} 3.15). 

\subsection{}\label{rrt15}
Let $L$ be a finite separable extension of $K$ and $r$ a rational number $\geq 0$. Then $\cG$ acts on 
$\Hom_K(L,\oK)$ via its action on $\oK$. We say that the ramification
of $L/K$ is {\em bounded} by $r$ (resp. by $r+$) if $\cG^{r}_{\log}$ (resp. $\cG^{r+}_{\log}$) 
acts trivially on $\Hom_K(L,\oK)$. We define 
the {\em conductor} $c$ of $L/K$ as the infimum of rational numbers $r> 0$ such that the ramification 
of $L/K$ is bounded by $r$. Then $c$ is a rational number, and 
the ramification of $L/K$ is bounded by $c+$ (\cite{as1} 9.5). If $c>0$, the ramification 
of $L/K$ is not bounded by $c$.

\begin{teo}[\cite{as2} Theorem 1]\label{rrt2}
For every rational number $r>0$, the group $\Gr^r_{\log}(\cG)$ is abelian 
and is contained in the center of the pro-$p$-group $\cP/\cG^{r+}_{\log}$. 
\end{teo} 

\begin{lem}[\cite{katz} 1.1] \label{rrt3}
Let $M$ be a $\mZ[\frac 1p]$-module on which $\cP$ acts through a 
finite discrete quotient, say by $\rho\colon \cP\rightarrow \Aut_\mZ(M)$. 
Then, 

{\rm (i)}\ $M$ has a unique direct sum decomposition 
\begin{equation}
M=\oplus_{r\in \mQ_{\geq 0}} M^{(r)}
\end{equation}
into $\cP$-stable submodules $M^{(r)}$, such that $M^{(0)}=M^{\cP}$ and  
for every $r>0$, 
\begin{equation}
(M^{(r)})^{\cG^r_{\log}}=0\ \ \ {\rm and}\ \ \ (M^{(r)})^{\cG^{r+}_{\log}}=M^{(r)}.
\end{equation}

{\rm (ii)}\ If $r>0$, then $M^{(r)}=0$ for all but the finitely many values of $r$
for which $\rho(\cG^r_{\log})\supsetneq  \rho(\cG^{r+}_{\log})$. 

{\rm (iii)}\ For variable $M$ but fixed $r$, 
the functor $M\mapsto M^{(r)}$ is exact. 

{\rm (iv)}\ For $M,N$ as above, 
we have $\Hom_{\cP{\rm \text-mod}}(M^{(r)},N^{(r')})=0$ if $r\not=r'$. 
\end{lem}

\begin{defi}\label{rrt4} 
The  decomposition $M=\oplus_{r\in \mQ_{\geq 0}} M^{(r)}$ of lemma \ref{rrt3} is called
the {\em slope decomposition} of $M$. The values $r\geq 0$ 
for which $M^{(r)}\not=0$ are called the {\em slopes} of $M$. 
We say that $M$ is {\em isoclinic} if it has only one slope. 
\end{defi}

These notions apply in particular to the case where $M$ is a $\mZ[\frac 1p]$-module on which $\cG$ acts through a 
finite discrete quotient.

\begin{lem}[\cite{katz} 1.4]\label{rrt5} 
If $A$ is a $\mZ[\frac 1p]$-algebra and $M$ is a left $A$-module on which $\cP$
acts $A$-linearly through a finite discrete quotient, then in the slope
decomposition $M=\oplus_r M^{(r)}$, each $M^{(r)}$ is an sub-$A$-module of $M$. 
For any $A$-algebra $B$, the slope decomposition of $B\otimes_AM$ is given 
by $B\otimes_AM=\oplus_r (B\otimes_AM^{(r)})$.
\end{lem}

\begin{lem}\label{rrt6}
Let $M$ be a $\Lambda$-module on which $\cP$ acts $\Lambda$-linearly 
through a finite discrete quotient, which is isoclinic of slope
$r>0$; so the subgroup $\cG^{r+}_{\log}$ acts trivially on $M$. 

{\rm (i)}\ Let $X(r)$ be the set of isomorphism classes of finite characters 
$\chi\colon \cG^r_{\log}/\cG^{r+}_{\log}\rightarrow \Lambda_\chi^\times$
such that $\Lambda_\chi$ is a finite \'etale $\Lambda$-algebra, generated 
by the image of $\chi$, and having a connected spectrum. 
Then $M$ has a unique direct sum decomposition 
\begin{equation}
M=\oplus_{\chi\in X(r)} M_\chi
\end{equation}
into $\cP$-stable sub-$\Lambda$-modules $M_\chi$, such that 
$\Lambda[\cG^r_{\log}]$ acts on $M_\chi$ through 
$\Lambda_\chi$ for every $\chi$. 

{\rm (ii)}\ There are finitely many characters $\chi\in X(r)$ for which $M_\chi\not=0$.
 
{\rm (iii)}\ For variable $M$ but fixed $\chi$, 
the functor $M\mapsto M_\chi$ is exact. 

{\rm (iv)}\ For $M,N$ as above, 
we have $\Hom_{\Lambda[\cP]}(M_\chi,N_{\chi'})=0$ if $\chi\not=\chi'$. 
\end{lem}
Let $P$ be a finite discrete quotient of $\cP/\cG^{r+}_{\log}$ 
through which $\cP$ acts on $M$ and let $C$ be the image of $\cG^r_{\log}/\cG^{r+}_{\log}$ in $P$. 
We know by \ref{rrt2} that $C$ is contained in the center of $P$. 
The connected components of $\Spec(\Lambda[C])$ correspond 
to the isomorphism classes of characters $\chi\colon C\rightarrow A^\times$,
where $A$ is a finite \'etale $\Lambda$-algebra, generated by the image 
of $\chi$, and having a connected spectrum.  
We obtain a set of orthogonal idempotents $e_\chi$ of $\Lambda[C]$,
indexed by such characters, whose sum is $1$, and 
are clearly central in $\Lambda[P]$. The lemma follows.

\begin{rema}\label{rrt7}
If $p^nC=0$ and $\Lambda$ contains a primitive 
$p^n$-th root of unity, then $\Lambda_\chi= \Lambda$ for every $\chi$
such that $M_\chi\not=0$.  
\end{rema}

\begin{defi}\label{rrt8}
The decomposition $M=\oplus_{\chi} M_\chi$ of lemma \ref{rrt6} is called
the {\em central character decomposition} of $M$. The characters 
$\chi\colon \cG^r_{\log}/\cG^{r+}_{\log}\rightarrow \Lambda_\chi^\times$ 
for which $M_\chi\not=0$ are called the {\em central characters} of $M$. 
\end{defi}

\begin{lem}\label{rrt9}
If $A$ is a $\Lambda$-algebra, and $M$ is a left $A$-module 
on which $\cP$ acts $A$-linearly through a finite discrete quotient,
which is isoclinic, 
then in the central character  decomposition $M=\oplus_\chi M_\chi$, 
each $M_\chi$ is an sub-$A$-module of $M$. 
For any $A$-algebra $B$, the central character decomposition 
of $B\otimes_AM$ is given by $B\otimes_AM=\oplus_\chi(B\otimes_AM_\chi)$. 
\end{lem}

This is clear from the proof of \ref{rrt6}. 

\subsection{}\label{rrt10}
For the rest of this section, we assume that $K$ has characteristic $p$ and that $F$ is of finite type over $k$. We 
denote by $\Omega^1_{\co_K}(\log)$ the $\co_K$-module of logarithmic $1$-differential forms
\begin{equation}
\Omega^1_{\co_K}(\log)=(\Omega^1_{\co_K}\oplus(\co_K\otimes_{\mZ}K^\times))/\fS,
\end{equation}
where $\fS$ is the sub-$\co_K$-module of $\Omega^1_{\co_K}\oplus(\co_K\otimes_{\mZ}K^\times)$ 
generated by elements of the form $(da,0)-(0,a\otimes a)$, for $a\in \co_K-\{0\}$. For every $a\in K^\times$, 
we denote by $d\log(a)\in \Omega^1_{\co_K}(\log)$ the class of $(0,1\otimes a)$.
Let $\pi$ be a uniformizer of $\co_K$. The morphism 
$\Omega^1_{\co_K}\oplus \co_K\rightarrow \Omega^1_{\co_K}(\log)$ that maps $(\omega,a)$ to
$\omega+ad\log(\pi)$ is surjective, and its kernel is generated by $(d\pi,0)-(0,\pi)$. 

We put $\Omega^1_{F}(\log)=\Omega^1_{\co_K}(\log)\otimes_{\co_K}F$ and denote by $d\log[\pi]$ the class of $d\log(\pi)$. 
It is easy to see that $\Omega^1_{F}(\log)$ is canonically isomorphic to the quotient of 
$\Omega^1_F\oplus(F\otimes_{\mZ}K^\times)$
by the sub-$F$-module generated by elements of the form $(d\oa,0)-(0,\oa\otimes a)$, for $a\in \co_K-\{0\}$,
where $\oa$ is the residue class of $a$ in $F$. Then we have an exact sequence
\begin{equation}
0\longrightarrow \Omega^1_F\longrightarrow \Omega^1_F(\log)\stackrel{\res}{\longrightarrow} F\longrightarrow 0,
\end{equation}
where $\res((0,a\otimes b))=a\cdot \ord(b)$ for $a\in F$ and $b\in K^\times$. In particular, $\Omega^1_F(\log)$ 
is an $F$-vector space of finite dimension. 

\subsection{}\label{rrt11}
Let $\co_{\oK}$ be the integral closure of $\co_K$ in $\oK$ and $\oF$ its residue field. 
For a rational number $r$, we put $\fm_{\oK}^r=\{x\in \oK\ ;\  \ord(x)\geq r\}$,  
$\fm_{\oK}^{r+}=\{x\in \oK\ ;\ \ord(x)> r\}$ and 
\begin{equation}\label{rrt11a}
\Theta^{(r)}_{\oF}=\Hom_F(\Omega^1_F(\log),\fm_{\oK}^{r}/\fm_{\oK}^{r+}),
\end{equation}
which is an $\oF$-vector space of finite dimension. We consider $\Theta^{(r)}_{\oF}$ 
as a smooth additive $\oF$-group-scheme. Let $\pi_1^\alg(\Theta^{(r)}_{\oF})$ be the quotient
of the abelianized fundamental group $\pi_1^\ab(\Theta^{(r)}_{\oF})$ of $\Theta^{(r)}_{\oF}$ classifying
étale isogenies; it is an abelian profinite group killed by $p$. Recall that Lang's isogeny 
$\mA^1_{\oF}\rightarrow \mA^1_{\oF}$, defined by $x\mapsto x^p-x$ (where $x$ is the canonical parameter 
of $\mA^1_{\oF}$), is a basis of the 
$\oF$-vector space $\Hom_{\mZ}(\pi_1^{\alg}(\mA^1_{\oF}), \mF_p)$. Therefore, we have a canonical
isomorphism 
\begin{equation}\label{rrt11b}
\Hom_{\mZ}(\pi_1^\alg(\Theta^{(r)}_{\oF}),\mF_p)\stackrel{\sim}{\rightarrow}
\Hom_{\oF}(\fm_{\oK}^{r}/\fm_{\oK}^{r+}, \Omega^1_F(\log)\otimes_F\oF).
\end{equation}

\begin{teo}[\cite{saito1} 1.24]\label{rrt12}
For every rational number $r>0$, there exists a canonical surjective homomorphism 
\begin{equation}\label{rrt12a}
\pi_1^{\alg}(\Theta^{(r)}_{\oF})\rightarrow \Gr^r_{\log}\cG.
\end{equation}
Consequently, the group $\Gr^r_{\log}\cG$ is killed by $p$, and we have a canonical injective homomorphism 
\begin{equation}\label{rrt12b}
\rsw\colon \Hom_{\mZ}(\Gr^r_{\log}\cG,\mF_p)\rightarrow 
\Hom_{\oF}(\fm_{\oK}^{r}/\fm_{\oK}^{r+}, \Omega^1_F(\log)\otimes_F\oF).
\end{equation}
\end{teo}

The homomorphism \eqref{rrt12b} is called the {\em refined Swan conductor}. 
This theorem has been recently extended to the unequal characteristic case by one of the authors (T.S.) \cite{saito2}.

\section{Ramification of Galois torsors}\label{ram-cov}
 
\subsection{}\label{ram-cov1}
In this section, we fix an snc-pair $(X,D)$ over $k$, and put $U=X-D$. We denote by $j\colon U\rightarrow X$
the canonical injection, by $X\Asterisk_kX$ the framed self-product of $(X,D)$ 
and by $\delta\colon X\rightarrow X\Asterisk_kX$ the framed diagonal \eqref{fram17}. 
We will no longer consider the logarithmic structure of $X\Asterisk_kX$;
only the underlying scheme will be of interest for us. 
We consider $X\Asterisk_kX$ as an $X$-scheme by the second projection. 
For any effective rational divisor $R$ on $X$ with support in $D$, we denote by $(X\Asterisk_kX)^{(R)}$
the dilatation of $X\Asterisk_kX$ along $\delta$ of thickening $R$ \eqref{log7}, by  
\begin{equation}\label{ram-cov1a}
\delta^{(R)}\colon X\rightarrow (X\Asterisk_kX)^{(R)}
\end{equation}
the canonical lifting of $\delta$, and by 
\begin{equation}\label{ram-cov1b}
j^{(R)}\colon U\times_kU\rightarrow (X\Asterisk_kX)^{(R)}
\end{equation}
the canonical open immersion. Then we have the following Cartesian diagram. 
\begin{equation}\label{ram-cov1c}
\xymatrix{
U\ar[r]^-(0.5){\delta_U}\ar[d]_j&{U\times_kU}\ar[d]^-(0.5){j^{(R)}}\\
X\ar[r]^-(0.5){\delta^{(R)}}&{(X\Asterisk_kX)^{(R)}}}
\end{equation}

\subsection{}\label{ram-cov2}
Let  $V$ be a Galois torsor over $U$ of group $G$ and $R$ an effective rational divisor on $X$ with support in $D$.
We consider $V\times_kV$ as a Galois torsor over $U\times_kU$
of group $G\times G$, and denote by $W$ the quotient of $V\times_kV$ by the group $\Delta(G)$, 
where $\Delta\colon G\rightarrow G\times G$ is the diagonal homomorphism.
The diagonal morphism $\delta_V\colon V\rightarrow V\times_kV$
induces a morphism $\varepsilon_U\colon U\rightarrow W$ above the diagonal morphism 
$\delta_U\colon U\rightarrow U\times_kU$.
Note that $W$ represents the sheaf of isomorphisms of $G$-torsors
from $U\times_kV$ to $V\times_kU$ over $U\times_kU$, and that $\varepsilon_U$ corresponds to the 
identity isomorphism of $V$ (identified with the pull-backs of $U\times_kV$ and $V\times_kU$
by $\delta_U$).
We denote by $Z$ the integral closure of $(X\Asterisk_kX)^{(R)}$ in $W$,
by $\pi\colon Z\rightarrow (X\times_kX)^{(R)}$ the canonical morphism and by $\varepsilon\colon X\rightarrow Z$
the morphism induced by $\varepsilon_U\colon U\rightarrow W$. We have $\pi\circ \varepsilon=\delta^{(R)}$.
\begin{equation}\label{ram-cov3a}
\xymatrix{
&W\ar[r]\ar[d]&Z\ar[d]^{\pi}&\\
U\ar@/^1pc/[ru]^{\varepsilon_U}\ar[r]^(0.4){\delta_U}&{U\times_kU}\ar[r]&{(X\Asterisk_kX)^{(R)}}&
X\ar@/_1pc/[ul]_{\varepsilon}\ar[l]_-(0.4){\delta^{(R)}}}
\end{equation}

\begin{defi}\label{ram-cov3}
We keep the assumptions of \eqref{ram-cov2} and let $x\in X$. 
We say that the {\em ramification of $V/U$ at $x$ is bounded by $R+$}
if the morphism $\pi\colon Z\rightarrow (X\Asterisk_kX)^{(R)}$ is étale at $\varepsilon(x)$,
and that the {\em ramification of $V/U$ along $D$ is bounded by $R+$}
if $\pi$ is étale over an open neighborhood of $\varepsilon(X)$. 
\end{defi}

\begin{lem}\label{ram-cov75}
Let $V$ be a Galois torsor over $U$, $R$ an effective rational divisor on $X$
with support in $D$ and $x\in X$.  
The ramification of $V/U$ at $x$ is bounded by $R+$ if and only if there exists an open neighborhood 
$X'$ of $x$ in $X$ such that if we put $U'=U\times_XX'$ and $V'=V\times_XX'$ and if we denote by $D'$
and $R'$ the pull-backs of $D$ and $R$ over $X'$, then the ramification of $V'/U'$ along $D'$ is bounded by $R'+$. 
\end{lem}

Only the necessity of the condition requires a proof. Assume that the ramification of $V/U$ at $x$ is bounded by $R+$.
Then with the notation of \eqref{ram-cov2}, there exists an open neighborhood $X'$ of $x$ in $X$ 
such that the morphism $\pi\colon Z\rightarrow (X\Asterisk_kX)^{(R)}$ is étale at each point $\varepsilon(x')$ for $x'\in X'$.  
It is clear that $X'$ satisfies the required property.

\begin{lem}\label{ram-cov38}
Let $V$ be a Galois torsor over $U$, $x\in X$ and $R$ and $R'$ effective rational divisors 
on $X$ with supports in $D$ such that $R'\geq R$. Then if the ramification of $V/U$ at $x$ (resp. along $D$) 
is bounded by $R+$, it is also bounded by $R'+$.
\end{lem}

We use the notation of \eqref{ram-cov2} both for $R$ and $R'$; we equip objects relative to $R'$ with a prime. 
There is a canonical morphism $u\colon (X\Asterisk_kX)^{(R')}\rightarrow (X\Asterisk_kX)^{(R)}$ \eqref{log4d}
that fits into the following commutative diagram.
\begin{equation}\label{ram-cov38a}
\xymatrix{
&{(X\Asterisk_kX)^{(R')}}\ar[d]^u&\\
X\ar[ru]^{\delta^{(R')}}\ar[r]^-(0.6){\delta^{(R)}}&{(X\Asterisk_kX)^{(R)}}&{U\times_kU}\ar[lu]\ar[l]}
\end{equation}
It induces a canonical morphism $v\colon Z'\rightarrow Z$ that fits into 
the following commutative diagram. 
\begin{equation}\label{ram-cov38b}
\xymatrix{
&Z'\ar[d]^v\ar[r]&{(X\Asterisk_kX)^{(R')}}\ar[d]^u\\
X\ar[ru]^{\varepsilon'}\ar[r]^-(0.5){\varepsilon}&Z\ar[r]&{(X\Asterisk_kX)^{(R)}}}
\end{equation}
Moreover, $Z'$ is the integral closure of $u^*(Z)$ in $W$. Let $\uZ$ (resp. $\uZ'$)
be the maximal open subscheme of $Z$ (resp. $Z'$) which is étale over $(X\Asterisk_kX)^{(R)}$
(resp. $(X\Asterisk_kX)^{(R')}$). Since $(X\Asterisk_kX)^{(R')}$ is normal \eqref{log5},
$u^*(\uZ)$ is normal. Therefore, we can identify $u^*(\uZ)$ with $v^{-1}(\uZ)\subset Z'$,
and we have $v^{-1}(\uZ)\subset \uZ'$, which implies the proposition. 

\begin{prop}\label{ram-cov31}
Let $V$ be a Galois torsor over $U$ of group $G$, $x\in X$, $H$ a normal subgroup of $G$ and 
$V'$ the quotient of $V$ by $H$. Then if the ramification of $V/U$ at $x$ (resp. along $D$) is bounded by $R+$,
the ramification of $V'/U$ at $x$ (resp. along $D$) is bounded by $R+$.
\end{prop}
It is enough to prove the proposition relative to $x$. 
We use the notation of \eqref{ram-cov2} for the Galois torsor $V$ over $U$.
We put $G'=G/H$ and denote by $\Delta'\colon G'\rightarrow G'\times G'$ the diagonal homomorphism,
by $W'$ the quotient of $V'\times_kV'$ by $\Delta'(G')$ and by $Z'$ the integral closure of $(X\Asterisk_kX)^{(R)}$
in $W'$. Let $H'$ be the subgroup of $G\times G$ of elements $(g,g')$ such that $g'g^{-1}\in H$
(i.e., the inverse image of $\Delta'(G')$ in $G\times G$). Then $W'$ is the quotient of $V\times_kV$
by $H'$. Since $\Delta(G)\subset H'$, there exists a canonical $(V\times_kV)$-morphism $W\rightarrow W'$, 
which induces an $(X\Asterisk_kX)^{(R)}$-morphism $\rho\colon Z\rightarrow Z'$. 
By assumption, there exists  an open neighborhood $Z_0$ of $\varepsilon(x)$ in $Z$, 
which is étale over $(X\Asterisk_kX)^{(R)}$. 
Then $Z_0$ is unramifed over $Z'$, and by (\cite{ega4} 18.10.1), $Z_0$ is étale over $Z'$; 
in particular, $Z_0$ is flat over $Z'$ and $\rho(Z_0)$ is an open subscheme of $Z'$. 
We conclude by (\cite{ega4} 17.7.7) that $\rho(Z_0)$ is étale over $(X\Asterisk_kX)^{(R)}$, which 
implies that the ramification of $V'/U$ at $x$ is bounded by $R+$.

\begin{prop}\label{ram-cov32}
Let $V$ be a Galois torsor over $U$, $R$ an effective rational divisor on $X$
with support in $D$, $f\colon (X',D')\rightarrow (X,D)$ a morphism of snc-pairs,
$U'=X'-D'$, $V'=V\times_UU'$, $R'=f^*(R)$, $x'\in X'$ and $x=f(x')$. Then~:

{\rm (i)}\ If the ramification of $V/U$ at  $x$ (resp. along $D$) is bounded by $R+$, 
the ramification of $V'/U'$ at $x'$ (resp. along $D'$) is bounded by $R'+$.

{\rm (ii)}\ Assume that $f$ is log-smooth \eqref{log1} and that the morphism
$X'\rightarrow X$ is flat at $x'$. Then the ramification of $V/U$ at $x$ is bounded by $R+$
if and only if the ramification of $V'/U'$ at $x'$ is bounded by $R'+$.

{\rm (iii)}\ Assume that $f$ is log-smooth and that the morphism
$X'\rightarrow X$ is faithfully flat. Then the ramification of $V/U$ along $D$ is bounded by $R+$
if and only if the ramification of $V'/U'$ along $D'$ is bounded by $R'+$.
\end{prop}

Let $G$ be the Galois group of $V/U$. We denote by $W$ (resp. $W'$) the quotient of $V\times_kV$ 
(resp. $V'\times_kV'$) by $\Delta(G)$ (cf. \ref{ram-cov2}). We have a commutative diagram \eqref{log71a}
\begin{equation}
\xymatrix{
X'\ar@{=}[r]\ar[d]_{\delta_{X'}}&X'\ar[r]^f\ar[d]^{\gamma_f}&X\ar[d]^{\delta_X}\\
{X'\Asterisk_kX'}\ar[r]^{f_1}&{X\Asterisk_kX'}\ar[r]^{f_2}&{X\Asterisk_kX}}
\end{equation}
where $f_1$ and $f_2$ are the morphisms defined in \eqref{fram17a},
$\delta_X$ and $\delta_{X'}$ are the framed diagonals and $\gamma_f$ is the framed graph of $f$. 
We put $\tf=f_2\circ f_1$. Let $U_1=f^{-1}(U)$ and let $U_1\Asterisk_kU_1$ be the framed self-product of $(U_1,D'|U_1)$.
Then we have $U_1\Asterisk_kU_1=(X'\Asterisk_kX')\times_{(X'\times_kX')}(U_1\times_kU_1)$. 
We put $W_1=W\times_{(U\times_kU)}(U\times_kU_1)$ and $\tW=W\times_{(U\times_kU)}(U_1\Asterisk_kU_1)$. 
We have the following commutative diagram (cf. \ref{ram-cov2}).
\begin{equation}\label{ram-cov32a}
\xymatrix{
U'\ar[r]\ar[d]&{U_1}\ar@{=}[r]\ar[d]&{U_1}\ar[r]\ar[d]&U\ar[d]\\
W'\ar[r]\ar[d]\ar@{}[rd]|{\Box}&{\tW}\ar[r]\ar[d]\ar@{}[rd]|{\Box}&{W_1}\ar[r]\ar[d]\ar@{}[rd]|{\Box}&W\ar[d]\\
{U'\times_kU'}\ar[r]\ar[dr]&{U_1\Asterisk_kU_1}\ar[r]\ar[d]\ar@{}[rd]|{\Box}&{U\times_kU_1}\ar[r]\ar[d]\ar@{}[rd]|{\Box}&
{U\times_kU}\ar[d]\\
&{X'\Asterisk_kX'}\ar[r]^{f_1}&{X\Asterisk_kX'}\ar[r]^{f_2}&{X\Asterisk_kX}}
\end{equation}
We denote by $(X'\Asterisk_kX')^{(R')}$ (resp. $(X\Asterisk_kX')^{(R')}$)
the dilatation of  $X'\Asterisk_kX'$ (resp. $X\Asterisk_kX'$) along $\delta_{X'}$ (resp. $\gamma_f$) of thickening $R'$. 
Let $Z$ (resp. $Z_1$, resp. $Z'$) be the integral closure of $(X\Asterisk_kX)^{(R)}$ in $W$
(resp. $(X\Asterisk_kX')^{(R')}$ in $W_1$, resp. $(X'\Asterisk_kX')^{(R')}$ in $W'$).
Since $W'$ is a dense open subscheme of $\tW$ and since the latter is regular,
$Z'$ is identified with the integral closure of $(X'\Asterisk_kX')^{(R')}$ in $\tW$.
Then we have the following commutative diagram \eqref{log71b}.
\begin{equation}\label{ram-cov32b}
\xymatrix{
X'\ar@{=}[r]\ar[d]_{\varepsilon'}&X'\ar[r]^h\ar[d]_{\varepsilon_1}&X\ar[d]^{\varepsilon}\\
{Z'}\ar[r]^{f'_1}\ar[d]_-(0.5){\pi'}&{Z_1}\ar[r]^{f'_2}\ar[d]_-(0.5){\pi_1}&Z\ar[d]^{\pi}\\
{(X'\Asterisk_kX')^{(R')}}\ar[r]^{f_1^{(R')}}&{(X\Asterisk_kX')^{(R')}}\ar[r]^{f_2^{(R)}}&{(X\Asterisk_kX)^{(R)}}}
\end{equation}
We put $\tf^{(R)}=f^{(R')}_2\circ f^{(R)}_1$ and $\tf'=f'_2\circ f'_1$.
Let $\uZ$ (resp. $\uZ'$, resp. $\uZ_1$) be the maximal open subscheme of $Z$ (resp. $Z'$, resp. $Z_1$) 
where $\pi$ (resp. $\pi'$, resp. $\pi_1$) is étale. 

(i) It is enough to prove the local proposition.
We denote by $Z^\dagger$ (resp. $\uZ^\dagger$) the base change of $Z$ (resp. $\uZ$) by $\tf^{(R)}$. 
Since $(X'\Asterisk_kX')^{(R')}$ is normal \eqref{log5}, $\uZ^\dagger$ is normal.
Therefore, we can identify $\uZ^\dagger$ with $\tf'^{-1}(\uZ)\subset Z'$,
and we have $\tf'^{-1}(\uZ)\subset \uZ'$. We have $\varepsilon(x)\in \uZ$ by assumption~; 
then $\varepsilon'(x')\in \uZ'$.

(ii) By (i), there is only one implication to prove~: we assume that $\pi'$ is étale at
$\varepsilon'(x')$ and prove that $\pi$ is étale at $\varepsilon(x)$. We proceed in three steps. 

(A) Assume first that the morphism $f\colon X'\rightarrow X$
is smooth (in the usual sense) and that $R'$ has integral coefficients.  
Since $f$ is log-smooth, the first condition is satisfied if for instance the morphism 
$X'_{\log D'}\rightarrow X_{\log D}$ is strict \eqref{log1}. 
By \ref{log73}(ii), the morphism $f_1^{(R')}$ is smooth. Then the diagram 
\begin{equation}\label{ram-cov32c}
\xymatrix{
Z'\ar[r]^{f'_1}\ar[d]_{\pi'}&{Z_1}\ar[d]^{\pi_1}\\
{(X'\Asterisk_kX')^{(R')}}\ar[r]^{f_1^{(R')}}&{(X\Asterisk_kX')^{(R')}}}
\end{equation}
is cartesian, and hence $\pi_1$ is étale at $\varepsilon_1(x')$.
On the other hand, $f_2$ induces an isomorphism \eqref{fram17b}
\begin{equation}
X\Asterisk_kX' \stackrel{\sim}{\rightarrow} (X\Asterisk_kX)\times_XX'.
\end{equation}
It follows from \ref{log52} that the diagram 
\begin{equation}\label{ram-cov32d}
\xymatrix{
{(X\Asterisk_kX')^{(R')}}\ar[r]^{f_2^{(R)}}\ar[d]&{(X\Asterisk_kX)^{(R)}}\ar[d]\\
X'\ar[r]^f&X}
\end{equation}
is cartesian. Then $f_2^{(R)}$ is smooth and the diagram 
\begin{equation}
\xymatrix{
{Z_1}\ar[r]^{f'_2}\ar[d]_{\pi_1}&{Z}\ar[d]^{\pi}\\
{(X\Asterisk_kX')^{(R')}}\ar[r]^{f_2^{(R)}}&{(X\Asterisk_kX)^{(R)}}}
\end{equation}
is cartesian. Hence $\pi$ is étale at $\varepsilon(x)$.

(B) Assume next that the following conditions are satisfied~:

(a) The irreducible components $D_1,\dots,D_m$ of $D$ are defined by equations 
$t_1,\dots,t_m\in \Gamma(X,\co_X)$.

(b) There exists integers $a_1,\dots,a_m\geq 1$ such that 
\begin{equation}\label{ram-cov32e}
X'=\frac{X[S_i,U_i^{\pm1}; 1\leq i\leq m]}{(U_iS_i^{a_i}-t_i; 1\leq i\leq m)},
\end{equation}
and that the divisor $D'$ is defined by the equation $\prod_{i=1}^m S_i$. 
Observe that $(X',D')$ is an snc-pair over $k$, that $(X',D')\rightarrow (X,D)$ is log-smooth and that $X'\rightarrow X$
is faithfully flat.

(c) $R'$ has integral coefficients. 

It follows from the first half of the proof of case (A) that $\pi_1$ is étale at $\varepsilon_1(x')$.
Consider the $X$-scheme
\begin{equation}\label{ram-cov32f}
X''=\frac{X[T_i,U_i^{\pm1}; 1\leq i\leq m]}{(U_iT_i-t_i; 1\leq i\leq m)},
\end{equation}
and denote by $h\colon X''\rightarrow X$ the structural morphism (which is smooth). 
Let  $g\colon X'\rightarrow X''$ be the finite $X$-morphism defined by $T_i\mapsto S_i^{a_i}$ $(1\leq i\leq m)$ and let
$x''=g(x')$. The morphism $f_2$ factors as 
\begin{equation}
X\Asterisk_kX'\stackrel{g_2}{\longrightarrow} X\Asterisk_kX''\stackrel{h_2}{\longrightarrow} X\Asterisk_kX.
\end{equation}
We put $R''=h^*(R)$, $U_2=h^{-1}(U)$, $W_2=W\times_{(U\times_kU)}(U\times_kU_2)$, and denote by 
$(X\Asterisk_kX'')^{(R'')}$ the dilatation of $X\Asterisk_kX''$ along the framed graph of $h$ of thickening $R''$,
and by $Z_2$ the integral closure of $(X\Asterisk_kX'')^{(R'')}$ in $W_2$. We have the following commutative diagram.
\begin{equation}\label{ram-cov32g}
\xymatrix{
X'\ar[r]^g\ar[d]_{\varepsilon_1}&X''\ar[r]^h\ar[d]_{\varepsilon_2}&X\ar[d]^{\varepsilon}\\
{Z_1}\ar[r]^{g'_2}\ar[d]_-(0.5){\pi_1}&{Z_2}\ar[r]^{h'_2}\ar[d]_-(0.5){\pi_2}&Z\ar[d]^{\pi}\\
{(X\Asterisk_kX')^{(R')}}\ar[r]^{g_2^{(R'')}}\ar[d]&{(X\Asterisk_kX'')^{(R'')}}\ar[r]^{h_2^{(R)}}\ar[d]&{(X\Asterisk_kX)^{(R)}}\ar[d]\\
X'\ar[r]^g&X''\ar[r]^h&X}
\end{equation}
By \ref{log52}, each irreducible component of $(X\Asterisk_kX')^{(R')}$
dominates an irreducible component of $(X\Asterisk_kX'')^{(R'')}$.
Then it follows from  \ref{rc2} that $\pi_2$ is étale at $\varepsilon_2(x'')$. 
Since $h$ is smooth, we deduce, as in the second half of the proof of case (A), 
that $\pi$ is étale at $\varepsilon(x)$.

(C) We consider finally the general case. We may assume that the morphism  
$X'\rightarrow X$ is flat and that the ramification of $V'/U'$ along $D'$ is bounded 
by $R'+$ \eqref{ram-cov75}. Let $D_1,\dots,D_n$ be the irreducible components of $D$ containing $x$. 
By \ref{rc1},  we may assume that $x'$ is contained in exactly $n$ irreducible components $D'_1,\dots,D'_n$ of $D'$
such that $D'_i$ dominates $D_i$ for all $1\leq i\leq n$. 
Moreover, we may assume that $D=\bigcup_{1\leq i\leq n}D_i$
and $D'=\bigcup_{1\leq i\leq n}D'_i$, and that for each $1\leq i\leq n$, 
$D_i$ is defined by an equation $t_i\in \Gamma(X,\co_X)$
and $D'_i$ is defined by a section $s_i\in \Gamma(X',\co_{X'})$. 
We write $R=\sum_{i=1}^nr_iD_i$ and $f^*(D_i)=e_iD'_i$ $(1\leq i\leq n)$; 
so we have $t_i=u_is^{e_i}_i$, where $u_i\in \Gamma(X,\co_X^\times)$. 
For each $1\leq i\leq n$, let $b_i$ be an integer $\geq 1$ such that $b_ir_i$ is an integer.
We put $a_i=b_ie_i$ $(1\leq i\leq n)$ and 
\begin{equation}
Y=\frac{X[S_i,U_i^{\pm1}; 1\leq i\leq n]}{(U_iS_i^{a_i}-t_i; 1\leq i\leq n)},
\end{equation}
that we equip with the simple normal crossing divisor $E$ defined by $\prod_{i=1}^nS_i$. 
Consider the logarithmic scheme \eqref{log1}
\begin{equation}
Y'=Y_{\log E}\timesl_{X_{\log D}}X'_{\log D'}.
\end{equation}
Since $e_i$ divides $a_i$ for all $1\leq i\leq n$, the morphism $Y'\rightarrow Y_{\log E}$
is strict, and since it is log-smooth, the morphism of the underlying schemes $Y'\rightarrow Y$ is smooth.
Hence, if $E'$ denotes the pull-back of $E$ on $Y'$, 
$(Y',E')$ is an snc-pair over $k$, and the logarithmic structure on $Y'$ is defined by $E'$. 
On the other hand, $Y'\rightarrow X'\times_XY$ is finite and dominant. Since $Y\rightarrow X$ is faithfully flat, 
there exists a point $y'\in Y'$ above $x'$.  We have the following commutative diagram of $\SNCP_k$.
\begin{equation}
\xymatrix{
{(Y',E')}\ar[r]^{\alpha}\ar[d]_{\gamma}&{(Y,E)}\ar[d]^{\beta}\\
{(X',D')}\ar[r]^f&{(X,D)}}
\end{equation}
By applying first (i) to the morphism $\gamma$ at the point $y'$,  then (ii) case (A) to the morphism $\alpha$ at 
the point $y'$, and finally (ii) case (B) to the morphism $\beta$ at the point $\alpha(y')$, 
we conclude that the ramification of $V/U$ at $x$ is bounded by $R+$.

(iii) It follows from (ii). 

\begin{cor}\label{ram-cov33}
Let $V$ be a Galois torsor over $U$ of group $G$, $I$ a subgroup of $G$ and
$R$ an effective rational divisor on $X$ with support in $D$. 
We denote by $U'$ the quotient of $V$ by $I$, by $X'$ the integral closure of $X$ in $U'$, 
by $f\colon X'\rightarrow X$ the structural morphism.
Let $X'_0$ be an open subscheme of $X'$ which is étale over $X$, $x'\in X'_0$ and $x=f(x')$. 
We put $U'_0=U'\times_{X'}X'_0$, $V_0=V\times_{U'}U'_0$ and denote by 
$D'_0$ and $R'_0$ the pull-backs of $D$ and $R$ over $X'_0$; so $(X'_0,D'_0)$ is an snc-pair over $k$.  
Then the ramification of $V/U$ at $x$ is bounded by $R+$ 
if and only if the ramification of $V_0/U'_0$ at $x'$ is bounded by $R'_0+$.
\end{cor}

By \ref{ram-cov32}(i), we may replace $X'_0$ by the maximal open subscheme of $X'$ which is étale over $X$.
So we may assume that $U'_0=U'\subset X'_0$ and $V_0=V$.  
We put $V'=V\times_UU'$ and consider it as a $G$ torsor over $U'$. We have a $U'$-isomorphism 
\begin{equation}\label{ram-cov33a}
\coprod_{I\backslash G}V\stackrel{\sim}{\rightarrow} V'.
\end{equation}
The action of $G$ on $V'$ induces an action on $\coprod_{I\backslash G}V$
defined, for $g,g'\in G$ and $x\in V$, by 
\begin{equation}
g(x_{Ig'})=(g(x))_{Ig'g^{-1}}.
\end{equation}
Let $\Delta\colon G\rightarrow G\times G$
be the diagonal homomorphism, $W$ the quotient of $V\times_kV$ by $\Delta(G)$, $W_1$ the quotient 
of $V\times_kV$ by $\Delta(I)$ and $W'$ the quotient of $V'\times_kV'$ by $\Delta(G)$. We denote by 
$\varpi\colon V\times_kV\rightarrow W_1$ the canonical morphism and by $\varepsilon \colon U\rightarrow W$, 
$\varepsilon_1 \colon U'\rightarrow W_1$
and $\varepsilon' \colon U'\rightarrow W'$ the sections induced by the diagonal morphisms 
$V\rightarrow V\times_kV$ and $V'\rightarrow V'\times_kV'$. 

 The isomorphism \eqref{ram-cov33a} induces an isomorphism 
\begin{equation}\label{ram-cov33b}
\coprod_{I\backslash G \times I\backslash G}V\times_kV\stackrel{\sim}{\rightarrow} V'\times_kV'.
\end{equation}
We denote also by $\Delta\colon I\backslash G\rightarrow I\backslash G\times I\backslash G$
the diagonal map. Then $\coprod_{\Delta(I\backslash G)}V\times_kV$ is an open subscheme of $V'\times_kV'$
stable under the action of $\Delta(G)$. Moreover, the diagonal morphism $V'\rightarrow V'\times_kV'$ 
is induced by the disjoint sum over $\Delta(I\backslash G)$ of the diagonal morphisms $V\rightarrow V\times_kV$. 
On the other hand, the morphism
\begin{equation}
\coprod_{\Delta(I\backslash G)}V\times_kV \rightarrow W_1
\end{equation}
sending $((x,y)_{Ig})$ to $\varpi((g(x),g(y)))$ makes $\coprod_{\Delta(I\backslash G)}V\times_kV$ as a $G$-torsor
over $W_1$. It follows that $W_1$ is an open and closed subscheme of $W'$ and that $\varepsilon'$ is induced by 
$\varepsilon_1$.
Therefore, the ramification of $V'/U'$ at $x'$ is bounded by $R'_0+$ if and only if the ramification 
of $V/U'$ at $x'$ is bounded by $R'_0+$. On the other hand, by \ref{ram-cov32}(ii), the ramification of $V/U$ 
at $x$ is bounded by $R+$ 
if and only if the ramification of $V'/U'$ at $x'$ is bounded by $R'_0+$, hence the proposition.

\begin{rema}\label{ram-cov34}
Let $R$ be an effective rational divisor on $X$ with support in $D$. Then there exists a log-smooth morphism 
of snc-pairs $f\colon (X',D')\rightarrow (X,D)$ such that the underlying morphism of schemes $X'\rightarrow X$
is faithfully flat and that $f^*(R)$ has integral coefficients. Indeed, let $x\in X$, $D_1,\dots,D_n$ be the irreducible 
components of $D$ containing $x$. The question being local on $X$, we may assume that $D=\bigcup_{1\leq i\leq n}D_i$
and that for each $1\leq i\leq n$, $D_i$ is defined by an equation $t_i\in \Gamma(X,\co_X)$. 
We write $R=\sum_{i=1}^nr_iD_i$. For each $1\leq i\leq n$, let $a_i$ be an integer $\geq 1$ such that $a_ir_i$
is an integer. Then 
\begin{equation}
X'=\frac{X[S_i,U_i^{\pm1}; 1\leq i\leq n]}{(U_iS_i^{a_i}-t_i; 1\leq i\leq n)}
\end{equation}
equipped with the normal crossing divisor $D'$ defined by $\prod_{i=1}^nS_i$ answers the question. 
\end{rema}

\subsection{}\label{ram-cov10}
Let $V$ be a Galois torsor over $U$ of group $G$ and let $\Delta\colon G\rightarrow G\times G$  
be the diagonal homomorphism. We denote by $W$ the quotient of $V\times_kV$ 
by the group $\Delta(G)$. The diagonal morphism $\delta_V\colon V\rightarrow V\times_kV$
induces a morphism $\varepsilon_U\colon U\rightarrow W$ above the diagonal morphism 
$\delta_U\colon U\rightarrow U\times_kU$. 

We claim that the quotient of $V\times_kV\times_kV$ 
by the diagonal action of $G$ is canonically isomorphic to $W\times_UW$. 
Indeed, the quotient of $V\times_kV\times_kV$ by $\Delta(G)\times G$
(resp. $G\times\Delta(G)$) is $W\times_kU$ (resp. $U\times_kW$).
\[
\xymatrix{
&{V\times_kV\times_kV}\ar[rd]^{G\times\Delta(G)}\ar[ld]_{\Delta(G)\times G}&\\
{W\times_kU}\ar[rd]&&{U\times_kW}\ar[ld]\\
&{U\times_kU\times_kU}&}
\]
Since $(\Delta(G)\times G)\cap (G\times\Delta(G))$ is the image of the diagonal homomorphism
$G\rightarrow G\times G\times G$, we deduce that the quotient of $V\times_kV\times_kV$ 
by the diagonal action of $G$ 
is canonically isomorphic to $(W\times_kU)\times_{(U\times_kU\times_kU)}(U\times_kW)$, and hence to $W\times_UW$. 

The morphism $\varpi_{13}\colon V\times_kV\times_kV\rightarrow V\times_kV$
defined by $\varpi_{13}(x_1,x_2,x_3)=(x_1,x_3)$ is equivariant for the diagonal actions of $G$ on both sides.
Taking quotients, we obtain a morphism 
\begin{equation}\label{ram-cov10a}
W\times_UW\rightarrow W
\end{equation}
that fits into the commutative diagram 
\begin{equation}\label{ram-cov10b}
\xymatrix{
{W\times_UW}\ar[rr]\ar[d]&& W\ar[d]\\
{(U\times_kU)\times_U(U\times_kU)}\ar@{=}[r]&{U\times_kU\times_kU}\ar[r]^-(0.5){\pr_{13}}&{U\times_kU}}
\end{equation}
where the morphism $\pr_{13}$ is defined by $\pr_{13}(x_1,x_2,x_3)=(x_1,x_3)$. 
If we consider $W$ as the $G$-torsor of isomorphisms of $G$-torsors
from $U\times_kV$ to $V\times_kU$ over $U\times_kU$, then the morphism \eqref{ram-cov10a} 
is the composition of isomorphisms. 

Let $R$ be an effective divisor on $X$ with support in $D$. 
We denote by $Z$ the integral closure of $(X\Asterisk_kX)^{(R)}$
in $W$, by $\pi\colon Z\rightarrow (X\Asterisk_kX)^{(R)}$ the canonical morphism, by $\varepsilon\colon X\rightarrow Z$
the morphism induced by $\varepsilon_U$ (cf. \ref{ram-cov2}) and (abusively) by
$\pr_1,\pr_2\colon Z\rightrightarrows X$ the morphisms induced by the  canonical projections $\pr_1$ 
and $\pr_2$ of $(X\Asterisk_kX)^{(R)}$.
We put $\mX=(X\Asterisk_kX)^{(R)}\times_X(X\Asterisk_kX)^{(R)}$ and denote by 
$\mY$ be the integral closure of $\mX$ in $W\times_UW$. 
Recall \eqref{log8d} that there is a canonical morphism $\mu\colon \mX\rightarrow (X\Asterisk_kX)^{(R)}$ extending 
$\pr_{13}$. Then diagram \eqref{ram-cov10b} induces a 
morphism $\nu\colon \mY\rightarrow Z$ that fits into the following commutative diagram. 
\begin{equation}\label{ram-cov10c}
\xymatrix{
{\mY}\ar[r]^-(0.5)\nu\ar[d]&{Z}\ar[d]\\
{\mX}\ar[r]^-(0.5)\mu&{(X\Asterisk_kX)^{(R)}}}
\end{equation}

Let $Z_0$ be the maximal open subscheme of $Z$ which is étale over $(X\Asterisk_kX)^{(R)}$. 
Observe that $\mY$ is the integral closure of $Z\times_XZ$ in $W\times_UW$. 
The canonical projections $\pr_1,\pr_2\colon (X\Asterisk_kX)^{(R)}\rightarrow X$ are smooth.
Then $\mX$ is smooth over $X$, and hence regular. 
Since $Z_0\times_{X}Z_0$ is étale over $\mX$, it is regular. 
Therefore, we can identify $Z_0\times_XZ_0$ with an open subscheme of $\mY$. We claim that
\begin{equation}\label{ram-cov10d}
Z_0\times_{X}Z_0\subset \nu^{-1}(Z_0).
\end{equation}
We denote by $\mu^*(Z)$ (resp. $\mu^*(Z_0)$) the base change of $Z$ (resp. $Z_0$) by $\mu$. 
Then $\nu$ induces  a finite $\mX$-morphism $\mY\rightarrow \mu^*(Z)$.
Since $\mu$ is smooth \eqref{log9}, $\mu^*(Z_0)$ is the maximal open subscheme of $\mu^*(Z)$
which is étale over $\mX$. Since $Z_0\times_XZ_0$ is étale
over $\mX$, it is unramified over $\mu^*(Z)$. Then by (\cite{ega4} 18.10.1),
$Z_0\times_{X}Z_0$ is étale over $\mu^*(Z)$; in particular, it is flat over $\mu^*(Z)$. 
Hence the inclusion \eqref{ram-cov10d} follows from (\cite{ega4} 17.7.7). 
By \eqref{ram-cov10d}, $\nu$ induces a morphism that we denote also by 
\begin{equation}\label{ram-cov10e}
\nu\colon Z_0\times_XZ_0\rightarrow Z_0.
\end{equation}

The automorphism $i$ of $V\times_kV$ switching the factors is equivariant for the diagonal action of $G$. 
Taking quotients, we obtain an automorphism $\iota_U$ of $W$ that lifts the automorphism of $U\times_kU$  
switching the factors. Then $\iota_U$ extends to an automorphism $\iota$ of $Z$ that fits into the commutative
diagram 
\begin{equation}\label{ram-cov10f}
\xymatrix{
Z\ar[r]^{\iota}\ar[d]_\pi&Z\ar[d]^\pi\\
{(X\Asterisk_k X)^{(R)}}\ar[r]^{\sigma}&{(X\Asterisk_k X)^{(R)}}}
\end{equation}
where $\sigma$ is the automorphism \eqref{log7d}. It is clear that $\iota(Z_0)=Z_0$; we denote 
also by $\iota$ the automorphism of $Z_0$ induced by $\iota$. Let
\begin{equation}
\alpha\colon Z_0\times_XZ_0\rightarrow Z_0\times_XZ_0
\end{equation}
be the morphism defined by $\alpha(x,y)=(\iota(y),\iota(x))$. It is well defined because of the following commutative 
diagram.
\begin{equation}
\xymatrix{
{Z_0\times_XZ_0}\ar[r]^-(0.5){\pr_2}\ar[d]_{\pr_1}&{Z_0}\ar[d]^{\pr_1}\ar[r]^-(0.5)\iota&{Z_0}\ar[dd]^{\pr_2}\\
{Z_0}\ar[r]^{\pr_2}\ar[d]_{\iota}&X\ar@{=}[rd]&\\
{Z_0}\ar[rr]^{\pr_1}&&X}
\end{equation}

\begin{lem}\label{ram-cov11}
We keep the assumptions of \eqref{ram-cov10}. 

{\rm (i)}\ The diagrams 
\begin{equation}\label{ram-cov11a}
\xymatrix{
{Z_0\times_XZ_0\times_XZ_0}\ar[rr]^-(0.5){\id\times \nu}\ar[d]_{\nu\times \id}&&{Z_0\times_XZ_0}\ar[d]^{\nu}\\
{Z_0\times_XZ_0}\ar[rr]^{\nu}&&{Z_0}}
\end{equation}
\begin{equation}\label{ram-cov11aa}
\xymatrix{
{Z_0\times_XZ_0}\ar[r]^\alpha\ar[d]_{\nu}&{Z_0\times_XZ_0}\ar[d]^\nu\\
Z_0\ar[r]^{\iota}&Z_0}
\end{equation}
are commutative.

{\rm (ii)}\ Assume moreover that the ramification of $V/U$ along $D$ is bounded by $R+$;
so we have $\varepsilon(X)\subset Z_0$. We denote also by $\varepsilon\colon X\rightarrow Z_0$ the morphism 
induced by $\varepsilon$.  Then the diagrams
\begin{equation}\label{ram-cov11b}
\xymatrix{
{Z_0}\ar[rr]^-(0.5){\id\times \varepsilon}\ar[d]_{\varepsilon\times \id}\ar[rrd]|{\id}&&{Z_0\times_XZ_0}\ar[d]^\nu\\
{Z_0\times_XZ_0}\ar[rr]^\nu&&{Z_0}}
\end{equation}
\begin{equation}\label{ram-cov11bb}
\xymatrix{
Z_0\ar[r]^{\pr_1}\ar[d]_{\id\times \iota}&X\ar[d]^{\varepsilon}&Z_0\ar[l]_{\pr_2}\ar[d]^{\iota\times \id}\\
{Z_0\times_XZ_0}\ar[r]^-(0.5)\nu&Z_0&{Z_0\times_XZ_0}\ar[l]_-(0.5){\nu}}
\end{equation}
are commutative.
\end{lem}

(i) The diagram 
\begin{equation}\label{ram-cov11c}
\xymatrix{
{V\times_kV\times_kV\times_kV}\ar[r]^-(0.5){\varpi_{124}}\ar[d]_-(0.5){\varpi_{134}}&
{V\times_kV\times_kV}\ar[d]^{\varpi_{13}}\\
{V\times_kV\times_kV}\ar[r]^{\varpi_{13}}&{V\times_kV}}
\end{equation}
where $\varpi_{1i4}(x_1,x_2,x_3,x_4)=(x_1,x_i,x_4)$ for $i\in \{2,3\}$, is commutative and equivariant 
for the diagonal actions of $G$. It induces by taking quotients the following commutative diagram.  
\begin{equation}\label{ram-cov11d}
\xymatrix{
{W\times_UW\times_UW}\ar[rr]^-(0.5){\id\times \nu_U}\ar[d]_-(0.5){\nu_U\times \id}&&
{W\times_UW}\ar[d]^{\nu_U}\\
{W\times_UW}\ar[rr]^{\nu_U}&&{W}}
\end{equation}
Indeed, if we denote by $\Delta\colon G\rightarrow G\times G$ and
$\nabla\colon G\rightarrow G\times G\times G$ the diagonal homomorphisms, then the 
quotient of $V\times_kV\times_kV\times_kV$ by $\Delta(G)\times G\times G$ (resp. by $G\times \nabla(G)$) 
is $W\times_kU\times_kU$ (resp. $U\times_k W\times_UW$).
\begin{equation}\label{ram-cov11e}
\xymatrix{
&{V\times_kV\times_kV\times_kV}\ar[ld]_{\Delta(G)\times G\times G}\ar[rd]^{G\times \nabla(G)}&\\
{W\times_kU\times_kU}\ar[rd]&&{U\times_k W\times_UW}\ar[ld]\\
&{U\times_kU\times_kU\times_kU}&}
\end{equation}
Since $(\Delta(G)\times G\times G)\cap (G\times\nabla(G))$ is the image of the diagonal homomorphism
$G\rightarrow G\times G\times G\times G$, it follows that $\varpi_{124}$ induces by quotient by the diagonal
actions of $G$ the morphism 
\[
\id\times\nu_U\colon W\times_UW\times_UW\rightarrow W\times_UW.
\] 
By switching the second and the third factors, we prove that 
$\varpi_{134}$ induces by quotient by the diagonal actions of $G$ the morphism 
$\nu_U\times \id\colon W\times_UW\times_UW\rightarrow W\times_UW$.
Therefore, the diagram \eqref{ram-cov11d} is commutative, and hence so is the diagram \eqref{ram-cov11a}. 

The diagram 
\begin{equation}
\xymatrix{
{V\times_kV\times_kV}\ar[r]^-(0.5)\beta\ar[d]_{\varpi_{13}}&{V\times_kV\times_kV}\ar[d]^{\varpi_{13}}\\
{V\times_kV}\ar[r]^-(0.5)i&{V\times_kV}}
\end{equation}
where $\beta(x,y,z)=(z,y,x)$ and $i(x,y)=(y,x)$, is commutative and equivariant 
for the diagonal actions of $G$. It induces by taking quotients the following commutative diagram.
\begin{equation}
\xymatrix{
{W\times_UW}\ar[r]^-(0.5){\alpha_U}\ar[d]_{\nu_U}&{W\times_UW}\ar[d]^{\nu_U}\\
{W}\ar[r]^-(0.5){\iota_U}&{W}}
\end{equation}
The latter proves that the diagram \eqref{ram-cov11aa} is commutative.

(ii) The diagram 
\begin{equation}\label{ram-cov11f}
\xymatrix{
{V\times_kV}\ar[r]^-(0.5){\id\times \delta_V}\ar[d]_-(0.5){\delta_V\times \id}\ar[rd]|{\id}&
{V\times_kV\times_kV}\ar[d]^{\varpi_{13}}\\
{V\times_kV\times_kV}\ar[r]^{\varpi_{13}}&{V\times_kV}}
\end{equation}
is commutative and equivariant 
for the diagonal actions of $G$. It induces by taking quotients the following commutative diagram.  
\begin{equation}\label{ram-cov11g}
\xymatrix{
{W}\ar[rr]^-(0.5){\id\times \varepsilon_U}\ar[d]_-(0.5){\varepsilon_U\times \id}\ar[rrd]|{\id}&&
{W\times_UW}\ar[d]^{\nu_U}\\
{W\times_UW}\ar[rr]^{\nu_U}&&{W}}
\end{equation}
The latter proves that the diagram \eqref{ram-cov11b} is commutative. 

The diagram 
\begin{equation}
\xymatrix{
{V\times_kV}\ar[r]^{\pr_1}\ar[d]_{\alpha}&V\ar[d]^{\delta_V}\\
{V\times_kV\times_kV}\ar[r]^-(0.5){\varpi_{13}}&{V\times_kV}}
\end{equation}
where $\alpha(x,y)=(x,y,x)$ is commutative and equivariant 
for the diagonal actions of $G$. It induces by taking quotients the following commutative diagram.
\begin{equation}
\xymatrix{
W\ar[r]^{\pr_1}\ar[d]_{\id\times \iota_U}&U\ar[d]^{\varepsilon_U}\\
{W\times_UW}\ar[r]^-(0.5){\nu_U}&{W}}
\end{equation}
The latter proves that the left square in the diagram \eqref{ram-cov11bb} is commutative. The same argument
shows that the right square is also commutative.

\begin{prop}\label{ram-cov41}
We keep the assumptions of \eqref{ram-cov10} (so $R$ has integral coefficients) and 
assume moreover that the ramification of $V/U$ along $D$ is bounded by $R+$. We put $F=Z_0\times_XR$ and
$E=(X\Asterisk_k X)^{(R)}\times_XR$, which is a vector bundle over $R$ \eqref{log7c}. Then~:

{\rm (i)}\ $F$ is a group scheme over $R$, and the morphism 
$\pi_R\colon F\rightarrow E$ induced by $\pi$
is a surjective étale morphism of group schemes over $R$. 

{\rm (ii)}\ For every geometric point $\ox$ of $R$, the neutral connected component $F^\circ_{\ox}$ of $F_\ox$
is isomorphic to a product of additive groups over $\ox$, the morphism $\pi_{\ox}\colon  F^\circ_{\ox}
\rightarrow E_{\ox}$ is finite étale and surjective and its kernel is an $\mF_p$-vector space of finite dimension. 
\end{prop}
 
(i) The closed subscheme $R\times_X(X\Asterisk_kX)^{(R)}$ of $(X\Asterisk_kX)^{(R)}$ is equal to $E$, 
and the canonical projections  $\pr_1,\pr_2\colon(X\Asterisk_kX)^{(R)}\rightrightarrows X$ 
induce the same morphism $E\rightarrow R$. 
Hence, the closed subscheme $R\times_XZ_0$ of $Z_0$ is equal to $F$, and the canonical morphisms 
$\pr_1,\pr_2\colon Z_0\rightrightarrows X$ induce the same morphism $F\rightarrow R$. Then it follows from \ref{ram-cov11} 
that $F$ is a group scheme over $R$. We deduce from \ref{log9} and the commutative diagram 
\begin{equation}
\xymatrix{
{Z_0\times_XZ_0}\ar[r]^-(0.5)\nu\ar[d]_{\pi\times \pi}&{Z_0}\ar[d]^{\pi}\\
{(X\Asterisk_kX)^{(R)}\times_X(X\Asterisk_kX)^{(R)}}\ar[r]^-(0.5)\mu&{(X\Asterisk_kX)^{(R)}}}
\end{equation}
that $\pi_R$ is a morphism of group schemes over $R$. By the definition of $Z_0$, $\pi_R$ is étale,
and it follows from (ii) that it is surjective. 

(ii) It follows from \ref{ram-cov9}.

\begin{cor}\label{ram-cov4}
We keep the assumptions of \eqref{ram-cov2} and assume moreover that $R$ has integral coefficients.
Then the following conditions are equivalent~:

{\rm (i)}\ The ramification of $V/U$ along $D$ is bounded by $R+$.

{\rm (ii)}\ There exists an open neighbourhood $Z_0$ of $\varepsilon(X)$ in $Z$ which is étale over $(X\Asterisk_kX)^{(R)}$
and such that $\pi(Z_0)$ contains $E=(X\Asterisk_kX)^{(R)}\times_XR$.
\end{cor}

\begin{rema}
We can deduce \ref{ram-cov4} from \eqref{ram-cov10d} by a shorter argument using only \ref{log9}. 
\end{rema}

\subsection{}\label{ramcodim1}
Let $R$ be an effective divisor on $X$ with support in $D$,
$\xi$ a generic point of $D$, $\oxi$ a geometric point of $X$ above $\xi$, 
$S$ the strict localization of $X$ at $\oxi$, $\eta$ the generic point of $S$,
$\eta'$ an integral finite étale extension of $\eta$ and $S'$ the integral closure of $S$ in $\eta'$. We put 
\begin{equation}\label{ramcodim1a}
(X\Asterisk_kS')^{(R)}=(X\Asterisk_kX)^{(R)}\times_XS'.
\end{equation}
(This notation could be justified by \ref{log52}). The morphism $\delta^{(R)}$ induces a section 
$\delta^{(R)}_{S'}\colon S'\rightarrow  (X\Asterisk_kS')^{(R)}$ of the canonical projection 
$(X\Asterisk_kS')^{(R)}\rightarrow S'$.

Let  $V$ be a Galois torsor over $U$ of group $G$, 
$Y$ the integral closure of $X$ in $V$, $f\colon Y\rightarrow X$ the canonical morphism, 
$\Delta\colon G\rightarrow G\times G$ the diagonal homomorphism and
$W$ the quotient of $V\times_kV$ by the group $\Delta(G)$. 
The diagonal morphism $\delta_V\colon V\rightarrow V\times_kV$
induces a morphism $\varepsilon_U\colon U\rightarrow W$ above the diagonal morphism 
$\delta_U\colon U\rightarrow U\times_kU$.
We denote by $Z$ the integral closure of $(X\Asterisk_kX)^{(R)}$ in $W$,
by $\pi\colon Z\rightarrow (X\times_kX)^{(R)}$ the canonical morphism and by $\varepsilon\colon X\rightarrow Z$
the morphism induced by $\varepsilon_U\colon U\rightarrow W$. We have $\pi\circ \varepsilon=\delta^{(R)}$
(cf. \eqref{ram-cov3a}). Let $\tZ_{S'}$ be the normalization of $Z_{S'}=Z\times_XS'$, or equivalently, the integral closure
of $(X\Asterisk_kS')^{(R)}$ in $W\times_U\eta'$.  The morphism $\varepsilon$ induces a section 
$\varepsilon_{S'}\colon S'\rightarrow \tZ_{S'}$ of the canonical morphism $\tZ_{S'}\rightarrow S'$, that lifts
$\delta^{(R)}_{S'}$. Observe that $\tZ_S=Z_S$. 

We denote by $Q_{S'}$ the normalization of $Y\times_X(X\Asterisk_kS')^{(R)}$, or equivalently,
the integral closure of $(X\Asterisk_kS')^{(R)}$ in $V\times_k\eta'$. If the canonical morphism 
$S'\rightarrow X$ factors as $S'\rightarrow Y\stackrel{f}{\rightarrow}X$, then we have canonical 
isomorphisms 
\begin{equation}\label{ramcodim1b}
V\times_k\eta'\simeq (V\times_kV)\times_V\eta'\simeq (W\times_UV)\times_V\eta'\simeq W\times_U\eta'.
\end{equation}
Hence, we obtain an isomorphism
\begin{equation}\label{ramcodim1c}
Q_{S'}\stackrel{\sim}{\rightarrow} \tZ_{S'}.
\end{equation}

\begin{lem}\label{ramcodim2}
We keep the assumptions of \eqref{ramcodim1} and let 
$\sigma\colon S'\rightarrow Q_{S'}$ be a section of the canonical morphism $Q_{S'}\rightarrow S'$
lifting the section $\delta^{(R)}_{S'}\colon S'\rightarrow (X\Asterisk_kS')^{(R)}$. 
If the canonical morphism $Q_{S'}\rightarrow (X\Asterisk_k S')^{(R)}$ is étale on an open neighborhood of
$\sigma$, then it is étale everywhere. 
\end{lem}

Observe first that $ (X\Asterisk_k S')^{(R)}$ is regular since it is smooth over $S'$. 
The group $G$ acts on $Q_{S'}$, and the quotient of $Q_{S'}$ by $G$ is $(X\Asterisk_kS')^{(R)}$. 
Therefore, the morphism $Q_{S'}\rightarrow (X\Asterisk_k S')^{(R)}$ is étale above an open neighborhood of
$\delta^{(R)}_{S'}(S')$. Then the assertion follows from Zariski-Nagata's purity theorem (\cite{sga2} X 3.4).

\begin{lem}\label{ramcodim3}
Under the assumptions of \eqref{ramcodim1}, the following conditions are equivalent~:

{\rm (i)}\ The canonical morphism $Z_{S}\rightarrow (X\Asterisk_kS)^{(R)}$ is étale on an open 
neighborhood of $\varepsilon_{S}(S)$. 

{\rm (ii)}\ For any integral finite étale extension $\eta'$ of $\eta$, 
the canonical morphism $\tZ_{S'}\rightarrow (X\Asterisk_kS')^{(R)}$ is étale on an open 
neighborhood of $\varepsilon_{S'}(S')$. 

{\rm (iii)}\ There exists an integral finite étale extension $\eta'$ of $\eta$ such that 
the canonical morphism $\tZ_{S'}\rightarrow (X\Asterisk_kS')^{(R)}$ is étale on an open 
neighborhood of $\varepsilon_{S'}(S')$. 

{\rm (iv)}\ For any connected component $T$ of $Y\times_XS$, 
the canonical morphism $Q_{T}\rightarrow (X\Asterisk_kT)^{(R)}$ is étale. 

{\rm (v)}\ There exists an integral finite étale extension $\eta'$ of $\eta$ such that 
the canonical morphism $Q_{S'}\rightarrow (X\Asterisk_kS')^{(R)}$ is étale. 
\end{lem}

Conditions (i), (ii) and (iii) are equivalent by \ref{rc2}. 
We have (ii)$\Rightarrow$(iv) by \eqref{ramcodim1c} and \ref{ramcodim2}. 
It is clear that we have (iv)$\Rightarrow$(v). If the condition (v) holds true, then it holds under the extra 
assumption that the canonical morphism $S'\rightarrow X$ factors as $S'\rightarrow Y\stackrel{f}{\rightarrow}X$. 
Then we have (v)$\Rightarrow$(iii) by \eqref{ramcodim1c}.

\begin{prop}\label{ram-cov8}
Let $V$ be a Galois torsor over $U$, $R$ a rational divisor 
on $X$ with support in $D$, $\xi$ a generic point of $D$, 
$\oxi$ a generic point of $X$ above $\xi$, $S$ the strict localization of $X$ at $\oxi$,
$K$ the fraction field of $\Gamma(S,\co_S)$ and $r$ the multiplicity of $R$ at $\xi$. 
We put $V\times_U\Spec(K)=\Spec(L)$, 
where $L=\prod_{i=1}^nL_i$ is a finite product of finite separable extensions of $K$. 
Then the following conditions are equivalent~:

{\rm (i)}\ The ramification of $V/U$ at $\xi$ is bounded by $R+$.

{\rm (ii)}\ For every $1\leq i\leq n$, the logarithmic ramification of $L_i/K$ is bounded by $r+$ \eqref{rrt15}. 
\end{prop}

By \ref{ram-cov32}(ii), \ref{ram-cov34} and (\cite{as2} 5.2), we may assume that $R$ has integral coefficients. 
We take again the notation of \ref{ramcodim1}. The condition (i) is equivalent to condition \ref{ramcodim3}(i). 
Condition (ii) is equivalent to condition \ref{ramcodim3}(v) by (\cite{saito1} 1.13 and the remark after its proof).  
Hence, the proposition follows from \ref{ramcodim3}. 

\begin{prop}\label{ram-cov6}
Let $V$ be a Galois torsor over $U$ of group $G$, $Y$ the integral closure of $X$ in $V$ and
$R$ an effective rational divisor on $X$ with support in $D$. Assume that the following conditions are satisfied~:
\begin{itemize}
\item[{\rm (i)}] $V/U$ has the property {\rm (NpS)} at every geometric point of $X$ \eqref{NS1}, 
that is, for every geometric point $\oy$ of $Y$, 
the inertia group $I_\oy\subset G$ of $\oy$ has a normal $p$-Sylow subgroup~;
\item[{\rm (ii)}] for every generic point $\xi$ of $D$, the ramification of $V/U$ at $\xi$ is bounded by $R+$.
\end{itemize}
Then the ramification of $V/U$ along $D$ is bounded by $R+$.
\end{prop}

It is enough to prove that for every $x\in X$, 
the ramification of $V/U$ at $x$ is bounded by $R+$. Let $\oy$ be a geometric point of $Y$ 
localized at a point $y\in Y$ above $x$. We denote by $U'$ the quotient of $V$ 
by the inertia group $I_\oy$ of $\oy$, 
by $X'$ the integral closure of $X$ in $U'$, by $f\colon X'\rightarrow X$
the structural morphism, by $X'_0$ the maximal open subscheme of $X'$ which is étale over $X$
and by $D'_0$ and $R'_0$ the pull-backs of $D$ and $R$ over $X'_0$; so $(X'_0,D'_0)$ is an snc-pair over $k$
and $U'\subset X'_0$. Since $x\in f(X'_0)$, it is enough  to prove that the ramification of $V/U'$
along $D'_0$ is bounded by $R'_0+$ \eqref{ram-cov33}. 
Replacing  $V/U$ and $(X,D)$ by $V/U'$ and $(X'_0,D'_0)$ (\cite{as1} 3.15)  (cf. \ref{ram-cov61}),
we may assume that $G$ has a normal $p$-Sylow subgroup. 
By \ref{ram-cov32}, \ref{ram-cov34}, \ref{ram-cov8} and (\cite{as1} 3.15), 
we may also assume that $R$ has integral coefficients. 

Assume next that $G$ has a normal $p$-Sylow subgroup $H$ and that $R$ has integral coefficients. 
It is enough to prove that for every $x\in X$, the ramification of $V/U$ at $x$ is bounded by $R+$. 
Let $U_1$ be the quotient of $V$ by $H$; so $U_1$ is a Galois torsor over $U$
which is tamely ramified along $D$. By Abhyankar's lemma (\cite{sga1} XIII 5.2 and 5.4), 
there exists a morphism of snc-pairs $h\colon (X',D')\rightarrow (X,D)$ satisfying the following properties~:

(P)\ The morphism $h$ is log-smooth, the morphism $X'\rightarrow X$ 
is flat and $x\in h(X')$. 

(Q)\ We put $U'=X'-D'$ and $U'_1=U_1\times_UU'$, and denote by $X'_1$
the integral closure of $X'$ in $U'_1$~; then $X'_1$ is étale over $X'$.

We put $V'=V\times_UU'$ and denote by $D'_1$ and $R'_1$ the pull-backs of $D$ and $R$ over $X'_1$. 
By  \ref{ram-cov32} and \ref{ram-cov33}, it is enough to prove that the ramification
of $V'/U'_1$ along $D'_1$ is bounded by $R'_1$. Hence, we are reduced to the case where $G$ is nilpotent
and $R$ has integral coefficients  (\cite{as1} 3.15) (cf. \ref{ram-cov61}).

Assume finally that $G$ is nilpotent and that $R$ has integral coefficients.
Let $\Delta\colon G\rightarrow G\times G$ be the diagonal homomorphism, 
$W$ the quotient of $V\times_kV$ by $\Delta(G)$ and $\varepsilon_U\colon U\rightarrow W$
the morphism induced by the diagonal morphism $\delta_V\colon V\rightarrow V\times_kV$.
We denote by $Z$ the integral closure of $(X\Asterisk_kX)^{(R)}$ in $W$
and by $\varepsilon\colon X\rightarrow Z$ the morphism induced by $\varepsilon_U$. 
Since $(X\Asterisk_kX)^{(R)}$ is smooth over $X$, it is regular \eqref{tub26}. 
Then the proposition follows from \ref{ram-cov5} applied 
to the open subscheme $U\times_kU$ of $(X\Asterisk_kX)^{(R)}$,
to the finite étale covering $W$ of $U\times_kU$ and to the closed subscheme $\varepsilon(X)$ of $Z$.

\begin{rema} \label{ram-cov61}
We keep the assumptions of \ref{ram-cov6} and consider the following conditions 

{\rm (i')}\ $G$ has a normal $p$-Sylow subgroup. 

{\rm (i'')}\ $G$ is nilpotent. 

Then we have (i'')$\Rightarrow$(i')$\Rightarrow$(i). Indeed, the first implication follows from (\cite{alg} chap.~I § 6.7 theo.~4)
and the second is a consequence of (\cite{alg} chap.~I § 6.6 cor.~3 of theo.~3).
\end{rema}

\begin{defi}\label{ram-cov17}
Let $V$ be a Galois torsor over $U$ of group $G$. We define the {\em conductor
of $V/U$ relatively to $X$} to be the minimum effective rational divisor $R$ on $X$ with support in 
$D$ such that for every generic point $\xi$ of $D$, the ramification of $V/U$ at $\xi$ is bounded by $R+$. 
\end{defi}

This terminology may be slightly misleading as the ramification of $V/U$ along $D$ may not be bounded by $R+$
in general. However, we have the following~:

\begin{prop}\label{ram-cov18}
Let $V$ be a Galois torsor over $U$ of group $G$. Assume that the following strong form of resolution
of singularities holds~:

\begin{itemize}
\item[{\rm (RS)}] For any $U$-admissible blow-up $Y$ of $X$, there exists an snc-pair $(Y',E')$ over $k$,
a morphism of pairs $(Y',E')\rightarrow (X,D)$ and a proper $X$-morphism $Y'\rightarrow Y$
inducing an isomorphism $Y'-E'\stackrel{\sim}{\rightarrow}U$. 
\end{itemize}

Then, there exists an snc-pair $(X',D')$ over $k$ and a proper morphism of snc-pairs $f\colon (X',D')\rightarrow (X,D)$
inducing an isomorphism $X'-D'\stackrel{\sim}{\rightarrow} U$ such that  
if we denote by $R'$ the conductor of $V/U$ relatively to $X'$,
the ramification of $V/U$ along $D'$ is bounded by $R'+$.
\end{prop}

By \ref{NS4}, there exists a $U$-admissible blow-up $Y\rightarrow X$ such that if we denote by 
$Y'$ the normalization of $Y$, $V/U$ has the property (NpS) at every geometric point of $Y'$.
By assumption (RS), there exists an snc-pair $(X',D')$ over $k$,
a morphism of pairs $(X',D')\rightarrow (X,D)$ and a proper $X$-morphism $X'\rightarrow Y$
inducing an isomorphism $X'-D'\stackrel{\sim}{\rightarrow}U$. 
Then $V/U$ has the property (NpS) at every geometric point of $X'$ 
\eqref{NS25}. Let $R'$ the conductor of $V/U$ relatively to $X'$. 
It follows from \ref{ram-cov6} that the ramification of $V/U$ along $D'$ is bounded by $R'+$.

\section{Ramification of $\ell$-adic sheaves}\label{ram}

\subsection{}\label{ram1}
In this section, we fix an snc-pair $(X,D)$ over $k$, and put $U=X-D$. 
We denote by $j\colon U\rightarrow X$ the canonical injection, 
by $X\Asterisk_kX$ the framed self-product of $(X,D)$ 
and by $\delta\colon X\rightarrow X\Asterisk_kX$ the framed diagonal map \eqref{fram17}. 
We will no longer consider the logarithmic structure of $X\Asterisk_kX$;
only the underlying scheme will be of interest for us. 
We consider $X\Asterisk_kX$ as an $X$-scheme by the second projection. 
For any effective rational divisor $R$ on $X$ with support in $D$, we denote by $(X\Asterisk_kX)^{(R)}$
the dilatation of $X\Asterisk_kX$ along $\delta$ of thickening $R$ \eqref{log7}, by  
\begin{equation}\label{ram1a}
\delta^{(R)}\colon X\rightarrow (X\Asterisk_kX)^{(R)}
\end{equation}
the canonical lifting of $\delta$, and by 
\begin{equation}\label{ram1b}
j^{(R)}\colon U\times_kU\rightarrow (X\Asterisk_kX)^{(R)}
\end{equation}
the canonical open immersion. Then we have the following Cartesian diagram.
\begin{equation}\label{ram1c}
\xymatrix{
U\ar[r]^-(0.5){\delta_U}\ar[d]_j&{U\times_kU}\ar[d]^-(0.5){j^{(R)}}\\
X\ar[r]^-(0.5){\delta^{(R)}}&{(X\Asterisk_kX)^{(R)}}}
\end{equation}

\begin{prop}\label{ram2}
Let $\cF$ be a locally constant constructible sheaf of $\Lambda$-modules on $U$, 
$R$ an effective rational divisor on $X$ with support in $D$, $x\in X$, $\ox$ a 
geometric point of $X$ above $x$, and $\cH(\cF)$ the sheaf on $U\times_kU$ defined in \eqref{not2a}.
Then the base change morphism
\begin{equation}\label{ram2a}
\alpha\colon \delta^{(R)*}j^{(R)}_*(\cH(\cF))\rightarrow j_*\delta_U^*(\cH(\cF))=j_*(\cEnd(\cF))
\end{equation}
relatively to the Cartesian diagram \eqref{ram1c} is injective. 
Furthermore, the following three conditions are equivalent~:

{\rm (i)}\ The stalk 
\begin{equation}\label{ram2aa}
\alpha_\ox\colon (\delta^{(R)*}j^{(R)}_*(\cH(\cF)))_{\ox}\rightarrow j_*(\cEnd(\cF)))_{\ox}
\end{equation}
of the morphism $\alpha$ at $\ox$ is an isomorphism. 

{\rm (ii)}\ The image of the identity endomorphism of $\cF$ in $j_*(\cEnd(\cF)))_{\ox}$ 
is contained in the image of $\alpha_\ox$ \eqref{ram2aa}. 

{\rm (iii)}\ There exists a Galois torsor $V$ over $U$ trivializing $\cF$ 
such that the ramification of $V/U$ at $x$ is bounded by $R+$.

Assume moreover that there exists a Galois torsor $V_0$ over $U$ of group $G_0$ satisfying the following condition~:

$(\Asterisk)$\ There exists a $\Lambda$-module
$M$ equipped with a faithful representation $\rho$ of $G_0$ and an isomorphism of sheaves (with Galois descent data)
$\cF|V_0\simeq M_{V_0}$, where $M_{V_0}$ is the constant sheaf on $V_0$ of value $M$ equipped with the 
descent datum defined by $\rho$. 

Then conditions {\rm (i)}, {\rm (ii)} and {\rm (iii)} are equivalent to the following condition~:

{\rm (iv)}\ The ramification of $V_0/U$ at $x$ is bounded by $R+$. 
\end{prop}

It follows from \ref{ram00} that $\alpha$ is injective. 
For the second proposition, we may assume that $X$ and hence $U$ are connected;
in particular, we may assume that there exists a Galois torsor $V_0$ over $U$ of group $G_0$ 
satisfying condition $(\Asterisk)$. We clearly have (i)$\Rightarrow$(ii) and (iv)$\Rightarrow$(iii). 

Let $V$ be a Galois torsor over $U$ of group $G$ trivializing $\cF$ and let $Y$ be the integral closure 
of $X$ in $V$. We denote by $\Delta\colon G\rightarrow  G\times G$ the diagonal homomorphism,
by $W$ the quotient of $V\times_kV$ by $\Delta(G)$, by $Z$ (resp. $\Sigma$)
the integral closure of $(X\Asterisk_kX)^{(R)}$ in $W$ (resp. $V\times_kV$) and by  
$j_W\colon W\rightarrow Z$  the canonical injection (cf. \ref{ram-cov2}).  
The diagonal morphism $\delta_V\colon V\rightarrow V\times_kV$ induces a morphism 
$\varepsilon_U\colon U\rightarrow W$ above the diagonal morphism $\delta_U\colon U\rightarrow U\times_kU$. 
Let $\varepsilon\colon X\rightarrow Z$ and $\sigma\colon Y\rightarrow \Sigma$ be 
the morphisms induced by $\varepsilon_U$ and $\delta_V$, respectively.
Then the following diagram is commutative.
\begin{equation}\label{ram2c}
\xymatrix{
V\ar[r]^-(0.5){\delta_V}\ar[d]&{V\times_kV}\ar[d]\ar[r]&\Sigma\ar[d]&Y\ar[l]_\sigma\ar[d]\\
U\ar[r]\ar[rd]_-(0.5){\delta_U}&W\ar[r]^{j_W}\ar[d]&Z\ar[d]&X\ar[l]_{\varepsilon}\ar[ld]^-(0.5){\delta^{(R)}}\\
&{U\times_kU}\ar[r]^-(0.45){j^{(R)}}&{(X\Asterisk_kX)^{(R)}}&}
\end{equation}
Each canonical projection $V\times_kV\rightarrow V$ induces a morphism $\Sigma \rightarrow Y$ 
for which $\sigma$ is a section. 
Let $\oy$ be a geometric point  of $Y$ above $\ox$,   
$I_{\oy}\subset G\times G$ the inertia group of $\sigma(\oy)$,
and $J_{\oy}\subset G$ the inertia group of $\oy$. 
Since $\sigma$ is a closed immersion, we have $J_\oy=\Delta^{-1}(I_\oy)$. 
Hence the following conditions are equivalent~:

(a) The ramification of $V/U$ at $x$ is bounded by $R+$.

(b) We have $I_{\oy}\subset \Delta(G)$.

(c)  We have $I_{\oy}=\Delta(J_{\oy})$.

Let $Y'$ be the connected component of $Y$ containing $\oy$ and  let
$V'=V\cap Y'$. The stabilizer $G'\subset G$ of $V'$ acts on $N=\Gamma(V',\cF)$. We have 
$I_{\oy}\subset G'$, and the morphism $\alpha_{\ox}$ \eqref{ram2aa} is canonically identified 
with the injective morphism
\begin{equation}\label{ram2d}
(\End(N))^{I_\oy}\rightarrow (\End(N))^{\Delta(J_\oy)}. 
\end{equation}
Hence, by the equivalence (a)$\Leftrightarrow$(c) above, we deduce the implication (iii)$\Rightarrow$(i). 

It remains to prove that (ii)$\Rightarrow$(iv). We keep the previous notation and assume moreover that $V=V_0$
and $G=G_0$. Consider the following commutative diagram with Cartesian squares.
\begin{equation}\label{ram2e}
\xymatrix{
U\ar[r]^j\ar[d]_{\varepsilon_U}&X\ar[d]_\varepsilon\ar@/^1pc/[dd]^{\delta^{(R)}}\\
W\ar[r]^{j_W}\ar[d]&Z\ar[d]\\
{U\times_kU}\ar[r]^-(0.45){j^{(R)}}&{(X\Asterisk_kX)^{(R)}}}
\end{equation}
By (\cite{egr1} 1.2.4(i)), the base change morphism $\alpha$ \eqref{ram2a} is composed of 
\begin{equation}
\xymatrix{
{\delta^{(R)*}(j_*^{(R)}(\cH(\cF))))}\ar[r]^{\varepsilon^*{\beta}}
&{\varepsilon^*(j_{W*}(\cH(\cF)|W))}\ar[r]^-(0.4)\gamma&{j_*(\cEnd(\cF))}},
\end{equation}
where $\beta\colon j^{(R)}_*(\cH(\cF))|Z\rightarrow j_{W*}(\cH(\cF)|W)$ and $\gamma$ are the base change 
morphisms relatively to  the lower and the upper squares of \eqref{ram2e} respectively.

If we equip $\End(M)$ with the canonical action of $G\times G$,
we deduce from the isomorphism $\cF|V\simeq M_V$ an isomorphism (of sheaves with Galois descent data)
$\cH(\cF)|(V\times_kV)\simeq \End(M)_{(V\times_kV)}$. 
Since the action of $G$ on $M$ is faithful, 
$\Delta(G)$ is the stabilizer of $\id\in \End(M)$ in $G\times G$. In particular, we may consider 
$\id$ as a section of $\Gamma(Z,j_{W*}(\cH(\cF)|W))$. Since $\gamma$ is injective by \ref{ram00}, 
condition (ii) is equivalent to the following condition~:

(ii') The image of $\id$ in $j_{W*}(\cH(\cF)|W)_{\varepsilon(\ox)}$ 
is contained in the image of $\beta_{\varepsilon(\ox)}$. 

The latter is equivalent to the fact that the ramification of $V/U$ at $x$ is bounded by $R+$ by \ref{ram0}.

\begin{defi}\label{ram3}
Let $\cF$ be a locally constant constructible sheaf of $\Lambda$-modules on $U$, 
$R$ an effective rational divisor on $X$ with support in $D$, $x\in X$ and
$\ox$ a geometric point of $X$ above $x$. 
We say that the {\em ramification of $\cF$ at $\ox$ is bounded by $R+$} 
if $\cF$ satisfies the equivalent conditions of \eqref{ram2}, 
and that the {\em ramification of $\cF$ along $D$ is bounded by $R+$} if
the ramification of $\cF$ at $\ox$ is bounded by $R+$ for every geometric
point $\ox$ of $X$.
\end{defi}

\begin{prop}\label{ram4}
Let $\cF$ be a locally constant constructible sheaf of free $\Lambda$-modules of rank one on $U$
and let $R$ be an effective rational divisor on $X$ with support in $D$.  
Then the ramification of $\cF$ along $D$ is bounded by $R+$ if and only if  
$j^{(R)}_*(\cH(\cF))$ is locally constant constructible 
over an open neighborhood of $\delta^{(R)}(X)$ in $(X\Asterisk_kX)^{(R)}$.
\end{prop}

We may assume that $X$ and hence $U$ are connected. 
Let $V$ be a Galois torsor over $U$ with abelian group $G$ and let $\chi\colon G\rightarrow \Lambda^\times$ 
be an injective homomorphism such that we have an isomorphism (of sheaves with Galois descent data)
$\cF|V\simeq \Lambda(\chi)_V$, where $\Lambda(\chi)_V$ is the constant sheaf on $V$ of value $\Lambda$
with the descent datum defined by $\chi$.  
We keep the same notation as in the proof of \ref{ram2}. 
Then $W$ is an abelian Galois torsor over $U\times_kU$ trivializing $\cH(\cF)$.  
Assume first that the ramification of $\cF$ along $D$ is bounded by $R+$;  
so $\pi\colon Z\rightarrow (X\Asterisk_kX)^{(R)}$ is étale over an open neighborhood of $\varepsilon(X)$
in $Z$ by \ref{ram2}. It follows that $\pi$ is étale over an open neighborhood $P$ of of $\delta^{(R)}(X)$ 
in $(X\Asterisk_k X)^{(R)}$. 
Therefore, $\cH(\cF)$ extends to a locally constant constructible sheaf on $P$,
and hence $j^{(R)}_*(\cH(\cF))$ is locally constant constructible over $P$ (\cite{sga4} IX 2.14.1). 
Conversely, if $j^{(R)}_*(\cH(\cF))$ is locally constant constructible 
over an open neighborhood of $\delta^{(R)}(X)$ in $(X\Asterisk_kX)^{(R)}$, then the base change 
morphism 
\begin{equation}\label{ram4a}
\alpha\colon \delta^{(R)*}j^{(R)}_*(\cH(\cF))\rightarrow j_*(\cEnd(\cF))=j_*(\Lambda)=\Lambda
\end{equation}
is an isomorphism since its restriction to $U$ is an isomorphism (\cite{sga4} IX 2.14.1).

\begin{lem}\label{ram71}
Let $\cF$ be a locally constant constructible sheaf of $\Lambda$-modules on $U$, 
$\ox$ a geometric point of $X$ and $R$ and $R'$ be rational divisors on $X$ with support in $D$ such that $R'\geq R$. 
If the ramification of $\cF$ at $\ox$ (resp. along $D$) is bounded by $R+$, 
then it is also bounded by $R'+$. 
\end{lem}
It follows from  \ref{ram-cov38} and \ref{ram2}.

\begin{prop}\label{ram7}
Let $\cF$ be a locally constant constructible sheaf of $\Lambda$-modules on $U$, 
$R$ an effective rational divisor on $X$ with support in $D$, 
$f\colon (X',D')\rightarrow (X,D)$ a morphism of snc-pairs over $k$,
$U'=f^{-1}(U)$,  $\cF'=\cF|U'$, $R'=f^*(R)$, $x'\in X'$, $\ox'$ a geometric point of $X'$ above $x'$ and $\ox=f(\ox')$. 

{\rm (i)}\ If the ramification of $\cF$ at $\ox$ (resp. along $D$) is bounded by $R+$, 
then the ramification of $\cF'$ at $\ox'$ (resp. along $D'$) is bounded by $R'+$.

{\rm (ii)}\ Assume that the morphism $X'\rightarrow X$ is étale at $x'$
and that the divisors $D'$ and $f^*(D)$ are equal on an open neighborhood of $x'$ in $X'$. 
Then the ramification of $\cF$ at $\ox$ is bounded by $R+$
if and only if the ramification of $\cF'$ at $\ox'$ is bounded by $R'+$.

{\rm (iii)}\ Assume that the morphism $X'\rightarrow X$ is étale 
and surjective and that $D'=f^*(D)$. Then the ramification of $\cF$ along $D$ is bounded by $R+$
if and only if the ramification of $\cF'$ along $D'$ is bounded by $R'+$.
\end{prop}

(i) It follows from  \ref{ram-cov32}(i) and \ref{ram2}. 

(ii) By (i), it is enough to prove that if the ramification of $\cF'$ at $\ox'$ is bounded by $R'+$, 
then the ramification of $\cF$ at $\ox$ is bounded by $R+$. 
Also by (i), we may assume that $X'\rightarrow X$ is étale and that $D'=f^*(D)$. 
The morphism $f$ induces morphisms \eqref{log71a}
\begin{equation}\label{ram7a}
\xymatrix{
{X'\Asterisk_kX'}\ar[r]^{f_1}&{X\Asterisk_kX'}\ar[r]^{f_2}&{X\Asterisk_kX}},
\end{equation}
from which we obtain the morphisms \eqref{log71b}
\begin{equation}\label{ram7b}
\xymatrix{
{(X'\Asterisk_kX')^{(R')}}\ar[r]^{f_1^{(R')}}&{(X\Asterisk_kX')^{(R')}}\ar[r]^{f_2^{(R)}}&{(X\Asterisk_kX)^{(R)}}}.
\end{equation}
We put $\tf=f_2\circ f_1$ and $\tf^{(R)}=f_2^{(R)}\circ f_1^{(R')}$. 
Since $X\Asterisk_kX'\simeq (X\Asterisk_kX)\times_XX'$,  $f_2^{(R)}$ is smooth by \ref{log52}. 
On the other hand, $f_1^{(R')}$ is étale by \ref{log73}(i). 
Therefore, $\tf^{(R)}$ is smooth. Consider the commutative diagram 
\begin{equation}\label{ram7c}
\xymatrix{
X'\ar[ddd]_f\ar[rd]_{\delta^{(R')}}\ar@{}[rddd]|*+[o][F-]{1}&&&U'\ar[lll]_{j'}\ar[ld]^{\delta_{U'}}\ar[ddd]^{f_U}\ar@{}[llld]|*+[o][F-]{4}\\
&{(X'\Asterisk_kX')^{(R')}}\ar[d]_{\tf^{(R)}}\ar@{}[rd]|*+[o][F-]{2}&{U'\times_kU'}\ar[l]_-(0.4){j'^{(R')}}\ar[d]&\\
&{(X\Asterisk_kX)^{(R)}}&{U\times_kU}\ar[l]_-(0.5){j^{(R)}}&\\
X\ar[ru]^{\delta^{(R)}}\ar@{}[urrr]|*+[o][F-]{5}&&&U\ar[lu]_{\delta_U}\ar[lll]_{j}\ar@{}[luuu]|*+[o][F-]{3}}
\end{equation}
where the square $\xymatrix{\ar@{}|*+[o][F-]{4}}$ is the analogue of the square $\xymatrix{\ar@{}|*+[o][F-]{5}}$ for $(X',D', R')$. 
We denote by $\xymatrix{\ar@{}|*+[o][F-]{6}}$ the face given by the exterior square. Then all squares are cartesian except $\xymatrix{\ar@{}|*+[o][F-]{1}}$ and $\xymatrix{\ar@{}|*+[o][F-]{3}}$. 
Observe that $\cH(\cF)|(U'\times_kU')\simeq \cH(\cF')$ and $\cEnd(\cF)|U'\simeq \cEnd(\cF')$. 
It follows from (\cite{egr1} 1.2.4(i)) that we have a commutative diagram
\begin{equation}\label{ram7d}
\xymatrix{
{f^*\delta^{(R)*}j^{(R)}_*(\cH(\cF))}\ar[r]^{f^*(\alpha)}\ar@{=}[d]&{f^*j_*\delta_U^*(\cH(\cF))}\ar@{=}[r]\ar[d]^{\gamma'}&
{f^*j_*(\cEnd(\cF))}\ar[dd]^{\gamma}\\
{\delta^{(R')*}\tf^{(R)*}j^{(R)}_*(\cH(\cF))}\ar[d]_{\beta}&{j'_*f_U^*\delta_U^*(\cH(\cF))}\ar@{=}[d]&\\
{\delta^{(R')*}j'^{(R')}_*(\cH(\cF'))}\ar[r]^{\alpha'}&{j'_*\delta_{U'}^*(\cH(\cF'))}\ar@{=}[r]&{j'_*(\cEnd(\cF'))}}
\end{equation}
where $\alpha\colon \delta^{(R)*}j^{(R)}_*(\cH(\cF))\rightarrow j_*\delta_U^*(\cH(\cF))$,
$\alpha',\beta,\gamma$ and $\gamma'$ are the base change morphisms. Since $f$ and $\tf^{(R)}$ are smooth,
$\beta$ and $\gamma$ are isomorphisms. 
Hence, we can identify the stalks $\alpha_{\ox}$ and $\alpha'_{\ox'}$ of $\alpha$ at $\ox$
and $\alpha'$ at $\ox'$ respectively, which implies the required assertion.

(iii) It follows from (ii).

\begin{prop}\label{ram-loc3}
Let $\cF$ be a locally constant constructible sheaf of $\Lambda$-modules on $U$, 
$R$ an effective rational divisor on $X$ with support in $D$ and
$\ox$ a geometric point of $X$. Then the following conditions are equivalent~:

{\rm (i)}\ The ramification of $\cF$ at $\ox$ is bounded by $R+$. 

{\rm (ii)}\ There exists an étale neighborhood $f\colon X'\rightarrow X$ of $\ox$ such that 
if we put $U'=f^{-1}(U)$, $D'=f^*(D)$ and $R'=f^*(R)$, then the ramification of $\cF|U'$ 
along $D'$ is bounded by $R'+$.

{\rm (iii)}\ Condition {\rm (ii)} holds for $f$ an open immersion.   
\end{prop}
 Indeed,  (i)$\Rightarrow$(iii) follows from \ref{ram-cov75} and \ref{ram2}, 
 (iii)$\Rightarrow$(ii) is obvious and (ii)$\Rightarrow$(i)  is a consequence of  \ref{ram7}.

\begin{prop}\label{ram-loc2}
Let $\cF$ be a locally constant constructible sheaf of $\Lambda$-modules on $U$, 
$R$ an effective rational divisor on $X$ with support in $D$, $\xi$  a generic point of $D$,
$\oxi$ a geometric point of $X$ above $\xi$, 
$X_{(\oxi)}$ the corresponding strictly local scheme, $\eta$ its generic point and
$r$ the multiplicity of $R$ at $\xi$. Then the following conditions are equivalent~:

{\rm (i)}\ The ramification of $\cF$ at $\oxi$ is bounded by $R+$.

{\rm (ii)}\ The ramification of $\cF|\eta$ is bounded by $r+$ in the sense of {\rm (\cite{saito1} 1.28)}.

{\rm (iii)}\ The sheaf $\cF|\eta$ is trivialized by a finite étale connected covering $\eta'$ of $\eta$ such that  
the logarithmic ramification of $\eta'/\eta$ is bounded by $r+$ \eqref{rrt15}. 
\end{prop}

Indeed, (ii)$\Leftrightarrow$(iii) is the definition and (i)$\Rightarrow$(iii) follows from 
\ref{ram-cov8} and \ref{ram2}. We prove (iii)$\Rightarrow$(i). We may assume that $\eta'$
is Galois over $\eta$. By \ref{ram-loc4}, there exists an étale morphism
$f\colon X'\rightarrow X$, a geometric point $\oxi'$ of $X'$ above $\oxi$
and a Galois torsor $V'$ over $U'=f^{-1}(U)$ trivializing $\cF|U'$ such that
if we identify the strictly local schemes $X'_{(\oxi')}$ and $X_{(\oxi)}$ by $f$, 
there exists  an $\eta$-isomorphism $V'\times_{U'}\eta\simeq \eta'$. 
It follows from \ref{ram-cov8} and \ref{ram2} that the ramification of $\cF|U'$ at $\oxi'$ is bounded by $R'+$.
Then the ramification of $\cF$ at $\oxi$ is bounded by $R+$ by \ref{ram7}(ii). 

\begin{prop}\label{ram5}
Let $\cF$ be a locally constant constructible sheaf of $\Lambda$-modules on $U$
and $R$ an effective rational divisor on $X$ with support in $D$. 
Assume that the following conditions are satisfied~:

{\rm (i)} For every geometric point $\ox$ of $X$, if we denote by $X_{(\ox)}$ the corresponding 
strictly local scheme and put $U_1=X_{(\ox)}\times_XU$, the sheaf $\cF|U_1$ is trivialized by 
a Galois torsor over $U_1$ whose group has a normal $p$-Sylow subgroup. 

{\rm (ii)}\  For every geometric point $\oxi$ of $X$ above a generic point of $D$, 
the ramification of $\cF$ at $\oxi$ is bounded by $R+$.

Then the ramification of $\cF$ along $D$ is bounded by $R+$. 
\end{prop} 

Let $\ox$ be a geometric point of $X$ and let $V_1$ be 
a Galois torsor over $U_1=X_{(\ox)}\times_XU$ trivializing $\cF|U_1$
whose group $G$ has a normal $p$-Sylow subgroup. 
It is enough to prove that the ramification of $\cF$ at $\ox$ is bounded by $R+$. 
By \ref{ram-loc4} and \ref{ram7}(ii), we may assume that there exists a Galois torsor over $U$
of group $G$ trivializing $\cF$. Hence, there is a minimal Galois torsor $V$
over $U$ trivializing $\cF$ whose group has a normal $p$-Sylow subgroup. 
Then it follows from \ref{ram2} and \ref{ram-cov6} 
that the ramification of $\cF$ along $D$ is bounded by $R+$.

\begin{defi}\label{ram8}
Let $\cF$ be a locally constant constructible sheaf of $\Lambda$-modules on $U$.

{\rm (i)}\ Let $\xi$ be a generic point of $D$, $\oxi$ a geometric point of $X$ above $\xi$, 
$X_{(\oxi)}$ the corresponding strictly local scheme and $\eta$ the generic point of $X_{(\oxi)}$. 
We define the {\em conductor of $\cF$ at $\xi$} to be the minimum of the set of rational numbers $r\geq 0$
such that the ramification of $\cF|\eta$ is bounded by $r+$ in the sense of {\rm (\cite{saito1} 1.28)}, 
that is, $\cF|\eta$ is trivialized by a finite étale connected covering $\eta'$ of $\eta$ such that  
the logarithmic ramification of $\eta'/\eta$ is bounded by $r+$ \eqref{rrt15}.

{\rm (ii)}\ We define the {\em conductor of $\cF$ relatively to $X$} to be 
the effective rational divisor on $X$ with support in $D$ 
whose multiplicity at any generic point $\xi$ of $D$ is the conductor of $\cF$ at $\xi$. 
By \ref{ram-loc2}, it is also the minimum of the set of effective rational divisors $R$ on $X$ with support in $D$ 
such that for every geometric point $\oxi$ of $X$ above a generic point of $D$,
the ramification of $\cF$ at $\oxi$ is bounded by $R+$. 
\end{defi}

The terminology in (ii) may be slightly misleading as the ramification of $\cF$ along $D$ may not be bounded by $R+$
in general. However, we have the following~: 

\begin{prop}\label{ram10}
Let $\cF$ be a locally constant constructible sheaf of $\Lambda$-modules on $U$.  
Assume that the strong form of resolution of singularities {\rm (RS)} in {\rm \ref{ram-cov18}} holds. 
Then, there exists an snc-pair $(X',D')$ over $k$ and a proper morphism of snc-pairs $f\colon (X',D')\rightarrow (X,D)$
inducing an isomorphism $X'-D'\stackrel{\sim}{\rightarrow} U$ such that  
if we denote by $R'$ the conductor of $\cF$ relatively to $X'$,
the ramification of $\cF$ along $D'$ is bounded by $R'+$.
\end{prop}

This follows from \ref{ram-cov18} and \ref{ram2}. 

\begin{lem}[\cite{saito1} 2.21]\label{ram9}
Let $\cF$ be a locally constant constructible sheaf of $\Lambda$-modules on $U$. 
Then the following conditions are equivalent~:

{\rm (i)}\ $\cF$ is tamely ramified along $D$.

{\rm (ii)}\ The conductor of $\cF$ vanishes. 

{\rm (iii)}\ The ramification of $\cF$ along $D$ is bounded by $0+$.
\end{lem}

\subsection{}\label{ram11}
Let $R$ be an effective divisor on $X$ with support in $D$. 
We know \eqref{tub26} that $(X\Asterisk_kX)^{(R)}$ is smooth over $X$ and that 
\begin{equation}\label{ram11a}
E^{(R)}=(X\Asterisk_k X)^{(R)}\times_XR
\end{equation}
is canonically isomorphic to the vector bundle $\bV(\Omega^1_{X/k}(\log D)\otimes_{\co_X}\co_{X}(R))\times_{X}R$
over $R$ \eqref{log7c}. We denote by $\chE^{(R)}$ the dual vector bundle.

Let $Y$ be an $X$-scheme separated of finite type over $X$. 
We put  $V=Y\times_XU$, 
$R_Y=R\times_XY$, $E^{(R)}_{Y}=E^{(R)}\times_XY$  and  $\chE^{(R)}_{Y}=\chE^{(R)}\times_XY$, and denote by
$j_Y\colon V\rightarrow Y$ and $j_Y^{(R)}\colon U\times_kV\rightarrow (X\Asterisk_kX)^{(R)}\times_XY$
the canonical injections. Consider the following commutative diagram with Cartesian squares.
\begin{equation}\label{ram11b}
\xymatrix{
{E^{(R)}_{Y}}\ar[r]\ar[d]&{(X\Asterisk_kX)^{(R)}\times_XY}\ar[d]^{\pr_2}&{U\times_kV}\ar[d]\ar[l]_-(0.4){j_Y^{(R)}}\\
{R_Y}\ar[r]&Y&V\ar[l]}
\end{equation}
Let $\cG$ be a sheaf of $\Lambda$-modules on $U\times_kV$. We 
call the sheaf over $E^{(R)}_{Y}$ defined by 
\begin{equation}\label{ram11c}
\nu_R(\cG, Y)=j^{(R)}_{Y*}(\cG)|E^{(R)}_{Y}
\end{equation}
the {\em $R$-specialization} of $\cG$ and denote it by $\nu_R(\cG, Y)$.

Let $f\colon Z\rightarrow Y$ be a separated morphism of finite type and $W=f^{-1}(V)$. We denote by 
\begin{eqnarray}
f^{(R)}\colon (X\Asterisk_kX)^{(R)}\times_XZ&\rightarrow& (X\Asterisk_kX)^{(R)}\times_XY\label{ram11d}\\
f^{(R)}_E\colon E^{(R)}_{Z}&\rightarrow& E^{(R)}_{Y}\label{ram11e}
\end{eqnarray} 
the morphisms induced by $f$. Then we have a base change morphism
\begin{equation}\label{ram11f}
f^{(R)*} j^{(R)}_{Y*}(\cG)\rightarrow j^{(R)}_{Z*}(\cG|(U\times_kW)), 
\end{equation}
from which we get the morphism 
\begin{equation}\label{ram11g}
f^{(R)*}_{E}(\nu_R(\cG, Y))\rightarrow \nu_{R}(\cG|(U\times_kW), Z).
\end{equation}

\begin{prop}\label{ram13}
Let $\cF$ be a locally constant constructible sheaf of $\Lambda$-modules on $U$,  
$R$ an effective divisor on $X$ with support in $D$, 
$f\colon Y\rightarrow X$ a separated morphism of finite type and $V=f^{-1}(U)$.
Assume that the ramification of $\cF$ along $D$ is bounded by $R+$.  
Then with the notation of \eqref{ram11}~:

{\rm (i)}\ The sheaves $j^{(R)}_{Y*}(\Lambda_U\boxtimes (\cF|V))$ 
and $\pr_2^*(j_{Y*}(\cF|V))$ are isomorphic over $(X\Asterisk_k X)^{(R)}\times_XY$.

{\rm (ii)}\ There exists an étale morphism $Z_0\rightarrow (X\Asterisk_kX)^{(R)}$ whose image contains $E^{(R)}$,
such that the pull-backs over $Z_0\times_XY$ of the sheaves
$j^{(R)}_{Y*}(\cF\boxtimes \Lambda_V)$ and $j^{(R)}_{Y*}(\Lambda_U\boxtimes (\cF|V))$ are isomorphic.

{\rm (iii)}\ There exists a canonical surjective morphism
\begin{equation}\label{ram13a}
f^{(R)*}_{E}(\nu_R(\cH(\cF),X))\otimes  \nu_R(\Lambda_U\boxtimes (\cF|V),Y)\rightarrow \nu_R(\cF\boxtimes \Lambda_V,Y).
\end{equation}

{\rm (iv)}\ Assume moreover that  
$Y$ is normal, that $V$ is dense in $Y$ and that $\cF|V$ is constant. Then $\nu_R(\cF\boxtimes \Lambda_V, Y)$ is 
locally constant and additive \eqref{add1}; in particular, its Fourier dual support is 
the underlying space of a closed subscheme of $\chE^{(R)}_{Y}$
which is finite over $R_Y$ \eqref{add5}.
\end{prop}

Note that statement (iv) will be extended and reinforced in \ref{ram155}.

(i) It follows from the smooth base change theorem as $\pr_2\colon (X\Asterisk_kX)^{(R)}\rightarrow X$
is smooth. 

(ii) Let $\tV$ be a Galois torsor over $U$ of group $G$ trivializing $\cF$
such that the ramification of $\tV/U$ is bounded by $R+$ \eqref{ram2}.
We denote by $\Delta\colon G\rightarrow G\times G$ the diagonal homomorphism,
by $W$ the quotient of $\tV\times_k\tV$ by $\Delta(G)$, by $Z$ the integral closure of $(X\Asterisk_kX)^{(R)}$
in $W$ and by $\pi\colon Z\rightarrow (X\Asterisk_kX)^{(R)}$ the canonical morphism.  
Let $Z_0$ be the maximal open subscheme of $Z$ which is étale over $(X\Asterisk_kX)^{(R)}$. 
We know by \ref{ram-cov41}(i) that $Z_0\times_XR$ is a group scheme and that
$\pi\times_XR\colon Z_0\times_XR\rightarrow E^{(R)}$
is a surjective étale morphism of group schemes over $R$.
Let $j'\colon W\rightarrow Z_0$ 
be the canonical injection~; so we have the following Cartesian diagram.
\begin{equation}\label{ram13b}
\xymatrix{
{W\times_UV}\ar[r]^{j'_Y}\ar[d]&{Z_0\times_XY}\ar[d]\\
{U\times_kV}\ar[r]^-(0.5){j_Y^{(R)}}&{(X\Asterisk_kX)^{(R)}\times_XY}}
\end{equation}
By the smooth base change theorem, the pull-back over $Z_0\times_XY$ of the sheaves
$j^{(R)}_{Y*}(\cF\boxtimes \Lambda_V)$ and $j^{(R)}_{Y*}(\Lambda_U\boxtimes (\cF|V))$ are isomorphic 
to $j'_{Y*}((\cF\boxtimes \Lambda_U)|(W\times_UV))$ and $j'_{Y*}((\Lambda_U\boxtimes \cF)|(W\times_UV))$ respectively.
On the other hand, we have an isomorphism $(\Lambda_{\tV}\boxtimes (\cF|\tV))\stackrel{\sim}{\rightarrow} 
((\cF|\tV)\boxtimes \Lambda_{\tV})$.  We deduce by Galois descent an isomorphism 
\begin{equation}\label{ram13c}
s\colon (\Lambda_U\boxtimes \cF)|W\stackrel{\sim}{\rightarrow} 
(\cF\boxtimes \Lambda_U)|W,
\end{equation}
and hence an isomorphism
\begin{equation}\label{ram13k}
j'_{Y*}((\Lambda_U\boxtimes \cF)|(W\times_UV))\stackrel{\sim}{\rightarrow} j'_{Y*}((\cF\boxtimes \Lambda_U)|(W\times_UV)).
\end{equation}
Note that the isomorphism $s$ can also be obtained by Galois descent from the universal isomorphism 
of $G$-torsors over $W$ \eqref{ram-cov2}
\begin{equation}
(U\times_kV)\times_{U\times_kU}W\stackrel{\sim}{\rightarrow}(V\times_kU)\times_{U\times_kU}W.
\end{equation}

(iii)  We have a morphism on $(X\Asterisk_kX)^{(R)}\times_XY$
\begin{equation}\label{ram13d}
j^{(R)}_{Y*}(\cH(\cF)|(U\times_kV))
\otimes j^{(R)}_{Y*}(\Lambda_U\boxtimes (\cF|V))\rightarrow j^{(R)}_{Y*}(\cF\boxtimes \Lambda_V)
\end{equation}
deduced by adjunction from the natural morphism on $U\times_kV$
\begin{equation}\label{ram13e}
(\cH(\cF)|(U\times_kV))\otimes(\Lambda_U\boxtimes (\cF|V))\rightarrow 
\cF\boxtimes \Lambda_V.
\end{equation}
We get from \eqref{ram13d} by pull-back to $E_{Y}$ a morphism 
\begin{equation}\label{ram13f}
\nu_R(\cH(\cF)|(U\times_kV),Y)\otimes \nu_R(\Lambda_U\boxtimes (\cF|V),Y)\rightarrow \nu_R(\cF\boxtimes \Lambda_V,Y).
\end{equation}
On the other hand, we have a canonical morphism \eqref{ram11g}
\begin{equation}\label{ram13g}
f^{(R)*}_E(\nu_R(\cH(\cF),X))\rightarrow \nu_R(\cH(\cF)|(U\times_kV),Y).
\end{equation}
We take for \eqref{ram13a} the morphism induced by \eqref{ram13f} and \eqref{ram13g}. 
We will prove that it is surjective. By the smooth base change theorem, 
the pull-back of the morphism \eqref{ram13d} over $Z_0\times_XY$ is the morphism 
\begin{equation}\label{ram13h}
j'_{Y*}(\cH(\cF)|(W\times_UV))\otimes j'_{Y*}((\Lambda_U\boxtimes \cF)|(W\times_UV))
\rightarrow j'_{Y*}((\cF\boxtimes \Lambda_U)|(W\times_UV))
\end{equation}
obtained by adjunction from the pull-back of the morphism \eqref{ram13e} over $W\times_UV$. 
We also have a canonical isomorphism
\begin{equation}\label{clean2h}
\Gamma(Z_0,j^{(R)}_*(\cH(\cF)))\simeq \Gamma(W,\cH|W).
\end{equation}
So we may consider the isomorphism $s$ \eqref{ram13c} as a section in $\Gamma(Z_0,j^{(R)}_*(\cH(\cF)))$.
It is clear that the pairing \eqref{ram13h} evaluated at the image of $f^{(R)*}(s)$
in $\Gamma(Z_0\times_XY, j'_{Y*}(\cH(\cF)|(W\times_UV)))$ induces the isomorphism \eqref{ram13k}.
Since $E^{(R)}\subset \pi(Z_0)$, we conclude that \eqref{ram13a} is surjective. 

(iv) Since $j_{Y*}(\cF|V)$ is constant (\cite{sga4} IX 2.14.1), the pull-back of 
$\nu_R(\cF\boxtimes\Lambda_V,Y)$ over $Z_0\times_XR_Y$ is constant by (i) and (ii).  
Hence $\nu_R(\cF\boxtimes\Lambda_V,Y)$ is locally constant.
For every geometric point $\ox$ of $R$, if we denote by $G_\ox$ the neutral connected component
of $Z_0\times_X\ox$, then the morphism $G_\ox\rightarrow E^{(R)}_\ox$ induced by $\pi$, is a finite étale surjective
morphism of group schemes over $\ox$ by \ref{ram-cov41}(ii).  
Therefore, $\nu_R(\cF\boxtimes\Lambda_V,Y)$ is additive by \ref{add16}. 
The last assertion follows from \ref{add11}. 

\begin{cor}[\cite{saito1}  2.25]\label{ram12}
Let $\cF$ be a locally constant constructible sheaf of $\Lambda$-modules on $U$ and let
$R$ be an effective divisor on $X$ with support in $D$ such that the ramification of $\cF$ 
along $D$ is bounded by $R+$. Then $\nu_R(\cH(\cF), X)$ is additive. 
\end{cor}

Let $V$ be a finite Galois torsor over $U$ trivializing $\cF$ such that the ramification of $V/U$ along $D$
is bounded by $R+$ \eqref{ram2}. We denote by $Y$ the integral closure of $X$ in $V$ and by $f\colon Y\rightarrow X$
the canonical morphism. In the following we take again the notation of \ref{ram13} and its proof (with $\tV=V$). 
Consider the following commutative diagram with Cartesian squares.
\begin{equation}\label{ram12a}
\xymatrix{
W\ar[r]^{j'}\ar[d]&{Z_0}\ar[d]^{\pi}\\
{U\times_kU}\ar[r]^-(0.5){j^{(R)}}\ar[d]&{(X\Asterisk_kX)^{(R)}}\ar[d]^{\pr_2}\\
U\ar[r]^j&X}
\end{equation}
By \ref{ram1-2}, the isomorphism $s$ \eqref{ram13c} induces isomorphisms 
\begin{equation}
\cH(\cF)|W\stackrel{\sim}{\rightarrow} \cEnd((\Lambda_U\boxtimes \cF)|W)
\stackrel{\sim}{\rightarrow}(\cEnd(\cF))|W.
\end{equation}
Since $\pi$ and $\pr_2$ are smooth, by the smooth base change theorem relatively to \eqref{ram12a}, we get an isomorphism
\begin{equation}
\pi^*(j^{(R)}_*(\cH(\cF)))\stackrel{\sim}{\rightarrow}\pi^*(\pr_2^*(j_*(\cEnd(\cF)))).
\end{equation} 
It follows that $\nu_X(\cH(\cF),X)$ is locally constant on the geometric fibers of $E^{(R)}\rightarrow R$. 
Since $\pi(Z_0)$ contains $E^{(R)}$, we conclude that $\nu_R(\cH(\cF),X)$ is locally constant on all geometric 
fibers of the vector bundle $E^{(R)}\rightarrow R$. 

On the other hand, the open immersion $j_Y^{(R)}$ is schematically dense. 
Hence the morphism \eqref{ram13g} is injective by \ref{ram00}. 
We fix a surjective morphism $\Lambda^n\rightarrow \cF|V$, from which we deduce 
an injective morphism $\cH(\cF)|(U\times_kV) \rightarrow \cF\boxtimes \Lambda^n_V$ \eqref{ram1-2}.
The latter induces an injective morphism
\begin{equation}
\nu_{R}(\cH(\cF)|(U\times_kV),Y)\rightarrow \nu_{R}(\cF\boxtimes \Lambda^n_V,Y).
\end{equation}
Composing with \eqref{ram13g}, we obtain an injective morphism 
\begin{equation}\label{ram12b}
f_E^{(R)*}(\nu_{R}(\cH(\cF),X))\rightarrow \nu_{R}(\cF\boxtimes \Lambda^n_V,Y).
\end{equation}
Then $f_E^{(R)*}(\nu_{R}(\cH(\cF),X))$ is additive by \ref{add8}.
Since $f$ is surjective, $\nu_R(\cH(\cF),X)$ is additive.

\begin{cor}\label{ram133}
Let $\cF$ be a locally constant constructible sheaf of $\Lambda$-modules on $U$,  
$R$ an effective divisor on $X$ with support in $D$, 
$f\colon Y\rightarrow X$ a separated morphism of finite type and $V=f^{-1}(U)$.
Assume that the ramification of $\cF$ along $D$ is bounded by $R+$.  
Then $\nu_R(\cF\boxtimes \Lambda_V, Y)$ is additive, and its Fourier dual support
is contained in the inverse image by the canonical projection $\chE_{Y}^{(R)}\rightarrow \chE^{(R)}$
of the Fourier dual support of $\nu_R(\cH(\cF),X)$.
\end{cor}

It follows from \ref{ram13}(i) that $\nu_R(\Lambda_U\boxtimes(\cF|V),Y)$ is constant on the fibers 
of the vector bundle $E^{(R)}_{Y}\rightarrow R_Y$. Hence by \ref{ram13}(ii), $\nu_R(\cF\boxtimes \Lambda_V,Y)$ is 
locally constant on the geometric fibers of $E^{(R)}_{Y}\rightarrow R_Y$. 
Then the proposition follows from \ref{add8}, \ref{ram13}(iii) and \ref{ram12}.

\begin{cor}\label{ram14}
We keep the assumptions of \eqref{ram133} and assume moreover that  
$Y$ is normal, that $V$ is dense in $Y$ and that $\cF|V$ is constant. Then
the Fourier dual support of $ \nu_R(\cF\boxtimes \Lambda_V,Y)$ 
is the inverse image by the canonical projection $\chE_{Y}^{(R)}\rightarrow \chE^{(R)}$
of the Fourier dual support of $\nu_R(\cH(\cF),X)$. 
\end{cor}

By \ref{ram133}, it is enough to prove that the inverse image of the Fourier  dual support of 
$\nu_R(\cH(\cF),X)$ by the canonical projection $\chE^{(R)}_{Y}\rightarrow \chE^{(R)}$
is contained in the Fourier dual support of $\nu_{R}(\cF\boxtimes \Lambda_V,Y)$.
This follows from the second part of the proof of \ref{ram12}. Indeed, 
since the sheaves $\nu_{R}(\cH(\cF),X)$ and $\nu_{R}(\cF\boxtimes \Lambda^n_V,Y)$ 
are additive by \ref{ram12} and \ref{ram13}(iv), the required assertion follows from the injective 
morphism \eqref{ram12b} by \ref{add8}.

\begin{cor}\label{ram15}
Let $\cF$ be a locally constant constructible sheaf of $\Lambda$-modules on $U$ and 
let $R$ be an effective divisor on $X$ with support in $D$ such that the ramification of $\cF$ 
along $D$ is bounded by $R+$.
Then the Fourier dual support of $\nu_R(\cH(\cF),X)$ is 
the underlying space of a closed subscheme of $\chE^{(R)}$
which is finite over $R$.
\end{cor} 

Let $V$ be a Galois torsor over $U$ trivializing $\cF$ and $Y$ the integral closure of $X$ in $V$.
It follows from \ref{ram14} that the Fourier dual support of $\nu_R(\cH(\cF),X)$
is the image by the canonical projection $\chE_{Y}^{(R)}\rightarrow \chE^{(R)}$
of the Fourier dual support of  $ \nu_R(\cF\boxtimes \Lambda_V,Y)$. Hence, 
the assertion follows from \ref{ram13}(iv). 

\begin{prop}\label{ram155}
Let $\cF$ be a locally constant constructible sheaf of $\Lambda$-modules on $U$,  
$R$ an effective divisor on $X$ with support in $D$, 
$f\colon Y\rightarrow X$ a separated morphism of finite type and $V=f^{-1}(U)$.
Assume that the ramification of $\cF$ along $D$ is bounded by $R+$.  
Then~:

{\rm (i)}\ $\nu_R(\cH(\cF), X)$ is additive and its Fourier dual support is 
the underlying space of a closed subscheme of $\chE^{(R)}$ which is finite over $R$.

{\rm (ii)}\  $\nu_R(\cF\boxtimes \Lambda_V, Y)$ is additive, and its Fourier dual support
is contained in the inverse image by the canonical projection $\chE_{Y}^{(R)}\rightarrow \chE^{(R)}$
of the Fourier dual support of $\nu_R(\cH(\cF),X)$.

{\rm (iii)}\ Assume moreover that $Y$ is normal, that $V$ is dense in $Y$ and that $\cF|V$ is constant. 
Then $\nu_R(\cF\boxtimes \Lambda_V, Y)$ is locally constant and additive, 
and its Fourier dual support is the inverse image by the canonical projection $\chE_{Y}^{(R)}\rightarrow \chE^{(R)}$
of the Fourier dual support of $\nu_R(\cH(\cF),X)$.
\end{prop}

This is a summary of results proved in \ref{ram13}, \ref{ram12}, \ref{ram133}, \ref{ram14} and \ref{ram15}. 

\subsection{}\label{cds1}
Let $f\colon (X',D')\rightarrow (X,D)$ be a {\em log-smooth} morphism of snc-pairs over $k$,  
$R$ an effective divisor on $X$ with support in $D$, $U_1=X'-D'$, $U'=f^{-1}(U)$ and $R'=f^*(R)$. 
The morphism $f$ induces morphisms \eqref{log71a}
\begin{equation}\label{cds1a}
\xymatrix{
{X'\Asterisk_kX'}\ar[r]^{f_1}&{X\Asterisk_kX'}\ar[r]^{f_2}&{X\Asterisk_kX}},
\end{equation}
from which we obtain the morphisms \eqref{log71b}
\begin{equation}\label{cds1b}
\xymatrix{
{(X'\Asterisk_kX')^{(R')}}\ar[r]^{f_1^{(R')}}&{(X\Asterisk_kX')^{(R')}}\ar[r]^{f_2^{(R)}}&{(X\Asterisk_kX)^{(R)}}}.
\end{equation}
We put 
\begin{eqnarray}
E^{(R)}&=&(X\Asterisk_kX)^{(R)}\times_XR,\label{cds1c}\\
E'^{(R')}&=&(X'\Asterisk_kX')^{(R')}\times_{X'}R'.\label{cds1d}
\end{eqnarray} 
We denote by $\chE^{(R)}$ and $\chE'^{(R')}$ the dual vector bundles.
For a sheaf of $\Lambda$-modules $\cG$ (resp. $\cG'$) on $U\times_kU$ (resp. $U'\times_kU'$), 
we denote by $\nu_{R}(\cG,X)$ (resp. $\nu'_{R'}(\cG',X')$) its $R$-specialization (resp. $R'$-specialization) 
in the sense of \eqref{ram11c} relatively to the snc-pair $(X,D)$ (resp. $(X',D')$). 

The morphism $f$ induces an exact sequence 
\begin{equation}\label{cds1e}
0\rightarrow f^*(\Omega^1_{X/k}(\log D))\rightarrow \Omega^1_{X'/k}(\log D')\rightarrow \Omega^1_{(X',D')/(X,D)}
\rightarrow 0,
\end{equation}
which is locally split. Hence, we have a surjective morphism of vector bundles over $R'$ 
\begin{equation}\label{cds1f}
\phi\colon E'^{(R')}\rightarrow E^{(R)}_{X'}.
\end{equation}
We denote by $\chphi\colon \chE^{(R)}_{X'}\rightarrow \chE'^{(R')}$ the dual morphism of $\phi$, 
which is a closed immersion.

\begin{prop}\label{cds2}
We keep the notation of \eqref{cds1} and let $\cF$ be a locally constant constructible sheaf 
of $\Lambda$-modules on $U$, $\cF_1=\cF|U_1$ and $\cF'=\cF|U'$. 
Assume that the ramification of $\cF$ along $D$ is bounded by $R+$ and that $f$ is log-smooth.
Then the Fourier dual support of $\nu'_{R'}(\cH(\cF_1),X')$ is the image by $\chphi$ of the inverse 
image of the Fourier dual support of  $\nu_R(\cH(\cF),X)$ by the canonical projection 
$\chE^{(R)}_{X'}\rightarrow \chE^{(R)}$.
\end{prop}

Let $V'$ be a Galois torsor over $U'$ trivializing $\cF'$, $V_1=V'\times_{U'}U_1$ and
$Y'$ the integral closure of $X'$ in $V'$.  
We denote by $\nu_R(\cF\boxtimes\Lambda_{V'},Y')$ the $R$-specialization of $\cF\boxtimes\Lambda_{V'}$
over $E^{(R)}_{Y'}$ in the sense of \eqref{ram11c}. It is an additive sheaf by \ref{ram13}(iv), and its 
Fourier dual support is the inverse image of the Fourier dual support of $\nu_R(\cH(\cF),X)$ 
by the canonical projection $\chE^{(R)}_{Y'}\rightarrow \chE^{(R)}$ \eqref{ram14}. 
On the other hand, the ramification of $\cF_1$ along $D'$ is bounded by $R'+$ by \ref{ram7}(i). 
We denote by $\nu'_{R'}(\cF_1\boxtimes\Lambda_{V_1},Y')$ the $R'$-specialization  
of $\cF_1\boxtimes\Lambda_{V_1}$ over $E'^{(R')}_{Y'}$ in the sense of \eqref{ram11c}  
relatively to the snc-pairs $(X',D')$. It is an additive sheaf and its Fourier dual support is the inverse image 
of the Fourier dual support of $\nu'_{R'}(\cH(\cF_1),X')$ by the canonical projection 
$\chE'^{(R')}_{Y'}\rightarrow \chE'^{(R')}$.
Since the canonical morphism $Y'\rightarrow X'$ is surjective, it is enough to prove that 
the Fourier dual support of $\nu'_{R'} (\cF_1\boxtimes\Lambda_{V_1},Y')$ is the image 
by $\chphi_{Y'}$ of  the Fourier dual support of $\nu_R(\cF\boxtimes\Lambda_{V'},Y')$.

On the one hand, $f_1^{(R')}$ is smooth by \ref{log73}(ii).
On the other hand, since the canonical morphism $X\Asterisk_kX'\rightarrow (X\Asterisk_kX)\times_{X}X'$ 
is an isomorphism \eqref{fram17b}, the morphism
\begin{equation}\label{cds2a}
(X\Asterisk_kX')^{(R')} \rightarrow (X\Asterisk_kX)^{(R)}\times_{X}X'
\end{equation}
induced by $f_2^{(R)}$ is an isomorphism by \ref{log52}. 
We denote by $U'\Asterisk_kU'$ the framed self-product of $(U',D'|U')$; 
so $U'\Asterisk_kU'=(X'\Asterisk_kX')\times_{(X'\times_kX')}(U'\times_kU')$. 
We have a commutative diagram 
\begin{equation}\label{cds3a}
\xymatrix{
{E'^{(R')}_{Y'}}\ar[r]\ar[d]_{\phi_{Y'}}&{(X'\Asterisk_kX')^{(R')}\times_{X'}Y'}
\ar[d]_{f_1^{(R')}\times_{X'}Y'}&{(U'\Asterisk_kU')\times_{U'}V'}\ar[l]\ar[d]&
{U_1\times_kV_1}\ar[l]_-(0.5)u\ar[ld]\ar@/_2pc/[ll]_-(0.5){j'^{(R')}_{Y'}}\\
{E_{Y'}^{(R)}}\ar[r]&{(X\Asterisk_kX')^{(R')}\times_{X'}Y'}&{U\times_kV'}\ar[l]_-(0.4){j_{Y'}^{(R)}}&}
\end{equation}
with Cartesian squares, where $u$, $j^{(R)}_{Y'}$ and $j'^{(R')}_{Y'}$ are the canonical injections. 

Since $U'\Asterisk_kU'$ is smooth over $U'$, $(U'\Asterisk_kU')\times_{U'}V'$ is normal and 
$u$ is dominant.  Therefore, the adjunction morphism 
\begin{equation}\label{cds3b}
(\cF\boxtimes \Lambda_{V'})|((U'\Asterisk_kU')\times_{U'}V')\rightarrow u_*(\cF_1\boxtimes \Lambda_{V_1})
\end{equation}
is an isomorphism by (\cite{sga4} IX 2.14.1). Then by the smooth base change theorem relatively 
to the Cartesian right square in \eqref{cds3a}, we have an isomorphism 
\begin{equation}\label{cds3c}
\phi_{Y'}^*(\nu_{R'}(\cF\boxtimes \Lambda_{V'},Y'))\stackrel{\sim}{\rightarrow} \nu'_{R'}(\cF_1\boxtimes \Lambda_{V_1},Y').
\end{equation}
Since $\chphi$ is a closed immersion, the required assertion follows from \eqref{cds3c} and \eqref{add3f}.

\begin{defi}\label{clean6}
Let $\cF$ be a locally constant constructible sheaf of $\Lambda$-modules on $U$. 

{\rm (i)}\ Let $\xi$ be a generic point of $D$, $X_{(\xi)}$ the henselization of $X$ at $\xi$, 
$\eta_\xi$ the generic point of $X_{(\xi)}$, 
$\oeta_\xi$ a geometric generic point of $X_{(\xi)}$ and $\cG_{\xi}$ the Galois group
of $\oeta_\xi$ over $\eta_\xi$. We say that $\cF$ is {\em isoclinic} at $\xi$ if the representation $\cF_{\oeta_\xi}$
of $\cG_{\xi}$ is isoclinic \eqref{rrt4}.

{\rm (ii)}\ We say that $\cF$ is {\em isoclinic along $D$} if it is isoclinic at all generic points of $D$.  
\end{defi}
 
\begin{defi}\label{clean7}
Let $\cF$ be a locally constant constructible sheaf of $\Lambda$-modules on $U$ which is isoclinic along $D$
and let $R$ be the conductor of $\cF$ relatively to $X$ \eqref{ram8}.  We say that $\cF$ is {\em clean} along $D$ 
if the following conditions are satisfied~:
\begin{itemize}
\item[{\rm (i)}] the ramification of $\cF$ along $D$ is bounded by $R+$; 

\item[{\rm (ii)}] there exists a log-smooth morphism of snc-pairs $f\colon (X',D')\rightarrow (X,D)$ over $k$
such that the morphism $X'\rightarrow X$ is faithfully flat, that $R'=f^*(R)$ has integral coefficients, 
and if we put $U'=X'-D'$ and $\cF'=\cF|U'$, 
that the $R'$-specialization $\nu'_{R'}(\cH(\cF'), X')$ of $\cH(\cF')$ in the sense of \eqref{ram11c}
relatively to $(X',D')$, is additive and non-degenerate \eqref{add5}. 
\end{itemize}
\end{defi}

Note that condition (i) implies that the ramification of $\cF'$ along $D'$ is bounded by $R'+$ by \ref{ram7}(i). 
Therefore, $\nu'_{R'}(\cH(\cF'), X')$ is additive by \ref{ram12}, and its Fourier dual support is the underlying space 
of a closed subscheme $S'$ of $\chE'^{(R')}=(X'\Asterisk_kX')^{(R')}\times_{X'}R'$, which is finite over $R'$ by \ref{ram15}.
Hence, the condition (ii) is equivalent to the fact that $S'$ does not meet the zero-section of $\chE'^{(R')}$.

\begin{prop}\label{clean71}
Let $\cF$ be a locally constant constructible sheaf of $\Lambda$-modules on $U$,
$R$ the conductor of $\cF$ relatively to $X$, $f\colon (X^\dagger,D^\dagger)\rightarrow (X,D)$ a log-smooth morphism 
of snc-pairs over $k$, $U^\dagger=X^\dagger-D^\dagger$, $\cF^\dagger=\cF|U^\dagger$ and $R^\dagger=f^*(R)$. 
Assume that $\cF$ is isoclinic and clean along $D$ and that $R^\dagger$ has integral coefficients. 
Then the $R^\dagger$-specialization $\nu^\dagger_{R^\dagger}(\cH(\cF^\dagger), X^\dagger)$ of $\cH(\cF^\dagger)$ 
in the sense of \eqref{ram11c} relatively to $(X^\dagger,D^\dagger)$, is additive and non-degenerate.
\end{prop}

By \ref{clean72}, there exists a commutative diagram of log-smooth morphisms of snc-pairs over $k$  
\begin{equation}
\xymatrix{
{(X^\ddagger,D^\ddagger)}\ar[r]^{f'}\ar[d]_{g^\dagger}&{(X',D')}\ar[d]^g\\
{(X^\dagger,D^\dagger)}\ar[r]^f&{(X,D)}}
\end{equation}
such that $X'\rightarrow X$ and $X^\ddagger\rightarrow X^\dagger$ are faithfully flat
and if we put $U'=X'-D'$, $\cF'=\cF|U'$ and $R'=g^*(R)$, that 
the $R'$-specialization $\nu'_{R'}(\cH(\cF'), X')$ of $\cH(\cF')$ 
in the sense of \eqref{ram11c} relatively to $(X',D')$, is additive and non-degenerate.
We put $U^\ddagger=X^\ddagger-D^\ddagger$, $\cF^\ddagger=\cF|U^\ddagger$ and $R^\ddagger=f'^*(R')$,
and denote by $\nu^\ddagger_{R^\ddagger}(\cH(\cF^\ddagger), X^\ddagger)$ 
the $R^\ddagger$-specialization of $\cH(\cF^\ddagger)$ 
in the sense of \eqref{ram11c} relatively to $(X^\ddagger,D^\ddagger)$.
We put
\begin{eqnarray}
E'^{(R')}&=&(X'\Asterisk_kX')^{(R')}\times_{X'}R',\\
E^{\dagger(R^\dagger)}&=&(X^\dagger\Asterisk_kX^\dagger)^{(R^\dagger)}\times_{X^\dagger}R^\dagger,\\
E^{\ddagger(R^\ddagger)}&=&(X^\ddagger\Asterisk_kX^\ddagger)^{(R^\ddagger)}\times_{X^\ddagger}R^\ddagger,
\end{eqnarray}
and denote by $\chE'^{(R')}, \chE^{\dagger(R^\dagger)}$ and $\chE^{\ddagger(R^\ddagger)}$ the dual vector bundles. 
Let $S'$ (resp. $S^\dagger$, resp. $S^\ddagger$) be the Fourier dual support of 
$\nu'_{R'}(\cH(\cF'), X')$ (resp. $\nu^\dagger_{R^\ddagger}(\cH(\cF^\dagger), X^\ddagger)$,
resp. $\nu^\ddagger_{R^\ddagger}(\cH(\cF^\ddagger), X^\ddagger)$) in 
$\chE'^{(R')}$ (resp. $\chE^{\dagger(R^\dagger)}$, resp. $\chE^{\ddagger(R^\ddagger)}$). 
The morphisms $f'$ and $g^\dagger$ induce, as in \eqref{cds1f}, surjective morphisms of vector bundles 
\begin{eqnarray}
\phi\colon E^{\ddagger(R^\ddagger)}&\rightarrow& E'^{(R')}_{X^\ddagger},\\
\psi\colon E^{\ddagger(R^\ddagger)}&\rightarrow& E^{\dagger(R')}_{X^\ddagger}.
\end{eqnarray}
Consider the diagram
\begin{equation}
\xymatrix{
{E'^{(R')}_{X^\ddagger}}\ar[r]^{\chphi}\ar[d]_{\pi'}&{E^{\ddagger(R^\ddagger)}}&
{E^{\dagger(R^\dagger)}_{X^\ddagger}}\ar[l]_{\chpsi}\ar[d]^{\pi^\dagger}\\
{E'^{(R')}}&&{E^{\dagger(R^\dagger)}}}
\end{equation}
where $\pi'$ and $\pi^\dagger$ are the canonical projections, 
and $\chphi$ and $\chpsi$ are the morphisms dual to $\phi$ and $\psi$ respectively. 
Note that $\chphi$ and $\chpsi$ are closed immersions. 
Then by \ref{cds2}, we have 
\begin{equation}
\chphi(\pi'^{-1}(S'))=S^\ddagger=\chpsi(\pi^{\dagger -1}(S^\dagger)).
\end{equation}
By assumption $S'$ does not meet the zero-section of $\chE'^{(R')}$. Hence, $\pi^{\dagger -1}(S^\dagger)$
does not meet the zero-section of $\chE^{\dagger(R^\dagger)}_{X^\ddagger}$. 
Since $X^\ddagger\rightarrow X^\dagger$ is surjective, we deduce that $S^\dagger$
does not meet the zero section of $E^{\dagger(R^\dagger)}$, which implies the required assertion. 

\begin{defi}\label{clean8}
Let $\cF$ be a locally constant constructible sheaf of $\Lambda$-modules on $U$
and $\ox$ a geometric point of $X$.  We say that $\cF$ is {\em clean} at $\ox$ 
if there exists an étale neighborhood $X'$ of $\ox$ in $X$ such that, 
if we put $U'=U\times_XX'$ and denote by $D'$ the pull-back of $D$ over $X$, 
there exists a finite decomposition 
\begin{equation}
\cF|U'=\oplus_{1\leq i\leq n} \cF'_i
\end{equation}
of $\cF|U'$ into a direct sum of locally constant constructible sheaves of $\Lambda$-modules $\cF'_i$ $(1\leq i\leq n)$ on $U'$
which are isoclinic and clean along $D'$ in the sense of \eqref{clean7}. 
We say that $\cF$ is {\em clean} along $D$ if it is clean at all geometric points of $X$. 
\end{defi}

We will prove in \ref{clean11} that  for isoclinic sheaves, definitions \ref{clean7} and \ref{clean8} are equivalent. 

\begin{lem}\label{clean9}
Let $\cF$ be a locally constant constructible sheaf of $\Lambda$-modules on $U$,
$R$ the conductor of $\cF$ relatively to $X$ and $\ox$ a geometric point of $X$. 
If $\cF$ is clean at $\ox$, then the ramification of $\cF$ at $\ox$ is bounded by $R+$. 
In particular, if $\cF$ is clean along $D$ then the ramification of $\cF$ along $D$ is bounded by $R+$. 
\end{lem}

By assumption, there exists an étale neighborhood $X'$ of $\ox$ in $X$ 
such that, if we put $U'=U\times_XX'$ and denote by $D'$ the pull-back of $D$ over $X$, 
there exists a finite decomposition $\cF|U'=\oplus_{1\leq i\leq n} \cF'_i$ 
of $\cF|U'$ into a direct sum of locally constant constructible sheaves of $\Lambda$-modules $\cF'_i$ $(1\leq i\leq n)$ on $U'$
which are isoclinic and clean along $D'$. For each $1\leq i\leq n$, let $R'_i$ be the conductor of $\cF'_i$ relatively to $X'$. 
The conductor $f^*(R)$ of $\cF|U'$ relatively to $X'$ is the maximum of the $R'_i$ $(1\leq i\leq n)$. 
Therefore, the ramification of $\cF|U'$ along $D'$ is bounded by $f^*(R)$. Then  by \ref{ram7}(ii), the 
ramification of $\cF$ at $\ox$ is bounded by $R+$.

\begin{prop}\label{clean11}
Let $\cF$ be a locally constant constructible sheaf of $\Lambda$-modules on $U$ which is isoclinic along $D$. 
Then $\cF$ is clean along $D$ in the sense of \eqref{clean7} if and only if it is clean along $D$ in the sense of \eqref{clean8}.
\end{prop}

We only need to prove that if $\cF$ is clean along $D$ in the sense of \eqref{clean8}, then it is 
clean along $D$ in the sense of \eqref{clean7}. 
Let $R$ be the conductor of $\cF$ relatively to $X$. We know by \ref{clean9}
that the ramification of $\cF$ along $D$ is bounded by $R+$. 
For every geometric point $\ox$ of $X$, there exists an étale neighborhood $X'$ of $\ox$ in $X$
such that, if we put $U'=U\times_XX'$ and denote by $D'$ the pull-back of $D$ over $X$, 
there exists a finite decomposition 
\begin{equation}
\cF|U'=\oplus_{1\leq i\leq n} \cF'_i
\end{equation}
of $\cF|U'$ into a direct sum of locally constant constructible sheaves of $\Lambda$-modules $\cF'_i$ $(1\leq i\leq n)$ on $U'$
which are isoclinic and clean along $D'$ in the sense of \eqref{clean7}.  Since $\cF$ is isoclinic along $D$, 
$\cF|U'$ is isoclinic along $D'$. 
Hence, for each $1\leq i\leq n$, the conductor of $\cF'_i$ is equal to the pull-back $R'$ of $R$ over $X'$. 
Then it follows from \ref{clean71} that $\cF$ is clean in the sense of \eqref{clean7}.

\end{document}